\documentclass[12pt,leqno]{amsart}

\usepackage{verbatim} 
  
\usepackage{amsmath,amscd,amsthm,amsxtra,amssymb}
\usepackage{epsfig,graphics,color,colortbl}
\usepackage{amssymb,latexsym}
\usepackage{mathrsfs}
\usepackage[poly,all]{xy}
\usepackage{marginnote}
\usepackage{xspace}
\usepackage[colorlinks=true, pdfstartview=FitV, linkcolor=blue,citecolor=blue,urlcolor=blue]{hyperref}

\usepackage{yfonts}
\usepackage{enumerate}

\usepackage[usenames,dvipsnames,svgnames,table]{xcolor}

\usepackage[normalem]{ulem}  

\usepackage[colorinlistoftodos]{todonotes}

\newdir{ >}{{}*!/-10pt/@{>}}

\allowdisplaybreaks[3]

\newlength{\mylength}
\renewcommand{\emptyset}{\varnothing}
\renewcommand{\le}{\leqslant}
\renewcommand{\ge}{\geqslant}

\theoremstyle{plain}
\newtheorem{thm}{\bf Theorem}[section]
\newtheorem{df}[thm]{\bf Definition}
\newtheorem{prop}[thm]{\bf Proposition}
\newtheorem{coro}[thm]{\bf Corollary}
\newtheorem{lem}[thm]{\bf Lemma}
\newtheorem{conj}[thm]{\bf Conjecture}
\newtheorem{sub}[thm]{\bf Sublemma}

\theoremstyle{definition}
\newtheorem{ex}[thm]{\bf Example}
\newtheorem{remark}[thm]{\bf Remark}
\newtheorem{definition}[thm]{\bf Definition}
\newtheorem{question}[thm]{\bf Qusetion}
\newtheorem*{convention}{\bf Convention}
\newtheorem*{prob}{\bf Problem}

\newcommand{\nc}{\newcommand}

\nc{\Prop}{\begin{prop}}
\nc{\enprop}{\end{prop}}
\nc{\Lemma}{\begin{lem}}
\nc{\enlemma}{\end{lem}}
\nc{\Ex}{\begin{ex}}
\nc{\enex}{\end{ex}}
\nc{\Th}{\begin{thm}}
\nc{\enth}{\end{thm}}
\nc{\Def}{\begin{definition}}
\nc{\edf}{\end{definition}}
\nc{\Conj}{\begin{conj}}
\nc{\enconj}{\end{conj}}
\nc{\Quest}{\begin{question}}
\nc{\enquest}{\end{question}}
\nc{\Rem}{\begin{remark}}
\nc{\enrem}{\end{remark}}
\nc{\Ans}{\Proof[{\bf Answer}]}
\nc{\enans}{\QED}
\nc{\Prob}{\begin{prob}}
  \nc{\enprob}{\end{prob}}
\nc{\Sublemma}{\begin{sub}}
\nc{\ensublemma}{\end{sub}}

\newenvironment{red}
{\relax\color{red}}
{\hspace*{.5ex}\relax}

\newcommand{\ber}{\begin{red}}
\newcommand{\er}{\end{red}}

\nc{\berm}{\ber{}\marginnote{\fbox{\scshape\lowercase{M}}}}

\nc{\bern}%
{\ber{}\marginnote{\fbox{\scshape{\lowercase{new}}}}}
\nc{\berMH}{\ber{}\marginnote{\fbox{\scshape\lowercase{MH}}}}
\nc{\berc}{\ber{}\marginnote{\fbox{\scshape\scriptsize{correction}}}}
\nc{\bera}{\ber{}\marginnote{\fbox{\scshape\scriptsize{additional}}}}
\nc{\on}{\operatorname}

\newcommand{\Q}{\mathbb {Q}}
\newcommand{\Z}{\ms{2mu}{\mathbb Z}}

\newcommand{\B}{{\mathbf{B}}}

\newcommand{\D}{\mathscr{D}\ms{1mu}}

\newcommand{\one}{{\bf{1}}}
\newcommand{\seteq}{\mathbin{:=}}

\newcommand{\hd}{{\mathrm{hd}}}      					 
\newcommand{\To}[1][{\hs{0.8ex}}]{\xrightarrow{\ms{7mu}{#1}\ms{7mu}}}

\newcommand{\g}{\ms{1mu}\mathfrak{g}\ms{1mu}}
\newcommand{\n}{\mathfrak{n}}

\newcommand{\Hom}{\operatorname{Hom}}
\newcommand{\HOM}{\on{\mathrm{H{\scriptstyle OM}}}}
\newcommand{\END}{\on{\mathrm{E\scriptstyle ND}}\ms{.1mu}}
\newcommand{\End}{\operatorname{End}}

\newcommand{\isoto}[1][]{\mathop{\xrightarrow%
[{\raisebox{.3ex}[0ex][.3ex]{$\scriptstyle{#1}$}}]%
{{\raisebox{-.6ex}[0ex][-.6ex]{$\mspace{2mu}\sim\mspace{3mu}$}}}}}

\newcommand{\tor}{{\on{tor}}}  

\newcommand{\Mod}{\on{Mod}}
\newcommand{\gmod}{\text{-}\mathrm{gmod}}
\newcommand{\gMod}{\text{-}\mathrm{gMod}}

\newcommand{\F}{\mathscr{F}}

\def\T{{\mathcal T}}

\newcommand{\conv}[1][]{
\underset{\raisebox{.5ex}{$\scriptstyle{#1}$}}{\mathbin{\scalebox{1.1}{$\mspace{1.5mu}\circ\mspace{1.5mu}$}}}}
\newcommand{\hconv}{\mathbin{\scalebox{.9}{$\nabla$}}}

\newcommand{\sconv}{\mathbin{\scalebox{.9}{$\Delta$}}}

\renewcommand{\Im}{\on{Im}}
\newcommand{\de}{\on{\textfrak{d}}}

\newcommand{\cmA}{\cartan}  
\newcommand{\wlP}{\mathsf{P}}   
\newcommand{\rlQ}{\mathsf{Q}}   
\newcommand{\sg}{\mathfrak{S}}   
\newcommand{\Po}{\wlP}

\nc{\qQ}{Q}
\newcommand{\bQ}{\overline{\qQ}}

\newcommand{\Zq}{{\Z[q,q^{-1}]}}  		

\newcommand{\wt}{\mathrm{wt}} 		
\newcommand{\bR}{\mathbf{k}} 		
\nc{\corp}{\bR}
\newcommand{\catC}{ \mathscr{C}}  	
\newcommand{\tcatC}{ \widetilde{\mathscr{C}}}  	
\newcommand{\tcatCs}{ \widetilde{\mathscr{C}}^*}  	
\newcommand{\catT}{ \mathcal{T}}  	
\newcommand{\lT}{ \widetilde{\mathcal{T}}}  	

\newcommand{\catTc}{ \mathcal{T}_{\mathrm{br}}}  	
\newcommand{\dM}{ \mathsf{M }}              
\newcommand{\gW}{\mathsf{W}}
\newcommand{\sgW}{\mathsf{W}^*}
\newcommand{\tE}{\widetilde{E}}  		
\newcommand{\ep}{\varepsilon}  		

\newcommand{\Ht}{\mathrm{ht}} 		
\newcommand{\coRl}{\mathrm{R}^{\mathrm{l}}} 	
\newcommand{\coRr}{\mathrm{R}^{\mathrm{r}}} 	
\newcommand{\La}{\Lambda} 			
\newcommand{\tLa}{\widetilde{\Lambda}} 			
\newcommand{\Res}{\mathrm{Res}} 			

\nc{\Ma}{{\ms{1.5mu}\mathsf{M}}}
\nc{\Ka}{{\ms{1.5mu}\mathsf{K}}}
\nc{\Na}{\mathsf{N}}
\nc{\Xa}{\mathsf{X}}
\nc{\Ya}{\mathsf{Y}}
\nc{\Laa}{\mathsf{L}}

\newcommand{\z}[1][{\Ma}]{{z_{#1}}}

\newcommand{\zN}{{z_\Na}}

\newcommand{\triv}{{\mathbf{1}}}   				
\newcommand{\id}{\ms{2mu}{\mathsf{id}}\ms{1mu}}   				
\newcommand{\dphi}{{\phi}}   				
\newcommand{\gH}{\mathrm{H}}   				

\newcommand{\lG}{\Gamma}   					

\newcommand{\opp}{\mathrm{opp}}

\nc{\be}{\begin{enumerate}}
\newcommand{\bnum}{\be[{\rm(i)}]}
\newcommand{\bna}{\be[{\rm(a)}]}

\newcommand{\rtl}{\rlQ}

\newcommand{\etens}{\boxtimes}

\newcommand{\rmat}[1]{\ms{1mu}{\mathbf{r}}_%
{\mspace{-2mu}\raisebox{-.6ex}{${\scriptstyle{#1}}$}}}

\newcommand{\shc}{{\ms{2mu}\mathcal{C}}}
\newcommand{\shs}{\mathscr{S}}

\newcommand{\Ob}{\on{Ob}}

\nc{\ms}{\mspace}
\nc{\cl}{\colon}
\nc{\ro}{{\rm (}}
\nc{\rf}{{\rm )}\xspace}
\nc{\noi}{\noindent}
\nc{\bl}{\bigl(}
\nc{\br}{\bigr)}

\newenvironment{myequationn}
{\relax\setlength{\arraycolsep}{1pt}\begin{eqnarray*}}
{\end{eqnarray*}}

\newenvironment{myequation}
{\relax\setlength{\arraycolsep}{1pt}\begin{eqnarray}}
{\end{eqnarray}}

\nc{\eq}{\begin{myequation}}
\nc{\eneq}{\end{myequation}}
\nc{\eqn}{\begin{myequationn}}
\nc{\eneqn}{\end{myequationn}}

\newenvironment{myarray}[1]{\relax\setlength{\arraycolsep}{1pt}
\begin{array}{#1}}{\end{array}\relax}

\newcommand{\ba}{\begin{myarray}}
\newcommand{\ea}{\end{myarray}}

\nc{\hs}{\hspace*}
\nc{\vs}{\vspace*}
\nc{\set}[2]{\left\{{#1}\mid{#2}\right\}}
\nc{\snoi}{\smallskip\noi}
\nc{\mnoi}{\medskip\noi}
\nc{\al}{\alpha}
\nc{\rmz}{\setminus\{0\}}
\nc{\tens}
[1][]{\mathbin{\otimes_{\raise1.5ex\hbox to-.1em{}#1}}}
\nc{\vphi}{\varphi}
\nc{\ee}{\end{enumerate}}
\nc{\la}{\lambda}
\nc{\bc}{\begin{cases}}
\nc{\ec}{\end{cases}}
\nc{\qtq}[1][and]{\quad\text{#1}\quad}
\nc{\qt}[1]{\quad\text{#1}}
\nc{\qqt}[1]{\quad\text{#1}}

\nc{\dual}{{\displaystyle{\ms{1mu}\star}}}
\nc{\wle}{\preceq}
\nc{\epito}{\twoheadrightarrow}
\nc{\epiTo}[1][]{\xymatrix@C=4ex{{}\ar@{->>}[r]^-{#1}&{}}}
\nc{\monoTo}[1][]{\xymatrix@C=3ex{\ar@{>->}[r]^-{{#1}}&}}
\nc{\monoto}{\rightarrowtail}
\nc{\Proof}{\begin{proof}}
\nc{\lan}{\langle}
\nc{\ran}{\rangle}
\nc{\ang}[1]{\lan{#1}\ran}
\nc{\QED}{\end{proof}}
\nc{\soplus}{\scalebox{.65}{\raisebox{.2ex}{$\displaystyle\bigoplus$}}}
\nc{\eps}{\varepsilon}
\nc{\supp}{\on{supp}}
\nc{\sct}{strongly commute\xspace}
\nc{\scts}{strongly commutes\xspace}
\nc{\bce}{\eta}			
\nc{\height}{\mathrm{ht}}
\nc{\braid}{{\ms{1mu}\mathrm{br}}}
\nc{\gp}{\mathfrak{p}}
\nc{\wtl}{\wlP}
\nc{\ra}{real and admits an affinization}
\nc{\ras}{real and admit affinizations}
\nc{\Cor}{\begin{coro}}
\nc{\encor}{\end{coro}}
\nc{\shf}{\mathcal{F}}
\nc{\Cw}[1][{w}]{\catC_{{#1}}}
\nc{\tCw}[1][{w}]{\widetilde{\catC}_{{#1}}} 
\nc{\akew}[1][2ex]{\rule[-1ex]{#1}{0ex}}
\nc{\ake}[1][2ex]{\rule[-1ex]{0ex}{#1}}
\nc{\akete}[1][-1ex]{\rule[{#1}]{0ex}{1ex}}
\nc{\tRm}{(R\gmod)\widetilde{\mbox{$\ake[2.5ex]\akew[.9ex]$}}}
\nc{\tX}{\widetilde{X}}
\nc{\corps}{\corp}
\nc{\tL}{\widetilde{L}}
\nc{\prtl}{\rtl_+}
\nc{\nrtl}{\rtl_-}
\nc{\tK}{\widetilde{K}}
\nc{\tep}{\widetilde\ep}
\nc{\teps}{\widetilde\ep}
\nc{\teta}{\widetilde\eta}
\nc{\ga}{\mathfrak{a}}
\nc{\scbul}{{\,\raise1pt\hbox{$\scriptscriptstyle\bullet$}\,}}
\nc{\bwr}{\mbox{\large$\wr$}}
\nc{\tR}{{\widetilde{R}}}
\nc{\lS}{\mathsf{S}}
\nc{\lZ}{\mathcal{Z}}
\nc{\prolim}[1][]{\mathop{\varprojlim}\limits_{{#1}}}
\nc{\sym}{\sg}
\newcounter{myc}
\newcounter{mycc}
\nc{\txi}{\tilde{\xi}}
\nc{\rl}{\rlQ}
\nc{\sfC}{\mathsf{C}}
\nc{\cor}{{\ms{1mu}\mathbf{k}\ms{1mu}}}
\nc{\Pro}{\on{Pro}}
\newcommand{\proolim}[1][]{\ms{-1mu}\mathop{\text{``}\ms{-.5mu}\varprojlim\ms{-4mu}\text{''}}\limits_{#1}}
\nc{\hM}{\widehat{\mathsf{M}}}
\nc{\aff}{\mathrm{aff}}
\nc{\rDa}{{\mathscr{D}_\aff}}
\nc{\st}[1]{\{{#1}\}}
\nc{\W}{\mathsf{W}}
\nc{\rt}{\Delta}
\nc{\pwtl}{\wtl_+}
\nc{\rev}{{\mathrm{rev}}}

\nc{\E}{\mathrm{E}}
\nc{\Es}{\mathrm{E}^*}

\nc{\Qt}{\mathrm{Q}}
\nc{\Qtw}[1][{s_iw_0}]{\mathscr{Q}_{#1}}
\nc{\Qtws}[1][{s_iw_0}]{\mathscr{Q}^*_{#1}}
\nc{\Ctr}{\mathsf{C}}
\nc{\Ctrs}{{\mathsf{C}^*}}
\nc{\Dynkin}{\Delta}
\nc{\cartan}{\mathsf{C}}
\nc{\sfc}{\mathsf{c}}
\nc{\sfa}{\mathsf{a}}
\nc{\SW}{\mathrm{K}}
\nc{\hSW}{\widehat{\mathrm{K}}}
\nc{\refl}{\mathscr{S}}
\nc{\Rre}{\mathrm{R}^{\mathrm{ren}}}
\nc{\Rpre}{\mathrm{R}^{\mathrm{ren}\;'}}
\nc{\bRre}{\ol{\mathrm{R}}^{\ms{2mu}\mathrm{ren}}}
\nc{\sfd}{\ms{1mu}\mathsf{d}\ms{1mu}}
\nc{\shm}{\mathcal{M}}
\nc{\sht}{\mathcal{T}}
\nc{\rank}{\mathrm{rank}}
\nc{\Da}{{\D}\ms{-2.8mu}\raisebox{-.35ex}{$\scriptstyle\mathrm{aff}$}}
\nc{\lDa}{\Da^{-1}}
\nc{\bchi}{{\scalebox{.9}{\mbox{$\mathscr{E}$}}}}
\nc{\bchis}{\bchi{}^{\ms{2mu}*}}
\nc{\Laf}{\mathscr{L}}
\nc{\tLaf}{\widetilde{\Laf}}
\nc{\wtaf}{\mathscr{W}\ms{-3mu}{\mathit{t}}}
\nc{\res}[1][]{\mathop\star\limits_{\raisebox{.4ex}{$\scriptstyle #1$}}\ms{2mu}}
\nc{\hchi}{\widehat{\chi}}
\nc{\convaff}{\mathop{\scalebox{1.1}{$\mspace{1.5mu}\circ\mspace{1.5mu}$}}\limits}
\nc{\cvb}{CVB\xspace}
\nc{\svelt}{essentailly samll\xspace}
\nc{\Proc}{\on{Pro}_{\mathrm{coh}}}
\nc{\sha}{\mathcal{A}}
\nc{\Ker}{\on{Ker}}
\nc{\Coker}{\on{Coker}}
\nc{\Aff}[1][z]{\on{Aff}_{\ms{1mu}#1}}
\nc{\scb}{\scalebox}
\nc{\afr}{affreal\xspace}
\nc{\epifrom}{\ms{-5mu}\xymatrix@C=3ex{{}&{}\ar@{->>}[l]}\ms{-5mu}}
\nc{\Mid}{\bigm|}
\nc{\ol}{\overline}
\nc{\bpsi}{\ol{\psi}}
\nc{\Rat}[1][z]{\on{Raff}_{\ms{1mu}#1}}
\nc{\indlim}{\varinjlim\limits}
\nc{\inddlim}{\mathop{\mbox{``{$\varinjlim$}''}}\limits}
\nc{\Rmat}{\ms{1mu}\mathrm{R}\ms{1mu}}
\nc{\Runi}{\mathrm{R}^{\mathrm{uni}}}
\nc{\Modg}{\mathrm{Modg}}
\nc{\KO}{quasi-rigid\xspace}
\nc{\hF}{\widehat{\F}}
\nc{\Modc}{\Mod_{\mathrm{coh}}}
\nc{\e}{\mathrm{e}}
\nc{\Idx}{\mathsf{\Lambda}}
\nc{\hA}{\widehat{A}}
\nc{\prood}{\mathop{\text{``}\prod\text{''}}\limits}

\nc{\hrefl}{\widehat{\mathscr{S}}}
\nc{\ev}{\mathrm{ev}}
\nc{\coev}{\mathrm{coev}}
\nc{\ihom}{\mathcal{H}om}
\nc{\tY}{\widetilde{Y}}
\nc{\tensz}{\tens[z]\ms{-3.5mu}}

\nc{\tRre}{\widetilde{\mathrm{R}}^{\mathrm{ren}}}
\nc{\dg}{\mathbf{\lambda}}
\nc{\can}{\mathrm{can}}
\nc{\dc}{\mathfrak{c}}
\nc{\htens}{\hconv}
\nc{\stens}{\sconv}
\nc{\Modgc}{\mathrm{Modg}_{\mathrm{coh}}}
\nc{\Aut}{\mathrm{Aut}}
\nc{\Rd}[1][\dg]{R_{#1}\gmod}
\nc{\nn}{\nonumber}
\nc{\Dual}{\mathrm{D}\ms{1mu}}
\nc{\DA}[1][A]{\ms{1mu}\mathrm{D}_{{#1}}}
\nc{\DmA}[1][A]{\ms{1mu}\Dual^{-1}_{{#1}}}

\nc{\DB}[1][B]{\ms{1mu}\D_{{#1}}}
\nc{\DmB}[1][B]{\ms{1mu}\D^{-1}_{{#1}}}

\nc{\tensa}{\tens[A]\ms{-3mu}}
\nc{\tensc}{\tens[{\ms{3mu}\cor}]\ms{-3mu}}
\nc{\bg}{{\ms{2mu}\mathrm{big}}}
\nc{\afn}{affine object\xspace}
\nc{\afns}{affine objects\xspace}
\nc{\subafn}{affine subobject\xspace}
\nc{\Afns}{Affine objects\xspace}
\nc{\dL}{\mathsf{L}}
\nc{\dX}{\mathsf{X}}
\nc{\dY}{\mathsf{Y}}
\renewcommand{\preceq}{\preccurlyeq}
\nc{\ltens}{\otimes\limits^{\mathbb{L}}}
\nc{\epl}{\epsilon}
\nc{\tI}{\widetilde{I}}
\nc{\tfC}{\widetilde{\sfC}}
\nc{\tfc}{\widetilde{\sfc}}
\nc{\Loc}{\mathscr{A}}
\nc{\Ri}{(R\gmod)_i}
\nc{\iR}{{}_i(R\gmod)}
\nc{\into}{\xymatrix@C=2.5ex{{}\ar@{^{(}->}[r]&{}}}
\nc{\tC}{\widetilde{C}}
\nc{\De}{\mathfrak{D}}
\nc{\Daf}{{\D}\ms{-2.8mu}\raisebox{-.35ex}{$\scriptstyle\mathrm{aff}$}}
\nc{\Qtp}{\Qt_+}
\nc{\Qtm}{\Qt_-}
\nc{\Di}{\D(\Qtp\ang{i})}
\nc{\Dip}{\D^{-1}(\Qtp\ang{i_+})}
\nc{\Diz}[1][z]{\Daf(\Qtp\ang{i}_{#1})}
\nc{\sDipz}{\D_{\scriptscriptstyle\mathrm{aff}}^{-1}(\Qtp\ang{i_+}_{z})}
\nc{\sDiz}[1][z]{\D_{\scriptscriptstyle\mathrm{aff}}(\Qtp\ang{i}_{#1})}
\nc{\Dipz}{\Daf^{-1}(\Qtp\ang{i_+}_{z})}
\nc{\Qipz}{\,\Qtp(\ang{i_+}_{z})}
\nc{\Dizi}{\Daf(\Qtp\ang{i}_{z_i})
}\nc{\iz}[1][z]{\ang{i}_{#1}}
\nc{\qiz}[1][z]{\Qtp(\ang{i}_{#1})}
\nc{\qi}{\Qtp(\ang{i})}
\nc{\tdM}{\widetilde{\mathsf{M}}}
\nc{\Locp}{\Loc_+}
\nc{\wh}[1]{\widehat{#1}}
\nc{\uprod}{\prod\limits^{\xleftarrow{}}}
\nc{\tcd}[1]{\textcircled{\scriptsize{#1}}}
\nc{\convz}{\conv[z]}
\nc{\Qi}[1][i]{\Qtp(\ang{#1})}
\nc{\Qiz}[1][z]{\Qtp(\ang{i}_{#1})}
\nc{\Ld}{\mathsf{L}}
\nc{\Psip}{\Psi_{+,-}}
\nc{\hcalA}{\widehat{\mathcal{A}}}
\nc{\TT}{\textbf{\textit{T}}}
\nc{\Qq}{\Q(q)}
\nc{\hK}{\widehat{K}}
\nc{\udset}[2]{\underset{\raisebox{2ex}{$#1$}}{#2}}
\nc{\rs}{r^*}
\nc{\ta}[2]{\tau_{#1}(#2)}
\nc{\rd}{\mathrm{rd}}
\nc{\ld}{\mathrm{ld}}
\nc{\Rc}{{R^{\mathrm{cyc}}}}
\nc{\Rh}{\mathrm{R}}
\nc{\At}{A\tens_\cor}
\nc{\Loco}[1][\pm]{(R_{#1}\gmod)^0}

\numberwithin{equation}{section}
\begin{document}

\title {Reflection functors on quiver Hecke algebras}

\author[M. Kashiwara]{Masaki Kashiwara}
\thanks{The research of M.\ Kashiwara
	was supported by Grant-in-Aid for Scientific Research (B)  23K20206,  
	Japan Society for the Promotion of Science.}
\address[M. Kashiwara]{%
Kyoto University Institute for Advanced Study, Research Institute
for Mathematical Sciences, Kyoto University, Kyoto 606-8502, Japan
}
\email[M. Kashiwara]{masaki@kurims.kyoto-u.ac.jp}

\author[M. Kim]{Myungho Kim}
\address[M. Kim]{Department of Mathematics, Kyung Hee University, Seoul 02447, Korea}
\email[M. Kim]{mkim@khu.ac.kr}
\thanks{The research of M.\ Kim was supported by the National Research Foundation of
Korea (NRF) Grant funded by the Korea government(MSIT)
(
NRF-2020R1A5A1016126).}

\author[S.-j. Oh]{Se-jin Oh}
\thanks{ The research of S.-j.\ Oh was supported by the National Research Foundation of
	Korea (NRF) Grant funded by the Korea government(MSIT) (NRF-2022R1A2C1004045).}
\address[S.-j. Oh]{ Department of Mathematics, Sungkyunkwan University, Suwon, South Korea}
\email[S.-j. Oh]{sejin092@gmail.com}

\author[E. Park]{Euiyong Park}
\thanks{The research of E.\ Park was supported by the National Research Foundation of Korea (NRF) Grant funded by the Korea Government(MSIT)(RS-2023-00273425 and NRF-2020R1A5A1016126).}
\address[E. Park]{Department of Mathematics, University of Seoul, Seoul 02504, Korea}
\email[E. Park]{epark@uos.ac.kr}

\makeatletter
\@namedef{subjclassname@2020}{\textup{2020} Mathematics Subject Classification}
\makeatother


\date{November 8, 2025}

\begin{abstract}
  In this paper, we construct the reflection functors
  for quiver Hecke algebras of an arbitrary 
 symmetrizable Kac-Moody type. These reflection functors categorify Lusztig's braid symmetries. 	
\end{abstract}

\maketitle

\tableofcontents

\section{Introduction}

The categorification using \emph{quiver Hecke algebras}  is a powerful tool to study quantum groups and their representations. 
Let $ \cmA$ be  a symmetrizable generalized  Cartan matrix and  $U_q(\g)$ the  corresponding quantum group. 
It was shown in \cite{KL09, KL11,R11} that the \emph{quantum unipotent coordinate ring} $A_q(\mathfrak n)$, which is the graded dual of  the negative part $U_q^-(\g)$,  is categorified by  the monoidal category $R\gmod$ of finite-dimensional graded modules over the quiver Hecke algebra $R$ associated with $\cmA$. Several interesting and important features of $A_q(\mathfrak n)$ have been lifted to the category $R\gmod$ and studied via the categorification. 
This categorical approach produced various successful results on quantum groups (see \cite{BK09, KK11, K^3, KKKO18, LV11, R11, VV09} and see also the references therein).

Lusztig's braid symmetries ${\mathsf S}_i$ (\cite{Lu93})  also have been interpreted and studied by the quiver Hecke algebras.  Here, ${\mathsf S}_i$ denotes the restriction of the automorphism of $U_q(\g)$ described in \cite{Lu93} to certain subalgebras ${\mathsf S}_i\cl A_q(\mathfrak{n})[i] \isoto A_q(\mathfrak{n})[i]^*$ (see, for example, \cite{ext}).

When $\cmA$ is of finite $ADE$ type, Kato constructed monoidal functors, called \emph{reflection functors}, at the level of the category $R\gmod$ and proved that these functors categorify the braid symmetries in a geometric manner (\cite{Kato14}). This geometric approach was extended to the cases of arbitrary symmetric type (\cite{Kato20}). On the other hand, an algebraic approach to construct the reflection functors was introduced by the authors using \emph{Schur-Weyl duality} of quiver Hecke algebras (\cite{refl}).
We assume that $\cmA$ is of \emph{arbitrary finite} type and denote by $\tcatC_{s_iw_0}$ (resp.\ $\tcatCs_{s_iw_0}$) the localization of the subcategory $\catC_{s_iw_0}$ (resp.\ $\catC^*_{s_iw_0}$) of $R\gmod$ (see \cite[Section 10]{refl} for details).
In this approach, the Schur-Weyl  duality datum arises from the \emph{affinizations} of certain determinantial modules and the right dual of the 1-dimensional module $\ang{i}$ in the localized categories. 
This duality datum provides the Schur-Weyl duality functor 
$$
\shs_i\cl \tcatC_{s_iw_0}\to\tcatCs_{s_iw_0},
$$
which gives the braid symmetry at the Grothendieck ring level. This functor is also called the \emph{reflection functor}.
Since the finite type condition is crucially used in the course of the proof, this approach cannot be applicable to the case of an
arbitrary symmetrizable Kac-Moody type.

In this paper, we construct reflection functors ${\F_i}$
for an \emph{arbitrary symmetrizable Kac-Moody type},  which categorify the braid symmetries. 
We thus complete the construction of reflection functors using the quiver Hecke algebras for all types. 
We develop a new algebraic approach using Schur-Weyl duality which produces desired functors. 
Unlike the previous approach in \cite{refl}, we no longer need the finite type  assumption.  
The key ingredient is to enlarge  the quiver Hecke algebra $R$ to larger algebras $R_\pm$ by adding a new index $i_\pm$ which is adjacent to $i$.
We  construct the localizations $\Loc_\pm$  of the enlarged categories $ R_\pm\gmod $ by special simple modules  $C_\pm$ (see \eqref{Eq: Cpm}). 
In the localization $\Loc_+$ (resp.\ $\Loc_-$), the simple object $\ang{i_+}$ (resp.\ $\ang{i_-}$) can be understood as the right (resp.\ left) dual of $\ang{i}$ up to a multiple of $C_+$ (resp.\ $C_-$).
These enlarged localizations allow us to construct exact monoidal functors $\F_i\cl  \Loc_- \To \Loc_+ $ and $\F^*_i\cl  \Loc_+ \To \Loc_- $ using the Schur-Weyl duality. 
We then restrict the functors $\F_i$ and $\F^*_i$ to the full subcategories $(R\gmod){}_i$ and ${}_i(R\gmod)$ to obtain the reflection functors
$$
\F_i\cl (R\gmod){}_i \longrightarrow {}_i(R\gmod) \quad \text{ and } \quad  \F^*_i\cl {}_i(R\gmod) \longrightarrow (R\gmod){}_i
$$
(see \eqref{Eq: RiiR} for the definition of $(R\gmod){}_i$ and ${}_i(R\gmod)$).
It should be noted that such a restriction was not successfully established in \cite{refl}.  
The functors $\F_i$ and $\F^*_i$ are  quasi-inverses  to each other
and they coincide with the braid symmetries at the Grothendieck ring level.

\medskip

Let us explain our results in more detail.
Let $R$ be the quiver Hecke algebra associated with an arbitrary generalized Cartan matrix $\cmA$. 
We add a new index $i_\pm$ to the index set $I$ to make  a bigger one
$I_\pm \seteq I \sqcup\{ i_\pm \}$ and define a new Cartan matrix $ \tfC_\pm$
as in \eqref{Eq: tC}, respectively.
The new added index $i_\pm$ is only adjacent to $i$ in the Dynkin diagram.
The new quiver Hecke algebra $R_\pm$ is defined by the polynomial parameters $ Q_{j,k}^\pm(u,v)$ (see \eqref{Eq: Qpm}) and the grade associator $\lambda_\pm$ (see \eqref{eq:new grading}).
We then consider the localization $\Loc_\pm$ of $R_\pm\gmod$ by the simple modules 
\begin{align*}
C_+\seteq\ang{i\,i_+}\in R_+\gmod \qtq
C_-\seteq\ang{i_-\;i}\in R_-\gmod
\end{align*}
with the left and right braider structure,  respectively (see Proposition \ref{prop:Cpm}). By the choice of the grade associator $\lambda_\pm$, the braidings  $\coRl_{C_+}(M)\cl C_+\conv M\to M\conv C_+$
and $\coRr_{C_-}(M) \cl M\conv C_-\to C_-\conv M$ ($M\in R_\pm\gmod$)  are
homogeneous of degree 0. 
Under the localization functors $ \Qt_\pm\cl R_\pm\gmod\to \Loc_\pm$, 
we show that $R\gmod$ is
embedded into $\Loc_\pm$ as a full subcategory
and stable by taking subquotients
(Propositions \ref{prop:RLocf} and \ref{prop:stab subq}), and investigate several compatibilities. 
Note that the duals $ \D^{\pm 1} (\Qt_\pm(\ang{i}))$
is isomorphic to $\Qt_\pm(\ang{i_\pm})$ 
up to a multiple of $\tC_\pm \seteq Q_\pm(C_\pm)$ in $\Loc_\pm$.

 As a Schur-Weyl duality datum for $i\in I$,  we choose real simple objects $\{ K_j^\pm \}_{j\in I_\mp} \subset \Loc_\pm$ defined in \eqref{Eq: Kj+} and \eqref{Eq: Kj-}
and their affinizations $\{ \tK_j^\pm \}_{j\in I_\mp}$ defined in \eqref{Eq: aKj+} and \eqref{Eq: aKj-}.
Note that $\tK_i^{\pm}$ and $\tK_{i_\mp}^{\pm}$ are defined as the duals of $\Qt_{\pm}(\ang{i}_z)$ and $\Qt_{\pm}(\ang{i_{\pm}}_z)$ in the category of affine objects and they are equipped with rational center structures (Proposition \ref{prop:affD}).
Proposition \ref{prop:LaSW} and Theorem \ref{thm:DeSW} say that 
$$
\La(K^\pm_j,K^\pm_k)=\la_{\mp}(\al_j,\al_k) \qtq \De(\tK^\pm_j,\tK^\pm_k)=Q^{\mp}_{j,k}(z_j,z_k)
$$
for $j,k\in I_\mp$ with $j\ne k$. Thus the family $\{ \tK_j^\pm \}_{j\in I_\mp}$  forms a Schur-Weyl duality datum (see Section \ref{Sec: sub SW duality}), which gives a duality functor 
$$
\Psi_{\pm,\mp}\cl  R_\mp  \gmod\to \Loc_\pm. 
$$
Lemma \ref{lem:Psi C} says that $\Psi_{\pm,\mp}(C_\mp)\simeq(\tC_\pm)^{\circ-1}$, which implies that  $\Psi_{\pm,\mp}$ factors through $\Loc_\mp$ by the universal property of localization. We  denote them  by $\F_i$ and $\F^*_i$, i.e., 
$$
\xymatrix{
R_- \gmod \ar[d]_{Q_-} \ar[drr]^{\Psi_{+,-}}	&&  && R_+\gmod \ar[d]_{Q_+}   \ar[drr]^{\Psi_{-,+}} && \\
\Loc_- \ar[rr]_{\F_i} && \Loc_+, && \Loc_+  \ar[rr]_{\F^*_i}  &&\Loc_-.
}
$$
The functors $ \F_i$ and $\F^*_i$ are quasi-inverses to each other (see Proposition \ref{prop:FF*}), which implies that $\Psi_{\pm,\mp}$ is  exact.

For any $\nu \in I^n$, we define an $R$-module $ \dM(\nu)\seteq R u(\nu)$ by the defining relations \eqref{Eq: Mnu}. 
We show that $\dM(\nu)$ is a projective object in $(R\gMod){}_i$ and  any objects in  $(R\gmod){}_i$ can be obtained  as a quotient of a direct sum of $\dM(\nu)$'s.
Since $\Psi_{+,-}(\dM(\nu)) \in {}_i(R\gMod) $ as seen in the proof of
Proposition \ref{prop:FiQtmRi}, we have  the restricted functor 
$$
\F_i\cl (R\gmod){}_i \longrightarrow {}_i(R\gmod).
$$ 
In a similar manner, one can obtain its quasi-inverse 
$$ 
 \F^*_i\cl {}_i(R\gmod) \longrightarrow (R\gmod){}_i.
$$
We finally show that the  homomorphism $K( (R\gmod){}_i) \longrightarrow  K( {}_i(R\gmod))$ induced from $\F_i$ coincides with the braid symmetry (Theorem \ref{th:saito}).

In the course of arguments,  Theorem \ref{th:J} takes an important technical part about homomorphisms involving $E_i$. The proof of Theorem  \ref{th:J} is given in Appendix \ref{sec:appendixB}.

This paper is organized as follows. 
In Section\;\ref{Sec: Preliminaries}, we review several results on
quiver Hecke algebras including localization and Schur-Weyl duality.  
In Section\;\ref{Sec: Apm}, we construct the localizations $\Loc_\pm$ of the enlarged algebras $R_\pm$, and Section\;\ref{Sec: SW duality} explains the Schur-Weyl duality functors in the setting of the localization $\Loc_\pm$.
In Section\;\ref{Sec: gen for iRgmod} and Section\;\ref{Sec: reflection}, we investigate the generators of the category ${}_i(R\gmod)$ and construct the reflection functors categorifying braid symmetries.
Appendix \ref{sec:appendixA},  \ref{sec:appendixB},  and \ref{sec:appendixC}  
are devoted to the proof of Proposition \ref{prop:affD},  the proof of Theorem \ref{th:J} and some complements to \cite{refl}  respectively.

\mnoi
{\bf Acknowledgments}
%
\quad The authors would like to thank the anonymous referee for numerous precious comments and even several simplifications of the proofs, which have greatly helped to improve the paper.

\begin{convention}  Throughout this paper, we use the following convention.
\bnum
\item For a statement $P$, we set $\delta(P)$ to be $1$ or $0$ depending on whether $P$ is true or not. In particular, $\delta_{i,j}=\delta(i=j)$ is the Kronecker delta.
  \item For a ring $A$, $A^\times$ denotes the group consisting of invertible elements of $A$.
\item For a monoidal category $\catT$ with the tensor product $\tens$,
  we denote by $\catT^\rev$ the monoidal category with the reversed tensor product $\tens_\rev$:
  $ X\tens_\rev Y\seteq Y\tens X$. \label{conv3}
\item  For elements $A,B$ in a vector space over a field $\cor$,  $A \equiv B$  means that $A =c B$ for some $c\in \cor^\times$. 
  \label{conv:4} 

\ee
\end{convention}

\section{Preliminaries} \label{Sec: Preliminaries}

\subsection{Quiver Hecke algebras and universal R-matrices}\label{subsec:QHA}
In this subsection, we review the notion of quiver Hecke algebras and  some  related topics  following \cite{KL09, Rouquier08, K^3}. 
\subsubsection{Quiver Hecke algebras}
A {\em Cartan datum} $ \bl\cmA,\wlP,\Pi,\Pi^\vee,(\cdot,\cdot) \br $  is a quintuple consisting of a generalized Cartan matrix $\cartan= (\sfc_{i,j})_{i,j\in I}$, a free abelian group $\wlP$ called the weight lattice, a set of simple roots  $\Pi = \{ \alpha_i \mid i\in I \} \subset \wlP$, a set of simple coroots $\Pi^{\vee} = \{ h_i \mid i\in I \} \subset \wlP^{\vee}\seteq\Hom( \wlP, \Z )$, and   
a $\Q$-valued  symmetric bilinear  $(\cdot,\cdot)$ form on $\wlP$
  satisfying the following conditions:
\begin{enumerate} [{\rm (a)}]
\item $\cmA = (\langle h_i,\alpha_j\rangle)_{i,j\in I}$,
\item  $(\alpha_i,\alpha_i)\in 2\Z_{>0}$ for $i\in I$,
\item $\langle h_i, \lambda \rangle =\dfrac{2(\alpha_i,\lambda)}{(\alpha_i,\alpha_i)}$ for $i\in I$ and $\lambda \in \Po$,
\item for each $i\in I$, there exists $\Lambda_i \in \wlP$
such that $\langle h_j, \Lambda_i \rangle = \delta_{ij}$ for any $j\in I$.
\end{enumerate}
We set $\sfd_i=\dfrac{(\al_i,\al_i)}{2} \in \Z_{>0}$ for $i\in I$.

The corresponding symmetrizable Kac-Moody algebra will be denoted by $\g\seteq\g(\cmA)$.
Let $\pwtl\seteq\st{\la\in\wlP\mid 
\text{$\ang{h_i,\la}\ge0$ for any $i\in I$}}$ be 
the set of dominant integral weights,
 $\rlQ\seteq\soplus_{i\in I} \Z\al_i$ the root lattice, $\rlQ_+\seteq\soplus_{i\in I} \Z_{\ge0} \al_i$, 
the positive root lattice,
and $\rlQ_-\seteq-\rlQ_+$ the negative root lattice. 
For $\beta=\sum_{i\in I} b_i \al_i\in\rtl$, we set $\height(\beta)\seteq \sum_{i\in I} |b_i|$.

\medskip
Let $\bR$ be a field.
A {\em quiver Hecke datum} associated with a Cartan matrix $\cmA$ is 
a family of polynomials
$\qQ_{i,j}(u,v) \in \bR[u,v]$ of the form
\begin{align}
\qQ_{i,j}(u,v) =\bc
                   \sum\limits
_{p(\alpha_i , \alpha_i) + q(\alpha_j , \alpha_j) = -2(\alpha_i , \alpha_j) } t_{i,j;p,q} u^pv^q &
\text{if $i \ne j$,}\\[3ex]
0 & \text{if $i=j$,}
\ec\label{eq:Q}
\end{align}
such that $t_{i,j;-a_{ij},0} \in  \bR^{\times}$ and
$$\qQ_{i,j}(u,v)= \qQ_{j,i}(v,u) \quad \text{for all} \ i,j\in I.$$

For $\beta\in \rlQ_+$  with $\height(\beta)=n$,  set
$$
I^\beta\seteq  \Bigl\{\nu=(\nu_1, \ldots, \nu_n ) \in I^n \bigm| \sum_{k=1}^n\alpha_{\nu_k} = \beta \Bigr\}.
$$
Then the symmetric group $\mathfrak{S}_n = \langle s_k \mid k=1, \ldots, n-1 \rangle$ acts  by place permutations on $I^\beta$.

\begin{df}
\ For $\beta\in\rlQ_+$ with $\height(\beta)=n$,
the {\em quiver Hecke algebra} $R(\beta)$ associated with a Cartan datum $ \bl\cmA,\Pi,\wlP,\Pi^\vee,(\cdot,\cdot) \br $ and a quiver Hecke datum $\st{\qQ_{i,j}(u,v)}_{i,j\in I}$
is the $\bR$-algebra generated by
$$
\{e(\nu) \mid \nu \in I^\beta \}, \; \{x_k \mid 1 \le k \le n \},
 \; \{\tau_l \mid 1 \le l \le n-1 \}
$$
satisfying the following defining relations:
\eqn
&& e(\nu) e(\nu') = \delta_{\nu,\nu'} e(\nu),\ \sum_{\nu \in I^{\beta}} e(\nu)=1,\
x_k e(\nu) =  e(\nu) x_k, \  x_k x_{k'} = x_{k'}x_k,\\
&& \tau_l e(\nu) = e(s_l(\nu)) \tau_l,\  \tau_k \tau_l = \tau_l \tau_k \text{ if } |k - l| > 1, \\[1ex]
&&  \tau_l^2 = \sum_{\nu\in I^\beta}\qQ_{\nu_l, \nu_{l+1}}(x_l, x_{l+1})e(\nu), \\[5pt]
&& \tau_l x_k - x_{s_l(k)} \tau_l =
\bl\delta(k=l+1)-\delta(k=l)\br
\sum_{\nu\in I^\beta,\ \nu_l=\nu_{l+1}}e(\nu),\\
&&\tau_{l+1} \tau_{l} \tau_{l+1} - \tau_{l} \tau_{l+1} \tau_{l}
=\sum_{\nu\in I^\beta,\ \nu_l=\nu_{l+2}}
\bQ_{\,\nu_l,\nu_{l+1}}(x_l,x_{l+1},x_{l+2}) e(\nu),
\eneqn
\end{df}
where
\begin{align*}
\bQ_{i,j}(u,v,w)\seteq\dfrac{ \qQ_{i,j}(u,v)- \qQ_{i,j}(w,v)}{u-w}\in \bR[u,v,w].
\end{align*}
For an $R(\beta)$-module $M$, we set $\wt(M)\seteq-\beta \in \rtl_-$. 

For $\mu=(\mu_1,\ldots,\mu_m) \in I^{\beta}$, $\nu = (\nu_1,\ldots,\nu_n)\in I^{\gamma}$,  we denote the concatenation $e(\mu_1,\ldots,\mu_m,\nu_1,\ldots,\nu_n)$  by $e(\mu,\nu)$.
For $\beta,\gamma \in \rtl_+$, we set $e(\beta,\gamma)\seteq \sum_{\mu \in I^{\beta}, \, \nu\in I^{\gamma}} e(\mu ,\nu)\in R(\beta+\gamma)$.
Then, we have an algebra embedding
$$R(\beta)\tens R(\gamma)\subset e(\beta,\gamma) R(\beta+\gamma)e(\beta,\gamma).$$

For an $R(\beta)$-module $M$ and an $R(\gamma)$-module $N$, their convolution product is defined by 
\eqn
M\conv N\seteq R(\beta+\gamma)e(\beta,\gamma) \otimes_{R(\beta)\tens R(\gamma)} (M\tens N).
\eneqn

\vskip 2em

\subsubsection{Restriction functors}
Let $i\in I$.
 For $\beta\in\prtl$, we set
$$e(i,*)=e(\al_i,\beta-\al_i)\in R(\beta)\qtq e(*,i)=e(\beta-\al_i,\al_i)\in R(\beta).$$
Then for an $R(\beta)$-module $M$,
$$\E_iM=e(i,*)M\qtq \E^*_iM=e(*,i)M$$
have an $R(\beta-\al_i)$-module  structure.
We have  morphisms
$$R(\al_i)\conv\E_iM\to M\qtq \E^*_iM\conv R(\al_i)\to M.$$
(See Theorem~\ref{th:J} for their dual form.) 

\subsubsection{Universal R-matrices} 
Let $\beta \in \rlQ_+$ with $m =  \Ht(\beta)$. For  $k=1, \ldots, m-1$,
the \emph{intertwiner} $\varphi_k \in R(\beta) $ is defined by 
\eq
\varphi_k e(\nu) =
\bc
 \bl\tau_k(x_k-x_{k+1})+1\br e(\nu) 
& \text{ if } \nu_k = \nu_{k+1}, 
 \\
 \tau_k e(\nu) & \text{ otherwise}
\ec\quad\qt{for $\nu \in I^\beta$.} \label{def:intertwiner}
\eneq
Since $\st{\vphi_k\mid  1\le k \le m-1}$ satisfies the braid relation,
$\vphi_w\in R(\beta)$ is well-defined for any  $w\in\sg_m$.

For $m,n \in \Z_{\ge 0}$, we set $w[m,n]$ to be the element of $\sg_{m+n}$ such that
$$
w[m,n](k) \seteq
\left\{
\begin{array}{ll}
 k+n & \text{ if } 1 \le k \le m,  \\
 k-m & \text{ if } m < k \le m+n.
\end{array}
\right.
$$

For an $R(\beta)$-module $M$ and an $R(\gamma)$-module $N$,
the map $M \otimes N \rightarrow e(\beta,\gamma)(N \conv M)$ defined by $$u \otimes v \mapsto \varphi_{w[\height{\gamma},\height{\beta}]}(v \etens u)$$
is $R(\beta)\otimes R(\gamma)$-linear and it
extends to an $R(\beta+\gamma)$-module homomorphism 
$$
\Runi_{M,N}\cl  M\conv N \longrightarrow N \conv M.
$$
We call $\Runi_{M,N}$ the {\em universal $R$-matrix} between $M$ and $N$.
Since the intertwiners satisfy the braid relations,  the universal R-matrices $\Runi_{M,N}$ satisfy the Yang-Baxter equation (\cite[(1.9)]{K^3}).

\subsection{Gradings on $R(\beta)$ and degrees of R-matrices} 
\hfill
In this subsection we will review the gradings on  $R(\beta)$ and degrees of R-matrices following mainly \cite{KKKO18, KP18}.
\subsubsection{Gradings on quiver Hecke algebras}
 A bilinear form $\dg\cl \rtl\times\rtl\to\Z$ is called
a {\em grade associator} if it satisfies 
\eq \label{eq: dg_condition}
\dg(\al,\beta)+\dg(\beta,\al)=-2(\al,\beta)\qt{for any $\al$, $\beta\in\rtl$. }
\eneq
Then, the algebra $R(\beta)$ is $\Z$-graded by
\eq \label{eq: dg_formula}
&&\deg(e(\nu))=0, \quad \deg(x_k e(\nu))= ( \alpha_{\nu_k} ,\alpha_{\nu_k}), \quad  \deg(\tau_l e(\nu))= \dg(\alpha_{\nu_{l+1}} , \alpha_{\nu_{l}}).
\label{eq:grading}
\eneq

\begin{remark}
Note that  \eqref{eq: dg_formula} is
  different from the definition of $\deg\bl \tau_l e(\nu)\br$ in \cite[(7.2)]{refl}, which was erroneous, but it did not affect the rest of  \cite{refl}. 
\end{remark}

We denote by $R_\dg(\beta)$ the graded algebra $R(\beta)$ with the grading 
\eqref{eq:grading}. 
The standard grading on $R(\beta)$ given by
$\dg_\can\seteq-(\;\cdot\;,\;\cdot\;)$ is introduced in \cite{KL09, Rouquier08}.
We will denote the graded algebra $R_{\dg_\can}(\beta)$ simply  by $R(\beta)$.

We denote by $\Modg(R_\dg(\beta))$ the $ \cor$-linear abelian category of
graded $R_{\dg}(\beta)$-modules whose morphisms are the  grading preserving $R_\dg(\beta)$-module homomorphisms.

The grading shift functor is  given by $(q M)_k=M_{k-1}$ $(k \in \Z)$.
Set
\eqn
\HOM_{\Modg(R_\dg(\beta))}(M,N) \seteq \bigoplus_{d\in\Z} \Hom_{\Modg(R_\dg(\beta))}(q^dM,N).
\eneqn
Note that if $M$ is finitely generated, then  $\HOM_{\Modg(R_\dg(\beta))}(M,N)$ is the set of $R$-linear homomorphisms from $M$ to $N$ not necessarily preserving the grading.  
We set $\Modg(R_\dg)=\soplus_{\beta\in\prtl}\Modg(R_\dg(\beta))$. Then $\Modg(R_\dg)$ is a graded monoidal category (see \eqref{eq:graded_monoidal}) 
whose tensor product is  the convolution product $\conv$.

For each $i\in I$,
there is a functor 
$\E_i \cl \Modg(R_\dg) \to \Modg(R_\dg)$  sending an $R_\dg(\beta)$-module $M$  to the $R_\dg(\beta-\al_i)$-module $e(\al_i,\beta-\al_i)M$. 
If $M\in \Modg(R_\dg)$ is simple, then  we denote by $\tE_i(M)$ the head of $\E_i(M) $, which is simple.

We denote by $R_\dg(\beta)\gmod$  the  full subcategory of $\Modg(R_\dg(\beta))$ consisting of  finite-dimensional modules over $\cor$. 
Set $R_\dg\gmod \seteq \soplus_{\beta\in\prtl} R_\dg(\beta)\gmod$.
We denote by $q$ the grading shift functor: $(qM)_n=M_{n-1}$ for
$M\in R_\dg\gmod$.
For $i\in I$, we write $q_i=q^{(\al_i,\al_i)/2}$. 

Let us denote by $\Modg\bl R_\dg(\beta)\br[q^{1/2}]$ the category of $ (\frac{1}{2} \Z)$-graded modules over $R_\dg(\beta)$.
Define a skew-symmetric bilinear form $\dc$ on $ \rtl$ as
$$\dc (\beta,\gamma)\seteq \dg(\beta ,\gamma)+(\beta,\gamma) \quad \text{for} \ \beta,\gamma \in \rtl.$$
For $\beta \in \rtl_+$ and $M \in \Modg(R(\beta))[q^{1/2}]$, set
$$(K_\dc(M))_n \seteq \bigoplus_{\nu\in I^\beta} e(\nu)M_{n-h(\nu),}$$
where $$h(\nu) \seteq \dfrac{1}{2}\sum_{1\le a< b \le \height(\beta)}\dc(\al_{\nu_a}, \al_{\nu_b}).$$

Then $K_\dc$ is an equivalence of categories from $\Modg(R(\beta))[q^{1/2}]$ to $\Modg(R_\dg(\beta))[q^{1/2}]$  which commutes with the grading shift functor (\cite[Lemma 1.5]{KP18}).  Moreover, 
$$K_\dc(M\conv N) \simeq q^{\frac{1}{2} \dc(\beta ,\gamma)}K_\dc(M) \conv K_\dc(N) \ \text{ for $M\in \Modg\bl R(\beta)\br$ and $N\in \Modg\bl R(\gamma)\br$.}$$

For $M\in \Modg(R_\dg(\beta))$, the graded dual vector space
$$M^\star \seteq  \HOM_{\bR}(M, \bR) = \soplus_{m\in \Z}\Hom_\cor(q^m M,\cor)$$
is an  $R_\dg(\beta)$-module via  the antiautomorphism  $\psi$ of $R(\beta)$ which fixes the generators
 $e(\nu)$, $x_k$, and $\tau_k$'s. 
Here the grading on $M^\star$ is given by
$$(M^\star)_d\seteq \soplus_{\nu \in I^\beta} \Hom_\cor\bl e(\nu)q^{d-2h(\nu)}M, \cor\br \quad \text{for} \ d\in\Z.$$
We have $(q M)^\star \simeq q^{-1} M^\star$ and  $K_\dc(M^\star) \simeq \bl K_\dc(M)\br^\star$ for $M\in \Modg(R(\beta))$.
 Then we have 
$$(M\conv N)^{\star} \simeq q^{(\wt(M),\wt(N))} N^\star  \conv M^\star \qt{in $\Modg(R_\dg)$}.$$

We say that $M \in R_\dg(\beta)\gmod$ is \emph{self-dual} if $M \simeq M^\star$ as $R_\dg(\beta)$-modules.
Note that a simple module  $M\in R(\beta)\gmod$ is self-dual if and only if $K_\dc(M)$ is a self-dual simple module in$R_\dg(\beta)\gmod [q^{1/2}]$.

For each $i\in I$ and $n\ge 1$, $R_\dg(n\al_i)$ has a unique self-dual simple module which is denoted by $L(i^n)$. Sometimes we will denote it by $\ang{i^n}$. 
Note that $\ang{i^n} =q_i^{n(n-1)/2}  \ang{i}^{\circ n}$ and we have $\dim_\cor(\ang{i^n})=n\ms{1mu}!$.

We shall denote by $\ang{k^cj^d}$ the simple head of $\ang{k^c}\conv \ang{j^d}$ for $j \neq k\in I$  and $c,d\in \Z_{\ge 0}$.
 Note that $\ang{k^c} \conv \ang{j^d}$ has a simple head, since $\ang{k^c}$ and $\ang{j^d}$ are unmixed (see section \ref{subsubsec:degree}).
Moreover $\ang{k^cj^d}$ is a self-dual simple in $R(c \al_k+d \al_j)\gmod$.
\vskip 1em

\subsubsection{Degrees of $R$-matrices}  \label{subsubsec:degree}
An ordered pair $(M,N)$ in $R_\dg(\beta)\gmod$ is called \emph{$\La$-definable} if $\HOM_{\Modg(R_\dg)} (M\conv N,N\conv M)=\cor \rmat{}$ for some non-zero $\rmat{}$.  
In this case we  call $\rmat{}$ the \emph{R-matrix} from $M\conv N$ to $N\conv M$, and
define $$\La(M,N)\seteq  \deg(\rmat{})$$ 
and we set
  $$\tLa(M,N)\seteq  \dfrac{1}{2}  \bl \La(M,N)-\dg(\wt(M),\wt(N))  \br.$$

For  $M\in R_\dg(\beta)\gmod$ and $N\in R_\dg(\gamma)\gmod$,  by applying $K^{-1}_\dc$, we obtain that 
\eqn&& 
\Hom_{R_\dg\gmod}(q^{d}M\conv N, N\conv M)  \\
&&\simeq \Hom_{R\gmod}(
q^{d-\frac{1}{2} \dc(\beta ,\gamma)}K^{-1}_\dc(M)\conv K^{-1}_\dc(N),q^{-\frac{1}{2} \dc(\gamma ,\beta)}  K^{-1}_\dc(N)\conv K^{-1}_\dc(M)) 
   \qt{for all $d\in \Z$. }
\eneqn

In particular, if $(M,N)$ is a $\La$-definable pair,  then we have
\eq \label{eq: dg_Lambda}
\La(M,N) = \La_{\dg_\can}(K^{-1}_\dc(M),K^{-1}_\dc(N))+\dc(\beta,\gamma),\eneq
where $\La_{\dg_\can}$ denotes the degree of $\rmat{}$ in the category $\Modg(R)$.
It follows that
\eqn
\text{ $\tLa(M,N)$ does not depend on the choice of $\dg$.}
\eneqn
If self-dual simple modules $M$ and $N$ in $R_\dg\gmod$ emph{strongly commute},  i.e., $M\conv N$ is simple,  then  $q^{\frac{1}{2}(\La(M,N)+(\wt(M),\wt(N)))} M\conv N$ 
is self-dual in $R_\dg\gmod$.

A pair $(M,N)$ is called \emph{$\de$-definable} if both $(M,N)$ and $(N,M)$ are $\La$-definable.  
In this case we define $$\de(M,N) \seteq \dfrac{1}{2}(\La(M,N)+\La(N,M)).$$
Then \eqref{eq: dg_Lambda} implies that 
\eqn
\text{ $\de(M,N)$ does not depend on the choice of $\dg$.}
\eneqn

\vskip 1em

We say that an ordered pair $(M,N)$ of $R$-modules is {\em unmixed}
if $$\sgW(M)\cap\gW(N)\subset\{0\}$$ 
where 
$$\text{$\gW(M)\seteq\set{\gamma\in \rtl_+}{e(\gamma,\beta-\gamma)M \neq 0}$ and $\sgW(M)\seteq\set{\gamma\in \rtl_+}{e(\beta-\gamma,\gamma)M \neq 0}$}$$ for an $R(\beta)$-module $M$.
If $(M,N)$ is an unmixed pair of simple modules, then $(M,N)$ is $\La$-definable (\cite[Proposition 2.2]{loc2}, \cite[Lemma 2.6]{TW16}). 
Indeed, we have
$\HOM_{\Modg(R_\dg)}(M\conv N,N\conv M) = \cor \rmat{}$ and 
the homomorphism $\rmat{}$ is given by
$$\rmat{} (w\tens v) = \tau_{w[n,m]} (v\tens w)\quad
\text{for $w \in M$ and $v \in N$,}$$
 where  $m=\height(\wt(M))$, $n=\height(\wt(N))$.
It follows that 
\eq \La(M,N)= \dg(\wt(M),\wt(N)) \qt{and hence} \quad \tLa(M,N) =0.\label{eq:unmixed}\eneq
For example, we have
\eqn
\La(\ang{j},\ang{k}) = \deg \tau_1 e(k,j) = \dg(\al_j,\al_k)
\qt{for $j,k \in I$ such that $j\not=k$.}
\eneqn

\vskip 1em
For $\beta \in \rlQ_+$ and $i\in I$, let
\begin{align} \label{Eq: def of p}
\mathfrak{p}_{i, \beta}  \seteq \sum_{\nu \in I^\beta} \Bigl(\hs{1ex}  \prod_{a \in \{1, \ldots, \Ht(\beta) \},\ \nu_a=i} x_a \Bigr) e(\nu)\in R_\dg(\beta).
\end{align}
Then $\mathfrak{p}_{i, \beta} $ belongs to the center of $R_\dg(\beta)$.

Let $M$ be a simple module in $R_\dg(\beta)\gmod$.  
 Assume that there exists   an $R_\dg(\beta)$-module $\Ma$ with an endomorphism $z_{\Ma}$ of $\Ma$
with degree $d_{\Ma} \in \Z_{>0}$ such that
\eq \label{eq:oldaff}
&&\hs{2ex}\parbox{75ex}{
\begin{enumerate}[\rm (i)]
\item $\Ma / z_{\Ma} \Ma \simeq M$,
\item $\Ma$ is a finitely generated free module over the polynomial ring $\bR[z_{\Ma}]$,
\item $\mathfrak{p}_{i,\beta} \Ma \ne 0$ for all $i\in I$.\label{it:nonzeroP}
\end{enumerate}
}\eneq
Then we say that $M$ admits a \emph{\ro strict\rf affinization} $\Ma$.
A simple module $M$ is called \emph{\afr} if $M$ is real and admits an affinization. 
For example, $R(\al_i)$ is an affinization of $L(i)\simeq R(\al_i)/R(\al_i)x_1e(i)$ with $z_{R(\al_i)} = x_1$ of degree $(\al_i,\al_i)$. 

Let $(M,N)$ be a pair of simple modules such that one of them is \afr.
Then $(M,N)$ is $\de$-definable (\cite[Proposition 2.10]{KP18}) and 
$\HOM_{\Modg(R_\dg)}(M\conv N ,N\conv M)=\cor \rmat{}$.
where $\rmat{}= z_{\Ma}^{-s} \Runi_{\Ma,N} \vert_{z_{\Ma}=0}$
when $M$ is \afr with an affinization $(\Ma,z_\Ma)$.  Here  $s\in\Z_{\ge 0}$ is the largest non-negative integer such that $ \Runi_{\Ma,N} (\Ma \conv N) \subset z_{\Ma}^{s} N \conv \Ma  $, called the \emph {order of zero of  $\Runi_{\Ma,N}$} In this case $\Rre_{\Ma,N}\seteq z_{\Ma}^{-s} \Runi_{\Ma,N} $ is called the \emph{renormalized R-matrix}.

\subsection{Braiders and localizations of monoidal categories}
In this subsection, we recall the notions of braiders and localization of monoidal categories
introduced in \cite{loc1}.

\subsubsection{Localizations of monoidal categories}
Let $\bR$ be a field. 
A $\bR$-linear monoidal category $\catT$ is  \emph{graded} if  there exist a $\cor$-linear auto-equivalence $q$ on $\catT$ and isomorphisms
$$q(X\tens Y)\simeq (qX)\tens Y\simeq X\tens(qY)\qt{functorial in $X$, $Y\in\catT$}$$ such that
the following diagrams commute for any $X,Y,Z\in\catT$:
\eq\hs{3ex}
\ba{l}
\xymatrix@R=3ex{q(X\tens Y\tens Z)\ar[r]\ar[d]&q(X\tens Y)\tens Z\ar[d]\\
  (qX)\tens (Y\tens Z)\ar[r]&(qX\tens Y)\tens Z,}\quad
\xymatrix@R=3ex{q(X\tens Y\tens Z)\ar[r]\ar[d]&q(X\tens Y)\tens Z\ar[d]\\
  X\tens q(Y\tens Z)\ar[r]&X\tens(q Y)\tens Z,}\\
\xymatrix@R=3ex{q(X\tens Y\tens Z)\ar[r]\ar[d]&(X\tens Y)\tens qZ\ar[d]\\
X\tens q(Y\tens Z)\ar[r]&X\tens (Y\tens qZ).}\\
\ea\label{eq:graded_monoidal}
\eneq 
By identifying $q$ with the invertible central object $q \one \in \catT$, we have $qX \simeq q \tens X$ for $X\in \catT$.
We write $q^n$ ($n\in\Z$) for $q^{\tens n}$ for the sake of simplicity. 

Let  $\Lambda$ be a $\Z$-module. A $\bR$-linear monoidal category $\catT$ is \emph{$\Lambda$-graded} if $\catT$ has a decomposition
$ \catT = \bigoplus_{\lambda \in \Lambda} \catT_\lambda $ such that $\triv \in \catT_0$ and $\otimes$ induces a bifunctor $\catT_{\lambda} \times \catT_{\mu} \rightarrow \catT_{\lambda+\mu}$
for any $\lambda, \mu \in \Lambda$.

\begin{df} \label{def:graded braider}
A \emph{graded left braider} is a triple $(C, \coRl_C, \dphi)$ consisting of an object $C$, a $\Z$-linear map $\dphi\cl  \Lambda \rightarrow \Z$ and a morphism
$$
\coRl_C(X) \cl  C \otimes X \longrightarrow q^{-\dphi(\lambda)} \otimes X \otimes C
$$
functorial in $X\in\catT_\la$  such that the diagrams 
$$
\xymatrix{
C \otimes X \otimes Y   \ar[rr]^{\coRl_C(X)\otimes Y}  \ar[drr]_{\coRl_C(X \otimes Y) \ \ }  &  &   q^{-\dphi(\lambda)} \otimes X \otimes C \otimes Y  \ar[d]^{ X \otimes \coRl_C(Y)}   \\
  & &   q^{ -\dphi(\lambda+\mu) } \otimes X \otimes Y \otimes C
}
\qtq 
\xymatrix{
C \otimes \triv   \ar[rr]^{R^l_C(\triv)}  \ar[drr]_{ \simeq }  &  &   \triv \otimes C   \ar[d]^{ \wr}   \\
& &    C
}
$$
commute for any $X \in \catT_\lambda$ and $Y \in \catT_{\mu}$.

A graded left braider $ ( C, \coRl_{C}, \phi  )$ is  a \emph{central object} if $\coRl_{C}(X)$ is an isomorphism for any $X \in \catT$.

\end{df}

 Recall that the category $\catTc^l$ consisting of graded left braiders in $\catT$ is a monoidal category with  a canonical faithful monoidal functor $\catTc^l \rightarrow \catT$.

 \vskip 2em
Let $I$ be an index set and let $ \st{(C_i, \coRl_{C_i}, \dphi_i )}_{i\in I} $ be a 
family of graded left braiders.
We say that $\st{(  C_i ,  \coRl_{C_i} ,  \dphi_i )}_{i\in I}$ is a \emph{real commuting family of graded left braiders} in $\catT$ if
\bna
\item  \label{Eq: 1 in grcf} $ C_i \in \catT_{\lambda_i}$  for some $\lambda_i \in \Lambda$, and
$\dphi_i(\lambda_i) = 0$, $\dphi_i( \lambda_j ) + \dphi_j( \lambda_i ) = 0$,
\item  $\coRl_{C_i}(C_i) \in \bR^\times \id_{C_i \otimes C_i}$ for  any$i\in I$,
\item   \label{Eq: 2 in grcf} $\coRl_{C_j}(C_i) \circ \coRl_{C_i}(C_j) \in \bR^\times \id_{C_i \otimes C_j}$ for any  $i,j\in I$.
\end{enumerate}

Set
$$
\lG \seteq  \Z^{\oplus I} \ \text{with the natural basis $\{ e_i \mid i\in I \} $, \ and } \quad \lG_{\ge 0} \seteq  \Z_{\ge 0}^{\oplus I}.
$$
For each $\al \in \Gamma_{\ge 0},$ there exists a graded left braider $C^\al=(C^\alpha, R_{C^\alpha}, \phi_\alpha )$, and  for $\alpha, \beta \in \lG_{\ge0}$, 
there exists an isomorphism $\xi_{\alpha, \beta}\cl  C^\alpha \otimes C^\beta \buildrel \sim\over \longrightarrow q^{H(\al,\beta)}\tens C^{\alpha+\beta}$ in $\catTc^l$  which satisfy several compatible relations, where $H(\cdot, \cdot)$ is a bilinear form on $\Gamma$ (see {\cite[Lemma 2.3, Lemma 1.16]{loc1}}).

\begin{thm}[{\cite[Theorem 1.4, Proposition 1.5]{loc1}}]\label{Thm: graded localization}
Let $ \st{C_i=(C_i, R_{C_i}, \dphi_i )}_{i\in I} $ be a 
real commuting family of graded left braiders  in $\catT$. 
Then there exist a monoidal category $\lT$, denoted also by $\catT[C_i^{ \otimes-1} \, \vert \,  i\in I]$, 
a monoidal functor 
$\Upsilon\cl \catT \to \lT$ 
and a real commuting family of graded braiders $\st{\tC_i=(\tC_i, R_{\tC_i}, \phi_i)}_{i\in I}$ in $\lT$ satisfy the following properties:

\bnum
\item
for $i\in I$, $\Upsilon(C_i) $ is isomorphic to $ \widetilde{C}_i$ and it is invertible in $(\lT)_{\,\mathrm{br}}$, 
\item 
for $i\in I$ and $X\in\catT_\la$, the diagram
$$
\xymatrix{
\Upsilon(C_i \otimes X)  \ar[r]^\sim  \ar[d]_{\Upsilon( R_{C_i} (X)  )\ms{10mu}}^-\bwr  & \widetilde{C}_i \otimes \Upsilon(X) \ar[d]_{ R_{\widetilde{C}_i} (\Upsilon(X)  )\ms{5mu}} ^-\bwr \\
\Upsilon(q^{-\phi_i(\la)}\tens X \otimes  C_i )  \ar[r]^\sim & q^{-\phi_i(\la)}\tens  \Upsilon(X)\otimes  \widetilde{C}_i
}
$$
commutes.
\setcounter{myc}{\value{enumi}}

\end{enumerate}

Moreover, the functor $\Upsilon$ satisfies the following universal property:
\bnum\setcounter{enumi}{\value{myc}}
\item  If there are another $\La$-graded monoidal category $\catT'$ 
with an invertible central object $q\in\catT'_0$
 and   a $\La$-graded  monoidal functor $\Upsilon'\cl  \catT \rightarrow \catT'$ 
such that  
\bna
\item  $\Upsilon'$ sends the central object $q\in\catT_0$ to $q\in\catT'_0$, 
\item   \label{Eq: loc 1}
$\Upsilon'(C_i) $ is invertible in $\catT'$ for any $i\in I$ and
\item 
for any $i\in I$ and $X\in\catT$, $\Upsilon'(R_{C_i}(X))\cl
\Upsilon'(C_i\tens X)\to\Upsilon'(q^{-\phi_i(\la)}\tens X\tens C_i)$ is an isomorphism,
\end{enumerate}
then there exists a monoidal functor $\mathcal F$, which is unique up to a unique isomorphism,  such that
the diagram
$$
\xymatrix{
\catT \ar[r]^{\Upsilon} \ar[dr]_{\Upsilon'}  & \lT \ar@{.>}[d]^{\mathcal F }\\
& \catT'
}
$$
commutes.

\item Assume further that $\catT$ is an abelian category and that the tensor product $\tens$ on $\catT$ is exact. Then 
the category $\lT$ is an abelian category with exact tensor product $\tens$, and the functor $\Upsilon\cl \catT \to \lT$ is exact. 

\end{enumerate}

\end{thm}

\begin{remark} \label{rem:changes}
\begin{enumerate}
\item
Recall that a \emph{graded right braider in $\catT$} is a triple $(C, \coRr_C, \phi)$ of an object $C$, a $\Z$-linear map $\phi\cl  \La \to \Z$, and a morphism functorial in $X \in \catT_\la$
\eqn
\coRr_{C} (X) \cl  X \tens C  \to q^{-\phi(\la)} C\tens X
\eneqn 
with analogous conditions to  the ones in Definition \ref{def:graded braider}. 
The theorem above also holds for a real commuting family of graded right braiders (see \cite[Theorem 2.4]{loc3}).
\item  We changed $\phi$  in \cite[Section 2.3]{loc1} with $-\phi$ in this paper in order to have $\phi(\la)=\deg(\coRl_C(X))$ when $X \in \catT_\la$. 
\end{enumerate}
\end{remark}

\vskip 2em 
\subsubsection{Non-degenerate braiders in $R\gmod$.}
Let $R$ be a quiver Hecke algebra and 
let $\dg$ be a grade associator.
Recall that  $R_\dg(\al_i)$ is an affinization of $L(i)= R_\dg(\al_i)/R_\dg(\al_i)x_1e(i)$ with $z_{R(\al_i)} = x_1$  for each $i\in I$.
\Prop [{\cite[Proposition 4.1]{loc1}}] \label{prop:nondeg braider}
Let $C$ be a simple $R_\dg(\beta)$-module,
Then there exists a graded left braider $\bl C, \coRl_C,\phi\br$ in $R_\dg\gmod$ such that 
\eqn 
\coRl_C(R_\dg(\al_i))=\Rre_{C,R_\dg(\al_i)} \ \quad \text{for every $i\in I$ up to a constant}.
\eneqn
Similarly, there exists a  graded right braider $\bl C, \coRr_C,\phi\br$ in $R_\dg\gmod$ such that
\eqn 
\coRr_C(R_\dg(\al_i))=\Rre_{R_\dg(\al_i),C} \ \quad \text{for every $i\in I$ up to a constant}.
\eneqn
\enprop

We call the braiders $\coRl_C$ and  $\coRr_C$ in the proposition above \emph{non-degenerate}, since $\coRl_C(L(i))$
and $\coRr_C(L(i))$ are  non-zero homomorphisms  for each $i\in I$. Note that  a non-degenerate graded left (respectively, right) braider associated with $C$ is essentially unique (\cite[Lemma 4.2]{loc1}).

Note that if the braider $\bl C, \coRl_C,\phi^\ell \br$ (respectively, $\bl C, \coRr_C,\phi^r \br$) is non-degenerate then we have $\phi^\ell(-\al_i) = \La(C, \ang{i})$ (respectively, $\phi^r(-\al_i) = \La(\ang{i},C)$) for $i\in I$.

\vskip 2em
Note that the category $\lT$ in Theorem \ref{Thm: graded localization} is explicitly given in \cite{loc1}.
We recall the category $\lT$ when $\T$ is a monoidal subcategory of $R_\dg\gmod$ and the real commuting family is a singleton $(C,\coRl_C,\phi)$ of non-degenerate left graded braider in $R_\dg\gmod$ for $C \in \T$:
Set $\wt(C)\seteq \eta \in\rtl_-$ and $\T_{-\beta} \seteq \T \cap R_\dg(\beta)\gmod$ for $-\beta\in \rtl_-$.
The objects and morphisms of $\lT$ are given as follows:
\begin{align*}
\Ob (\lT) &\seteq  \Ob(\catT) \times \Z, \\
\Hom_{\lT}( (X, m), (Y, n) ) &\seteq\hs{-5ex}   \indlim_{  \substack{l \in \Z_{\ge 0} \\   l+m, \ l+n\ge 0, \   \lambda + m\eta = \mu + n\eta  }     }\hs{-5ex}\Hom_{\catT}(  C^{l+ m}\conv X,  q^{ \gH(l, n-m) - (l+n)\dphi(\mu)  } Y \conv C^{ l + n} )
\end{align*}
where $X \in \catT_{\lambda}$ and $Y \in \catT_{\mu}$,  and $H(m,n)\seteq-\tLa(C^{\circ m},C^{\circ n})=- \dfrac{mn(\wt(C),\wt(C))}{2}\in \Z$.
Note that $q^{-H(m,n)}C^m \conv C^n\simeq C^{m+n}$ is self-dual in $R_\dg\gmod$ for any $m,n\in\Z_{\ge0}$  if $C$ is self-dual. 

For the composition of morphisms in $\lT$ and the associativity of the composition, see \cite[Section 2.2, Section 2.3]{loc1}.

 By the construction, we have the decomposition
$$
\lT = \bigoplus_{\mu \in \rtl_-} \lT_{\mu}, \qquad \text{where }  \lT_{\mu} \seteq  \{ (X, m) \mid X \in \catT_\lambda, \ \lambda + m \eta =\mu  \}.
$$

 For $m, n \in \Z $, $X \in \catT_\lambda$ and $Y \in \catT_\mu$,
we have
$$
(X, m) \conv(Y, n) \seteq  ( q^{\gH(m, n)+n\dphi(\lambda)}  X \conv Y,m+n  ).
$$ 
For the tensor product of morphisms,  see \cite[Section 2.2, Section 2.3]{loc1}.

The localization of $\T\subset R_\dg\gmod$ via a non-degenerate graded right braider $(C, \coRr_C,\phi)$ has an analogous description (see, \cite[Section 2.1]{loc3}).


The statement (ii) in the proposition below is proved in
\cite[Theorem 4.12]{loc1}.
Since the statement (i) is similarly proved, we do not repeat the proof.

\Prop [{\cite[Theorem 4.12]{loc1}}] \label{prop: right dual}
Let $(C, \coRl_C,\phi)$ be a non-degenerate graded {\em real} left braider in $R_\dg\gmod$ and let $Q\cl R_\dg\gmod \to R_\dg\gmod[C^{\circ -1}  ] $ be the localization functor and set $\tC\seteq Q(C)$.
Let $i\in I$.
For  $\ell\in\Z_{>0}$, set
$L_{\ell}(i) \seteq L(i)_z / z^\ell L(i)_z$ where $L(i)_z$ denotes the affinization of $L(i)$ of degree $2\sfd_i$. 
\bnum
\item
Assume that $\eps_i(C)=1$.
For $\ell \in \Z_{>0}$, we set
$K_\ell \seteq \E_i(C^{\circ \ell})$.
Then there is an epimorphism 
$$L_{\ell}(i) \circ K_\ell \epito C^{\circ \ell}. $$
Moreover, $Q(K_\ell) \circ \tC^{\circ -\ell}$ is a right dual of $Q(L_\ell(i))$ in $R_\dg\gmod[C^{\circ -1}]$.

\item
Assume that $\eps^*_i(C)=1$.
For $\ell \in \Z_{>0}$, we set
$K_\ell \seteq \E^*_i(C^{\circ \ell})$.
Then there is an epimorphism 
$$K_\ell\circ L_{\ell}(i)\epito C^{\circ \ell}, $$

Moreover, $\tC^{\circ -\ell}\conv Q(K_\ell)$ is a left dual of $Q(L_\ell(i))$ in $R_\dg\gmod[C^{\circ -1}]$.
\ee
\enprop

 Note that, in case (i), $C$ and $\ang{i}$ commute by
\cite[Proposition 4.11 (i)]{loc1}.

\subsection{Affinizations and Schur-Weyl duality functors}
In this subsection, we recall  results in \cite{refl}.
\subsubsection{Affine objects in an abelian category}
Let $\cor$ be a base field. 
Throughout this subsection, we consider a $\cor$-linear category $\shc$ satisfying
\eqn&&\hs{2ex}\left\{
\parbox{75ex}{\be[{$\bullet$}]
\item $\shc$ is abelian,
\item $\shc$ is $\cor$-linear graded with an auto-equivalence $q$, 
\item any object has a finite length,
\item  any simple object $S$ is absolutely simple, i.e.\ $\cor\isoto\END_\shc(S)
\seteq\HOM_\shc(S,S)$.\ee
}\right.
\label{cond:fcat}
\eneqn
Here $\HOM_\shc(M,N)\seteq\soplus_{n\in\Z}\HOM_{\shc}(M,N)_n$
with $\HOM_{\shc}(M,N)_n\seteq\Hom_{\shc}(q^nM,N)$.

Let $\Pro(\shc)$ be the category of pro-objects of $\shc$ (see \cite[Section 6.1]{KS} for the definition). Then $\Pro(\shc)$ is  equivalent to the opposite category of the category of $\cor$-linear  left exact functors from $\shc$ to the category of $\cor$-vector spaces. 
In particular, $\Pro(\shc)$ is $\cor$-linear and abelian such that the canonical fully-faithful functor  $\shc \to \Pro(\shc)$ is $\cor$-linear and exact.
The category $\Pro(\shc)$ admits small codirected projective limits, denoted by $\proolim$ (\cite[Proposition 6.1.8]{KS}), and every object in $\Pro(\shc)$ is isomorphic to the projective limit $\proolim[i \in I]X_i$ of a small codirected projective system $\{X_i\}_{i\in I}$ in $\shc$.

Let $A=\soplus_{k\in \Z} A_k$ be a  finitely generated,  graded commutative $\cor$-algebra such that 
 $A_0\simeq\cor$, and $A_k=0$ for any $k<0$.
For each $m\in \Z$, set 
$A_{\ge m} \seteq\soplus _{k\ge m} A_k.$

Let $\Modg(A,\Pro(\shc))$ be the category of graded $A$-modules in $\Pro(\shc)$, and
let $\Proc(A, \shc)$ be the full subcategory of $\Modg(A,\Pro(\shc))$ consisting of objects $\Ma$ such that
$\Ma /A_{> 0}  \Ma \in \shc$ and $\Ma \isoto \proolim[k] \Ma/ A_{\ge k }\Ma$, where $A_{\ge k }\Ma \seteq \Im(A_{\ge k} \tens_A \Ma \to A\tens_A \Ma)$.

Let $z$ be an indeterminate of homogeneous degree $d \in \Z_{> 0}$. 
Let  $\Aff(\shc)$ be the full subcategory of $\Proc(\cor[z],\shc)$ consisting of objects $\Ma$ such that $z\in \END_{\Pro(\shc)}(\Ma)$ is a monomorphism.  We call the pair $(\Ma,z)$ of $\Aff(\shc)$ and $z \in \END_{\Pro(\shc)}(\Ma)_d$ an \emph{affine object} in $\shc$.

\Prop\label{prop:qutaff}
Let $\dX\in\Proc(\cor[z],\shc)$ such that 
$\dX/z\dX$ has a simple head $L$.
Assume that $q^nL$ does not appear in $\Ker(\dX/z\dX\to L)$ for any $n\in\Z_{>0}$.
Then an epimorphism 
$\dX\epito \dL$
in $\Proc(\cor[z],\shc)$
for $\dL\in\Aff(\shc)$ such that $\dL/z\dL$ is simple
is unique \ro if it exists\rf.
\enprop
\Proof
Set $\dX_\tor=\bigcup_{n\in\Z_{\ge0}}\Ker(z^n)$.
Replacing $\dX$ with $\dX/\dX_\tor$, we may assume from the beginning
that $\dX\in\Aff(\shc)$. 
Let $0\to \dY_k\to\dX\To[f_k] \dL_k\to 0$
be an exact sequence where $\dL_k\in \Aff(\shc)$ with
$\dL_k/ z \dL_k\simeq L$ ($k=1,2$).
It is enough to show that $\dY_1=\dY_2$.

If $f_1(\dY_2)=0$, then $f_1$ factors as $\dX\To[f_2 ]\dL_2\to\dL_1$ and there is an monomorphism 
$Y_2 \to Y_1$.  
The epimorphism $\dL_2\to \dL_1$ induces an epimorphism $\dL_2/z \dL_2 \to \dL_1/z \dL_1$ in $\shc$.  Since $\dL_k/z \dL_k\simeq L$ for $k=1,2$,  $\dL_2/z \dL_2 \to \dL_1/z \dL_1$ is an isomorphism  and hence  $\dL_2\to \dL_1$ is an isomorphism by \cite[Lemma 3.5]{refl}. 
Thus $Y_2 \to Y_1$ is an isomorphism.

Hence we may assume that $f_1(\dY_2)\not=0$.
Let $n_1$ be the largest non-negative integer such that $f_1(\dY_2) \subset z^{n_1}\dL_1$. 
Then we obtain
$\dY_2/z\dY_2 \To[f_1\vert_{z=0}] z^{n_1}\dL_1/z^{1+n_1}\dL_1
\simeq q^{\deg(z)n_1}L$.
Since $L$ is simple , $f_1\vert_{z=0}$ is an epimorphism. Then by \cite[Lemma 3.5]{refl},  $f_1\cl \dY_2 \to z^{n_1}\dL_1$ is an epimorphism.  
Hence we have
$f_1(\dY_2)=z^{n_1}\dL_1$.
Since $\dY_2/z\dY_2\simeq\Ker(\dX/z\dX\to\dL_2/z\dL_2)$, 
$q^{n_1}L$ appears in $\Ker(\dX/z\dX\to L)$ and hence
we have $n_1=0$. 
Hence we obtain $\dX=\dY_1+\dY_2$.
On the other hand, since $\dX /z\dX$ has a simple head $L$,
$R\seteq\Ker(\dX   \to\dL_k/z\dL_k)=\dY_k+z\dX$ does not depend on $k$.
Thus we have 
$\dY_1+\dY_2\subset R$, which is a contradiction. 
\QED
A simple object $M \in \shc$ is called \emph{real} if $M\tens M$ is simple.

\Ex
Let $\shc=\cor[x]\gmod$ with $\deg(x)=d>0$,
$\dX=\cor[z,x]/\cor[z,x](z^2-x^2)\in\Aff(\shc)$, with $\deg(z)=d$.
Then, $\dX/z\dX\simeq\cor[x]/\cor[x]x^2$ has a simple head $L=\cor[x]/\cor[x]x$.
Let $\dL_\pm=\cor[z,x]/\cor[z,x](z\pm x)\in\Aff(\shc)$.
Then,
we have $\dX\epito \dL_\pm$.
Note that $[\dX/z\dX]=(1+q^d)[L]$.
\enex

\vskip 1em
We denote by $\Rat(\shc)$ the category 
with $\Ob\bl\Proc(\cor[z],\shc)\br$  
as the set of objects and with the morphisms defined as follows. 
For $\Ma,\Na\in\Ob\bl\Rat(\shc)\br$,
\eqn\Hom_{\,\Rat(\shc)}(\Ma,\Na)&&=\indlim_{k\in\Z_{\ge0}}
\Hom_{\,\Proc(\cor[z],\shc)}\bl(z^k\cor[z])\tens[{\cor[z]}]\Ma,\Na\br\\
&&\simeq
\indlim_{k\in\Z_{\ge0}}
\Hom_{\,\Proc(\cor[z],\shc)}(\Ma,\cor[z]z^{-k}\tens_{\cor[z]}\Na).\eneqn
Then we have
$$\HOM_{\Rat(\shc)}(\Ma,\Na)\simeq
\cor[z,z^{-1}]\tens_{\cor[z]}\HOM_{\Proc(\cor[z],\shc)}(\Ma,\Na)
\qt{for $\Ma,\Na\in\Proc(\cor[z],\shc)$}$$
so that every object of $\Rat(\shc)$ is a $\cor[z,z^{-1}]$-module, i.e.,
$z\in\END_{\Rat(\shc)}(\Ma$) is invertible for any $\Ma\in\Rat(\shc)$.

\Def\label{def:ratisom}
We say that a morphism
$f\cl \Ma\to \Na$ in $\Proc(\cor[z],\shc)$ or $\Rat(\shc)$ is a
{\em rational isomorphism} if it is an isomorphism in $\Rat(\shc)$.
\edf
 Hence, $f\in\Hom_{\Proc(\cor[z],\shc)}(\Ma,\Na)$ is a rational isomorphism if and only if there exist $n\in\Z_{\ge0}$ and
a morphism
$g\cl \Na\to \Ma$ in $\Proc(\cor[z],\shc)$ such that
$g\circ f=z^n\id_{\Ma}$ and $f\circ g=z^n\id_{\Na}$. 

\subsubsection{Rational centers and affinizations}
From now on, we assume further  that $\shc$ is a monoidal category with a tensor product $\tens$ such that 
\eq \label{eq:shc mon}
&&\left\{\parbox{70ex}{
\bna
\item $\tens$ is $\cor$-bilinear and bi-exact,
\item the unit object $\one$ satisfies $\END_\shc(\one)\simeq\cor$, 
\item  $\shc$ is a graded monoidal category with $q$ as a grading
  shift functor (see\eqref{eq:graded_monoidal}). 
\setcounter{mycc}{\value{enumi}}
\ee}\right.
\label{cond:exactmono}
\eneq
Then $(\Pro(\shc), \tens) $ is a monoidal category on which $\tens$ is bi-exact.
Note that there are  bi-exact functors 
$$\tens\;\cl \shc\times\Rat(\shc)\to \Rat(\shc)
\qtq \tens\;\cl \Rat(\shc)\times\shc\to \Rat(\shc).$$

In the rest of this section we assume further that
\eqn
&&\hs{1ex}\parbox{70ex}{
\bna
\setcounter{enumi}{\value{mycc}}
\item $\shc$ has a decomposition $\shc=\soplus_{\la\in\Idx}\shc_\la$
 where $\Idx$ is an abelian monoid such that
$\one\in \shc_\la$ with $\la=0$,
$\shc_\la\tens\shc_\mu\subset\shc_{\la+\mu}$ for any $\la,\mu\in\Idx$,  
and $\shc_\la$ is stable by $q$. 
\ee
} \label{cond:g}
\eneqn
\Def\label{def:rational center}
A {\em rational center} in $\shc$ is a triple $(\Ma,\phi,\Rmat_\Ma)$
of $\Ma\in\Rat(\shc)$, an additive map $\phi\cl\Idx\to\Z$ and an isomorphism 
$$\Rmat_\Ma(X)\cl q^{\tens\,\phi(\la)}\tens\Ma\tens X\isoto X\tens\Ma$$
in $\Rat(\shc)$ functorial in $X\in\shc_\la$ such that
$$
\xymatrix@C=15ex
{q^{\tens\,\phi(\la+\mu)}\tens\Ma\tens X\tens Y\ar[r]^{\Rmat_\Ma(X)\tens Y}\ar[dr]_{\Rmat_\Ma(X\tens Y)}&q^{\tens\,\phi(\mu)}\tens X\tens\Ma\tens Y\ar[d]^{X\tens\Rmat_\Ma(Y)}\\
&X\tens Y\tens\Ma}$$
and $$\xymatrix@C=12ex{\Ma\tens\one\ar[r]^{\Rmat_\Ma(\one)}\ar[dr]^-{\sim}&\one\tens \Ma\ar[d]^\bwr\\
&\Ma}
$$
commute in $\Rat(\shc)$ for any $X\in\shc_\la$ and $Y\in\shc_\mu$ ($\la,\mu\in\Idx$).
\edf
{\em In the sequel we neglect grading shift.}

\Def\label{def:affinization}
An {\em affinization} $\Ma$ of $M\in\shc$ is an \afn
$(\Ma,\z)$ with a rational center $(\Ma,\Rmat_\Ma)$ and an isomorphism
$\Ma/\z\Ma\simeq M$.

We call a simple object $M$ in $\shc$ \emph{\afr} if $M$ is real, and $M$ admits an affinization $\Ma$.

\edf

\begin{remark}
Let $M$ be a simple module in $R_\dg(\beta)\gmod$. Assume that 
an $R_\dg(\beta)$-module  $(\Ma,z)$  satisfies \eqref{eq:oldaff}. 
Then $(\Ma,\Runi_\Ma)$ is a rational center  in $R_\dg\gmod$ (\cite[Lemma 2.9]{KP18}, see also  \cite[Proposition 7.3]{refl}). 
Hence $\Ma$ is an affinization of $M$.  We call $(\Ma,\Runi_\Ma)$ a \emph{canonical affinization} of $M$ if $(\Ma,z)$  satisfies the conditions \eqref{eq:oldaff}. Note that there are affinizations in $R_\dg\gmod$ which are not canonical; see \cite[Remark 7.4]{refl}.
\end{remark}

\Th  [{\cite[Theorem 6.10]{refl}}]\label{th:ren_r_matrix}
Let $(\Ma,\Rmat_{\Ma})$ be an affinization of $M\in\shc$, and
let $(\Na,\zN)$ be an \afn of $N\in\shc$.
Assume that
$\dim \HOM_{\shc}(M\tens N,N\tens M)=1$.
Then there exist 
a homogeneous $f(\z,\zN)\in \bl\cor[\z,\z^{-1}]\br[[\zN]]$ and a morphism
$\Rre_{\Ma,\Na}\cl \Ma\tens \Na\to \Na\tens \Ma$ in $\Proc(\cor[\z,\zN],\shc)$
such that 
\bna
\item$\HOM_{\cor[\z,\zN]}(\Ma\tens\Na,\Na\tens \Ma)
=\cor[\z,\zN]\Rre_{\Ma,\Na}$,
\item
as an element of 
$\HOM_{\Rat(\shc)}\bl \Ma\tens(\Na/\zN^k\Na),\;(\Na/\zN^k\Na)\tens\Ma\br$,
we have $\Rre_{\Ma,\Na}\vert_{\Na/\zN^k\Na}=f(\z,\zN)\Rmat_\Ma(\Na/\zN^k\Na)$
 for any $k\in\Z_{>0}$, 
\item $\Rre_{\Ma,\Na}\bigm|_{\z=\zN=0}$ does not vanish,
\item $f(\z,\zN)\vert_{\zN=0}$ is a monomial of $\z$.
\ee
Moreover such $\Rre_{\Ma,\Na}$ and $f(\z,\zN)$ are unique up to constant multiples.
We call $\Rre_{\Ma,\Na}$ a \emph{renormalized R-matrix}.
\enth

\begin{definition}  
Let $M$ and $N$ be simple objects in $\shc$. 
\begin{enumerate}[(i)]
\item 
  We say that $(M,N)$ is {\em $\La$-definable} if
  $\dim\HOM_\shc(M\tens N,N\tens M)=1$.
  Then we define 
\eq
&&\La(M,N)\seteq\deg(\rmat{M,N})\in\Z,
\eneq where $\HOM_\shc(M\tens N,N\tens M)=\cor\rmat{M,N}$. 
\item  We say that $(M,N)$ is {\em $\de$-definable} if $(M,N)$ and $(N,M)$ are
  $\La$-definable.
  Then, we define
\eq
&&\de(M,N) \seteq  \dfrac{1}{2} \bl\La(M,N)+\La(N,M)\br\in \dfrac{1}{2}\Z_{\ge 0}. \label{eq:de}
\eneq 
\item 
Let $(\Ma, \Rmat_{\Ma})$ and  $(\Na, \Rmat_{\Na})$
be affinizations of $M$ and $N$, respectively. Assume that  $(M,N)$ is $\de$-definable and $\dim \END_\shc(M\tens N)=1$.
Then define $\De(\Ma,\Na) \in \cor[\z,z_{\Na}]  $ by
\eq
 \Rre_{\Na,\Ma} \circ \Rre_{\Ma,\Na}  = \De(\Ma,\Na) \id_{\Ma\tens \Na}. \label{eq:Daf}
\eneq
\ee
Note that  such $\De(\Ma,\Na)$ exists by \cite[Proposition 6.2]{refl},
and it is well-defined up to a constant multiple in $\cor^\times$.
If moreover $\dim\END_\shc(N\tens M)=1$, then we have $ \De(\Ma,\Na) = \De(\Na,\Ma)$.
\end{definition}

\Lemma\label{lem:LaXY}
Let $(M,N)$ be a $\La$-definable pair of simple objects such that
$X$ has a right dual $\D(X)$.
Then $(Y,\D(X))$ is $\La$-definable, and
$$\La(X,Y)=\La\bl Y,\D(X)\br.$$
\enlemma
\Proof
It immediately follows from
$\HOM(X\tens Y,Y\tens X)\simeq\HOM\bl Y\tens\D(X),\D(X)\tens Y\br$.
\QED

\subsubsection{Quasi-rigidity}
We denote by $M\hconv N$ and $M\sconv N$ the head and the socle  of $M\tens N$, respectively.

Recall that 
a monoidal category $\sha$ is \emph{quasi-rigid}
 if it satisfies:
\bna
\item $\sha$ is abelian and $\tens$ is bi-exact,
\item for any $L,M,N\in\sha$, $X\subset L\tens M$ and
$Y\subset M\tens N$ such that $X\tens N\subset L\tens Y\subset L\tens M\tens N$, there exists $K\subset M$ such that
$X\subset L\tens K$ and $K\tens N\subset Y$,
\item for any $L,M,N\in\sha$, $X\subset M\tens N$ and
$Y\subset L\tens M$ such that $L\tens X\subset Y\tens N$, there exists $K\subset M$ such that
$X\subset K\tens N$ and $L\tens K\subset Y$.
\ee
If $\sha$ is a rigid monoidal category with bi-exact tensor product, then it is quasi-rigid (\cite[Lemma 6.17]{refl}). 
Note that any monoidal subcategory of a quasi-rigid category  which is closed under taking subquotients is quasi-rigid.

One of the most useful properties of quasi-rigidity is the following
\Prop[{\cite[Proposition 6.19]{refl}}] \label{prop:simplehd}
Assume  that $\shc$ is \KO.
Let $M$ and $N$ be simple objects of $\shc$ and assume that one of them is \afr.
Then,
$M\tens N$ has a simple head and a simple socle.
Moreover, $(M,N)$ is $\de$-definable and
$$M\hconv N\simeq\Im(\rmat{M,N})\simeq N\sconv M$$
up to grading shifts. 
\enprop

By Proposition \ref{prop:qutaff},   Proposition \ref{prop:simplehd}, and \cite[Proposition 6.20]{refl}, we have
\Cor\label{cor:qutaff}
Assume that  $\shc$ is quasi-rigid.
Let $M$ and $N$ be simple objects of $\shc$ such that one of them is \afr.
\bnum
\item
Then $M\hconv N$ appears once in $M\tens N$ and 
$q^nM\hconv N$ does not appear in $M\tens N$ if $n\not=0$.
\item Let $\dX\in\Proc(\cor[z],\shc)$ such that 
$\dX/z\dX\simeq M\tens N$.
Then an epimorphism $\dX\epito \dL$
in $\Proc(\cor[z],\shc)$
for $\dL\in\Aff(\shc)$ such that $\dL/z\dL$ is simple
is unique \ro if it exists\rf. \ee
\encor

The category $\Modg(\cor[z],\Pro(\shc))$ is a monoidal category with a tensor product $\tensz$ defined by the following universal property: 
\eqn  
\hs{2ex}&&\ba{l}
\Hom_{\Modg(\cor[z], \Pro(\shc))}(\Ma\tensz\Na,\Laa)
\\ \hs{4ex}
\simeq \{f\in\Hom_{\Pro(\shc)}(\Ma\tens\Na,\Laa)\mid
(z\id_{\Laa})\circ f=f\circ(z\id_{\Ma}\tens\id_\Na)
=f\circ(\id_\Ma\tens z\id_\Na)
\}.
\ea\label{Eq: ten_A}
\eneqn
Note that $\tensz$ is right exact.
Moreover,  $(\Proc(\cor[z],\shc),\tensz\,)$ and $(\Aff(\shc),\tensz\,)$  become  monoidal categories with the unit object $\cor[z]$ (\cite[Lemma 5.1, Lemma 5.4]{refl}).
Note that there is an epimorphism $\Ma \tens \Na \epito \Ma\tensz \Na$ in $\Pro(\shc)$.

\Prop[{cf.\ \cite[Proposition 6.29 (ii) ]{refl}}]\label{prop:defactor}
Assume that $\shc$ is quasi-rigid.
Let $(\Ma,z)$, $(\Na,z)$ and $(\Laa,z)$ be an affinization with $\deg(z)=d$ of 
a real simple $M$, $N$  and $L$ in $\shc$, respectively.
Assume that 
\bna
\item $\de(M,N)>0$,
\item there exists an epimorphism 
$\Ma\tensz \Na\epito\Laa$  in $\Aff(\shc)$.\label{it:reg}
\ee
Then we have
$\De(\Ma,\Na)\in\cor[\z,z_{\Na}](\z-z_{\Na})$. 
Here $\z=z\vert_\Ma$ and $\zN=z\vert_{\Na}$.
\enprop

\Proof
Set $z=\z$, $w=\zN$ and $f(z,w)=\De(\Ma,\Na)$. 
Assume that  $f(z,w)$ is a homogeneous function of degree $r$ in $z,w$
(counting the degrees of $z,w$ as one),
i.e., $\de(M,N)=dr/2$.
Then we have
$\Rre_{\Ma,\Na}\cl\Ma\conv\Na\to\Na\conv\Ma$
and
$\Rre_{\Na,\Ma}\cl\Na\conv\Ma\to\Ma\conv\Na$
such that
$$\Rre_{\Na,\Ma}\circ\Rre_{\Ma,\Na}=f(z,w)\id_{\Ma\tens\Na}.$$
Hence, the composition
$$\Ma\tensz\Na\To[\Rre_{\Ma,\Na}\vert_{z=w}]\Na\tensz\Ma
\To[\Rre_{\Na,\Ma}\vert_{z=w}]\Ma\tensz \Na$$
is $f(1,1)z^r$.

Assume that $f(1,1)\neq 0$. 
Then
\eqn 
\phi : \Na \tensz \Ma \To[\bRre_{\Na,\Ma}\vert_{z=w}] \Ma\tensz \Na \To \Laa
\eneqn
is nonzero,
since $f(1,1)  z^r$ is a monomorphism. 
Hence 
\eqn
\phi = z^a \psi  \qt{for some } a \in \Z \  \text{and} \  \psi \in  \HOM_{\Proc(\cor[z],\shc)}(\Na\tensz \Ma, \Laa)  \ \text{such that} \ \psi\vert_{z=0} \neq 0.
\eneqn
It follows that 
\eqn
N\hconv M   \simeq  \Im(\psi \vert_{z=0})  \simeq L  \simeq M \hconv N
\eneqn
which contradicts $\de(M,N)>0$.
\QED


The following lemma is used in the proof of Lemma~\ref{lem:bijSW}.
\Lemma\label{lem:ijc}
Assume that $\shc$ is quasi-rigid.
Let $M$ and $N$ be \afr simple objects of $\shc$ and $\rmat{M,N}\cl M\tens N$ be
an R-matrix.
Then for any positive integer $n$,
the composition
$$f_n\cl M\tens N^{\otimes n}\To[\rmat{M,N}\tens N^{\otimes(n-1)}]
N\tens M\tens N^{\otimes(n-1)}\To\cdots\To  N^{\otimes n}\tens M$$
does not vanish,
and $f_n$ is equal to $\rmat{M,N^{\otimes n}}$
\enlemma
\Proof
We argue by the induction on $n$.
By the induction hypothesis, $f_{n-1}$ does not vanish and hence
the commutative diagram
$$\xymatrix@C=15ex{M\tens N^{\otimes n}\ar[d]^\bwr\ar[rr]^{f_n}
  &&N^{\otimes n}\tens M\ar[d]^\bwr\\
M\tens  N^{\otimes (n-1)}\tens N\ar[r]^{f_{n-1}}
& N^{\otimes (n-1)}\tens M\tens N\ar[r]^{ N^{\otimes (n-1)}\tens \rmat{M,N}}&
N^{\otimes (n-1)}\tens N\tens M
}
$$
implies the desired result.
Note that the composition in the bottom row is non-zero, since $M$ is simple.
\QED

\subsubsection{Rational center structures on duals in $\Aff(\shc)$}

For $M \in \shc$, we denote by $\D(M)$ (respectively, $\D^{-1}(M)$) the right dual (respectively,  the left dual) of $M$ if it exists.
The following proposition can be proved similarly to
\cite[Proposition 5.6]{refl}.

\Prop\label{prop:Dual}
Let $\Ma\in \Aff(\shc)$.
\bnum
\item
Assume that $\Ma/z^m\Ma$  admits a right dual  in $\shc$ for any $m\in \Z_{\ge 0}$.  
Then the right dual $\Da(\Ma)$ of $\Ma$ exists in the category $( \Aff(\shc),\tensz \,)$ and
\eqn
\Da(\Ma)/z^m\Da(\Ma)\simeq \D\bl(\cor[z]/z^m \cor[z])^*\tens_{\cor[z]} \Ma\br.
\eneqn
\item
Similarly if $\Ma/z^m\Ma$ admits a left dual  in $\shc$ for any $m\in \Z_{\ge 0}$,
then the left dual $\Da^{-1}(\Ma)$ exists in the category $( \Aff(\shc),\tensz \,)$ and
\eqn
\Da^{-1}(\Ma)/z^m\Da^{-1}(\Ma)\simeq \D^{-1}\bl(\cor[z]/z^m \cor[z])^*\tens_{\cor[z]} \Ma\br.
\eneqn
\ee
\enprop

Note that $(\cor[z]/z^{m} \cor[z])^*\simeq
z^{1-m}\cor[z]/z\cor[z]$.

\Prop[{\cite[Lemma 6.6]{refl}}] \label{prop:rational dual}
Assume that $\shc$ is rigid.
Let $(\Ma,\Rmat_{\Ma})$ be a rational center in $\shc$. 
Then $(\Da^{\pm1}(\Ma),\Rmat_{\Da^{\pm1}(\Ma)})$  is a rational center.

Here,$\Rmat_{\Da^{\pm1}(\Ma)}(X)\seteq \Da^{\pm1}(\Rmat_{\Ma}(\D^{\mp1}(X)))$.
\enprop

\subsubsection{Schur-Weyl dualities} \label{Sec: sub SW duality}
Let $R$ be the quiver Hecke algebra 
associated with  a Cartan matrix $\cartan$ and a set of polynomials  $\st{Q_{i,j}(u,v)}$ as in \eqref{eq:Q}.
Let $\st{(\hSW_i,z_i)}_{i\in I}$ be a family of affinizations in $\shc$
such that 
\eq  \label{eq: duality daum}
&&\left\{\parbox{70ex}{
\bna
\item $\deg(z_i)=(\al_i,\al_i)$,
\item $\SW_i\seteq\hSW_i/z_i\hSW_i$ is real simple in $\shc$ for any $i\in I$,
\item $\dim \HOM_\shc(\SW_i\tens \SW_j,\SW_j \tens \SW_i) =  \HOM_\shc(\SW_i\tens \SW_j, \SW_i \tens \SW_j)=1$ for any $i,j\in I$, \label{eq: one dim hom}
\item  $\De(\hSW_i,\hSW_j)\equiv Q_{i,j}(z_i,z_j)$ for $i\not=j$
(see Convention~\eqref{conv:4}).  \label{eq:DeQ}
\ee}\right.\label{eq:SW}
\eneq
We define the bilinear form  $\dg\cl\rtl\times\rtl\to\Z$ by
$$\dg(\al_i,\al_j)=
\bc\La(\SW_i,\SW_j)&\text{for $i,j\in I$ such that $i\not=j$.}\\
-(\al_i,\al_j)&\text{if $i=j\in I$.}\ec
$$
Then  $\dg$ is a grade associator, i.e., it satisfies \eqref{eq: dg_condition}.

We normalize $\{\Rre_{\hSW_i,\hSW_j}\}_{i,j\in I}$ such that
\eq
&&\left\{
\parbox{70ex}{
\bna
\item
$\Rre_{\hSW_j,\hSW_i}\circ \Rre_{\hSW_i,\hSW_j}
=Q_{i,j}(z_i,z_j)\id_{\hSW_i\tens \hSW_j}$ for $i\not=j$,
\item for any $i\in I$,
there exists $T_{i}\in\END_{\Pro(\shc)}(\hSW_i\tens \hSW_i)$
such that
$$\Rre_{\hSW_i,\hSW_i}-\id_{\hSW_i\tens \hSW_i}=(z_i\tens 1-1\tens z_i)\circ
T_{i}.$$
\ee}\right.
\label{eq:normR}
\eneq

We call $\bl\st{(\hSW_i,z_i)}_{i\in I},\st{\Rre_{\hSW_i,\hSW_j}}_{i,j\in I}\br$ a {\em duality datum}.

We denote by $\Modgc(R_\dg)$ the abelian category of finitely generated $R_\dg$-modules. It can be embedded into $\Pro(R_\dg)$ as a full subcategory. 

\Prop[{\cite[Proposition 7.6]{refl}}]  \label{prop:SWfunctor}
Let $\bl\st{(\hSW_i,z_i)}_{i\in I},\st{\Rre_{\hSW_i,\hSW_j}}_{i,j\in I}\br$ be a duality datum.
\bnum
\item
There exists a  right exact monoidal functor 
$$\hF\cl\Modgc(R_\dg)\to\Pro(\shc)$$
such that
\eq&&\hF(\tL(i)_{z_i})\simeq\hSW_i\qt{and}\quad  \hF(L(i))\simeq \SW_i,  \label{eq:SW1}\\
&&\hF(e(i,j)\vphi_1)=\Rre_{\hSW_i,\hSW_j}\in \HOM_{\cor[z_i,z_j]}(\hSW_i\tens \hSW_j,\hSW_j\tens \hSW_i). \label{eq:SW2}
\eneq
Here $\tL(i)_{z_i}$ is the affinization $\bl R(\al_i),z_i\br=(\cor[x_1],x_1)$
of $L(i)=R(\al_i)/R(\al_i)x_1$, 
$\vphi_1$ is given in \eqref{def:intertwiner}, and
$e(i,j)\vphi_1\cl R_\dg(\al_i+\al_j)e(i,j)
\to R_\dg(\al_i+\al_j)e(j,i)$ is the morphism by the right multiplication.

Moreover, such a right exact monoidal functor $\hF$ is unique up to an isomorphism.
\item
The functor $\hF$  can be restricted to a monoidal functor
$$\F\cl R_\dg\gmod \to\ \shc.$$
\item
Assume further that 
the functor $\hF$ is exact.  
\bna
\item If $(\Ma,\z)$ is an \afn of $M$ in $R_\dg\gmod$, then $(\hF(\Ma),\hF(\z))$ is an \afn of $\F(M)\in\shc$.

\item If $\shc$ is quasi-rigid, then  $\F$ sends a simple module in $R_\dg\gmod$  to a  simple objects in $\shc$ or $0$.
.\end{enumerate}

\ee
\enprop

\begin{remark}
If $\shc$ is quasi-rigid, then 
 the assumption  \eqref{eq: duality daum} \eqref{eq: one dim hom} holds. See \cite[Lemma 6.16, Proposition 6.19, and Corollary 6.21]{refl} for the consequence of quasi-rigidity assumption.
\end{remark}

\Prop [{\cite[Proposition 4.32]{Murata}}] \label{prop:F proj}
Let $\hF$ be the functor in  {\rm Proposition \ref{prop:SWfunctor}~(i)}.
The full subcategory $\Modg_{{\rm f.p.d}}(R_\dg)$ consisting of  $M$ in $\Modgc(R_\dg)$ with finite projective dimension is $\hF$-projective in the sense of \cite[Definition 13.3.4]{KS}.
\enprop

\section{Construction of $\Loc_\pm$} \label{Sec: Apm}

Let $\cartan = (\cartan_{i,j})_{i,j \in I}$ be the Cartan matrix of a symmetrizable Kac-Moody algebra $\g$ and
let $R$ be the quiver Hecke algebra associated with the quiver Hecke datum $\st{Q_{j,k}(u,v)}_{j,k\in I}$.

Let $i\in I$.
Define 
\begin{equation} \label{Eq: RiiR}
\begin{aligned}
\Ri&=\st{M\in R\gmod\mid \E_iM\simeq0}\qtq\\
\iR&=\st{M\in R\gmod\mid \Es_iM\simeq0}.
\end{aligned}
\end{equation}
One of the goals of this paper is to construct the reflection functor
$$\F_i\cl \Ri\isoto\iR.$$
As the first step of the goal,   in this section we will consider quiver Hecke algebras $R_\pm$ which are obtained by adding new vertices $i_\pm$ to the Dynkin diagram of $\g$,  and we will construct the monoidal categories $\Loc_\pm$ which are localizations of $R_\pm\gmod$ via a  left graded braider and a right graded braider, respectively. 

We add vertices $i_+$ and $i_-$ to $I$
to form $I_\pm \seteq I\sqcup\st{i_\pm}$, respectively. 
Let $\tfC_\pm=(\tfc_{i,j})_{i,j\in I_\pm}$,  be the Cartan matrix indexed by $I_\pm$,  where
\begin{align} \label{Eq: tC}
\tfc_{j,k}\seteq
\bc \sfc_{j,k}&\text{if $j,k\in I$,}\\
2&\text{if $j=k$,}\\
-1&\text{if $\st{j,k}=\st{i,i_\pm}$,}\\
0&\text{if $(j,k)\in\bl\st{i_\pm}\times(I\setminus\st{i})\br\cup
  \bl(I\setminus\st{i}\times\st{i_\pm}\br$.}
\ec
\end{align}
The inner product on the associated root lattice
$\rtl(R_\pm)\seteq\rtl\oplus\Z\al_{i_\pm}$
is given by
$$(\al_j,\al_k)=
\bc (\al_j,\al_k)&\text{if $j,k\in I$,}\\
-(\al_i,\al_i)/2&\text{if $\st{j,k}=\st{i,i_\pm}$,}\\
(\al_i,\al_i)&\text{if $j=k=i_\pm$,}\\
0&\text{otherwise.}
\ec$$

Let $R_\pm$ be the quiver Hecke algebra
associated with the quiver Hecke datum $\st{Q^\pm_{j,k}(u,v)}_{j,k\in I_\pm}$,
where
\begin{align} \label{Eq: Qpm}
Q^\pm_{j,k}(u,v)=
\bc Q_{j,k}(u,v)&\text{if $j,k\in I$,}\\
u-v&\text{if $j=i_\pm$, $k=i$.}\\
1&\text{if $j=i_\pm
 $ and $k\in I\setminus\st{i}$.}
\ec
\end{align}

We give a grading on $R_\pm$ as follows: 
\eq \label{eq:new grading}
&&\deg\bl e(\nu)\br=0,\quad\deg\bl x_ke(\nu)\br=2\sfd_{\nu_k},
\quad\\
&&\deg\bl e(\nu)\tau_k \br=\la_\pm(\al_{\nu_k},\al_{\nu_{k+1}})
\qtq[or]
  \deg\bl \tau_ke(\nu)\br=\la_\pm(\al_{\nu_{k+1}},\al_{\nu_{k}}),
\eneq
where $\la_\pm$ is a grade associator of $R_\pm$ defined by
  \eq \label{eq:degree modified}
 && \ba{ll}
  \la_\pm(\al_j,\al_k)=-(\al_j,\al_k)\hs{2ex}&\text{if $j,k\in I$,}\\
  \la_{\pm}(\al_{i_\pm},\al_{i_\pm})=-(\al_i,\al_i) = -(\al_{i_\pm}, \al_{i_\pm}),\\
  \la_+(\al_i,\al_{i_+})=
    \la_-(\al_{i_-},\al_i)=0,\\
    \la_+(\al_{i_+},\al_j)=\la_-(\al_j,\al_{i_-})=
    (\al_i,\al_j)&\text{if $j\in I$,}\\
    \la_+(\al_j,\al_{i_+})=\la_-(\al_{i_-},\al_j)=-(\al_i,\al_j)\hs{1ex}&\text{if $j\in I\setminus\st{i}$.}
  \ea
  \eneq
  Note that
  \eqn
  &&\la_\pm(\al,\beta)+\la_\pm(\beta,\al)=-2(\al,\beta)\qt{for any $\al,\beta\in\rtl(R_\pm)$, } 
\eneqn
so that the grading  \eqref{eq:new grading} on $R_\pm(\beta)$ is well-defined.
From now on, we shall denote by $R_\pm(\beta)$ the $\Z$-graded algebra $(R_\pm)_{\la_\pm}(\beta)$ with the grading \eqref{eq:new grading} for $\beta \in \rtl(R_\pm)$.

  By the definition, we have  in $R_\pm(\beta) \gmod$ 
  \eqn
\La(\ang{j},\ang{k})= \delta(j\neq k) \la_\pm(\al_j,\al_k)\qt{for any
  $j,k\in I_\pm$.}\eneqn

\begin{remark}
  We modified the grading so that
 the  braider $R_{C_\pm}$ considered in Proposition \ref{prop:Cpm} below has homogeneous degree $0$ 
  and
  the functor $\Psi_{\pm,\mp}$
defined  below is compatible with the grading.
\end{remark}

Since $ \la_\pm(\cdot,\cdot) = - (\cdot, \cdot)$ on $\rtl$, we have an embedding, which is compatible with the grading by \eqref{eq:degree modified}, 
  $$R\gmod\into R_\pm\gmod,$$
From now on, we regard $R\gmod$ as a full subcategory of $R_\pm\gmod$ via this embedding.

  We set
 \begin{align} \label{Eq: Cpm} 
  C_-\seteq\ang{i_-\;i}\in R_-\gmod\qtq
 C_+\seteq\ang{i\,i_+}\in R_+\gmod. 
 \end{align}
Then $(C_+,M)$ and $(M,C_-)$ are unmixed pairs for any simple $M \in R\gmod$, since
$\gW^*(C_+)=\st{0, \al_{i_+},\al_i+\al_{i_+}}$, and $\gW(C_-)=
\st{0,\al_{i_-},\al_i+\al_{i_-}}$.

Note that for $j,k \in I_\pm$ such that $(\al_j,\al_k)<0$,  we have 
\eqn 
&&\La(\ang{kj},\ang{j}) = \La(\ang{k},\ang{kj}) = \dg_\pm(\al_k,\al_j),  \\
&& \La(\ang{kj},\ang{k}) = \dg_\pm(\al_j,\al_k) -(\al_k,\al_k), \qt{and} \\
&&\ \La(\ang{j}, \ang{kj}) = \dg_\pm(\al_j,\al_k) - (\al_j,\al_j).
\eneqn
Hence by \eqref{eq:degree modified}, we have  
 \eq \label{eq:LaiC+}
&&\La(i,C_\pm)=\La(C_\pm,i)=
  \La(i_\pm,C_\pm)=\La(C_\pm,i_\pm)=0.
\eneq

  Note that we have an equivalence of graded monoidal categories
  $$(R_-\gmod)^\rev\isoto R_+\gmod$$
  by $\ang{j}\mapsto\ang{j}$ for $j\in I$ and
  $\ang{i_-}\mapsto \ang{i_+}$  (see Convention~\eqref{conv3}). 
  By this equivalence, $C_-$ is sent to $C_+$.
  Hence the statements on $R_\pm$ have counterparts for $R_\mp$,
  and we sometimes give the proofs only for $R_+$.
  
  \Prop \label{prop:Cpm} The objects $C_\pm$ have the following properties.
  \bnum
\item    We have $\La(C_+,\ang{j})=0$ for any $j\in I_+$ and 
  $\La(\ang{j},C_-)=0$ for any $j\in I_-$.
  In particular,  $\La(C_+,M)\le0$ for any simple $M\in R_+\gmod$
and    $\La(M,C_-)\le0$ for any simple $M\in R_-\gmod$.
\item
$C_+$ gives a non-degenerate  graded real  left braider in $R_+\gmod$, where
$\coRl_{C_+}(M)\cl C_+\conv M\to M\conv C_+$ 
has homogeneous degree $0$ for any $M\in R_+\gmod$.
\item
$C_-$ gives a non-degenerate graded real right braider in $R_-\gmod$, where
$\coRr_{C_-}(M)\cl M\conv C_-\to C_-\conv M$
has homogeneous degree $0$ for any $M\in R_-\gmod$.
\ee
\enprop

\Proof 
(i) 
By \eqref{eq:LaiC+}, we may assume that $j\neq i, i_+$.
Then we have
\eqn
\La(C_+,\ang{j}  ) = \La(\ang{i}, \ang{j} )+     \La(\ang{i_+}, \ang{j} )
=  \dg_+(\al_i,\al_j)+ \dg_+(\al_{i_+},\al_j)  =-(\al_i,\al_j)+(\al_i,\al_j) =0. 
\eneqn
The second assertion follows from the fact that for every simple $M\in R_+ \gmod$, there exists $m\in \Z_{\ge 0}$,  $N\in R_+\gmod$ and $j\in I_+$ such that $M=\ang{j^m} \hconv N$.  

(ii) and (iii) follows from Proposition \ref{prop:nondeg braider} and (i).
The reality follows from
$\La(C_+,C_+)=\La(C_+,\ang{i})+\La(C_+,\ang{i_+})$.
\QED

In the sequel, we normalize the braiding $\Rmat_{C_\pm}$ so that it satisfies
\eq
\Rmat_{C_\pm}(C_\pm)=\id_{C_\pm\circ C_\pm}.
\eneq

Let us denote by $\Loc_\pm$
the localization
$R_\pm\gmod[C_\pm^{\circ -1}]$ of $R_\pm\gmod$ by the left and right braider $C_\pm$, respectively.

Let $\Qt_\pm\cl R_\pm\gmod\to \Loc_\pm$ be the localization functor. 
It is an exact monoidal functor by Theorem \ref{Thm: graded localization}.

\Lemma[{\cite[Proposition 4.8]{loc1}}] \label{lem: Q simple}
For any simple $M\in R_+\gmod$,
the object  $\Qt_+(M)$ is simple or zero. 
 The object $\Qt_+(M)$ is simple
if and only  if  $\La(C_+,M)=0$.
\enlemma

\Lemma\label{lem:jC+}
For any $j\in I$ such that $(\al_i,\al_j)<0$, we have
 $$\La(C_+,\ang{j\,i\,i_+})=-(\al_i,\al_i)<0.$$ 
Here, $\ang{jii_+}$ is the one-dimensional
simple head of $\ang{j} \conv  \ang{ii_+}$. 
In particular, we have
$\Qt_+(\ang{j\,i\,i_+})\simeq0$.
\enlemma
\Proof
 Since $(C_++, \ang{j\,i\,i_+})$ is an unmixed pair, we have (see \eqref{eq:unmixed})
$$\La(C_+,\ang{j\,i\,i_+})=\dg_+\bl(\wt(C_+),\wt(\ang{j\,i\,i_+}\br
=\dg_+(\al_i+\al_{i_+},\al_j+\al_i+\al_{i_+})=-(\al_i,\al_i).$$ 
\QED

  We define
  \eqn
  \tC_\pm&&\seteq\Qt_\pm(C_\pm)\in\Loc_\pm.
    \eneqn

  \Prop\label{prop:RLocf}
  The restriction
  $$\Qt_\pm\vert_{R\gmod}\cl R\gmod\to\Loc_\pm$$
  is fully faithful.
\enprop
\Proof
In the proof we ignore the grading.
We prove the assertion only for $\Qt_+$.
Let $M,N\in R\gmod$.
Then, we have
$$\Hom_{\Loc_+}(\Qt_+(M),\Qt_+(N)\br
=
\indlim_{n\in\Z_{\ge0}}\Hom_{R_+}\bl C_+^{\circ n}\conv M,N\conv C_+^{\circ n}\br.$$
On the other hand,
the pair $\bl C_+^{\circ n},N\br$ is unmixed so that  
$\Res(N\conv C_+^{\circ n})\simeq C_+^{\circ n}\tens N$.
Hence, we have
\eqn
\Hom_{R_+}\bl C_+^{\circ n}\conv M,N\conv C_+^{\circ n}\br&&\simeq
\Hom_{R_+ \tens R_+ }\bl C_+^{\circ n}\tens M, \Res(N\conv C_+^{\circ n})\br\\
&&\simeq\Hom_{ R_+\tens R_+}\bl C_+^{\circ n}\tens M, C_+^{\circ n}\tens N\br
\simeq\Hom_{ R_+ }(M,N).
\eneqn
Thus, we obtain
$\Hom_{\Locp}\bl\Qt_+(M),\Qt_+(N)\br\simeq\Hom_{R}(M,N)$.
\QED

\Rem
In order to define the monoidal category $\Loc_+$, we have used the non-degenerate left braider associated with $C_+$.
If we use the non-degenerate right braider associated with $C_+$ instead, then Proposition~\ref{prop:RLocf} fails.

Indeed, if $\Qt\cl R_+\gmod\to \Loc_+$
is the localization functor by the right braider, then
$\Qt(\ang{12})\simeq0$ when $\g=A_2$  and $C_+=\ang{11_+}$.
\enrem

  \Prop \label{prop:stab subq}
The category $R\gmod$ is stable by taking subquotients
  as a full subcategory of $\Loc_\pm$.
  \enprop
  \Proof
  We shall prove only for $\Loc_+$.
  We shall show that $R\gmod$ is stable by taking quotients.
  Let $M\in R\gmod$ and $f\cl\Qt_+(M)\epito X$ is an epimorphism
  in $\Loc_+$. Then there exist $N\in R_+\gmod$ and a morphism
  $h\cl C_+^{\circ n} \conv M\to N$ 
  such that $\Qt(h)$ is isomorphic to $C_+^{\circ n} \conv f$. 
In particular,  $\Qtp(N) \simeq \tC^{\circ n}_+ \conv X$. 
  Replacing $N$ with the image of $h$, we may assume that
  $h$ is an epimorphism from the beginning.
  Set $\al=-\wt(C_+^{\circ n})$ and $\beta=-\wt(M)$.
  Since $(C_+^{\circ n}, M)$ is unmixed, we obtain an epimorphism
  \eqn
  C_+^{\circ n}\tens M\simeq\Res_{\al,\beta}(C_+^{\circ n} \conv M)\epito
  \Res_{\al,\beta}(N).
  \eneqn
  Hence there exists $L\in R\gmod$ such that
  $\Res_{\al,\beta}(N)\simeq C_+^{\circ n}\tens L$,
  and we obtain an epimorphism
  $M\epito L$ in $R\gmod$. 
  By the functoriality,
  we have
  $$C_+^{\circ n} \conv M\epito C_+^{\circ n} \conv L\to N
  \to L\conv C_+^{\circ n}.$$ 
Since the composition $C_+^{\circ n} \conv L\to N
\to L\conv C_+^{\circ n}$ is an isomorphism
after taking $\Qt_+$, we obtain
$\Qt_+(N)\simeq \Qt_+(C_+^{\circ n} \conv L)$, which implies that
$X\simeq\Qt_+(L)$.
\QED

\Def Let $\Loco$ be the full subcategory of
$R_\pm\gmod$ which consists of modules $M\in R_\pm\gmod$ which satisfy the following equivalent conditions:
\bna
\item
  any simple subquotient $S$ of $M$ commutes with $C_\pm$ and satisfies
  $\La(C_+,S)=0$ (resp.\ $\La(S,C_-)=0$).
\item
  $ \coRl_{C_+}(M)\cl C_+\conv M\to M\conv C_+$
  \ro resp.\ $\coRr_{C_-}(M)\cl M\conv C_-\to  C_-\conv M$\rf is an isomorphism.
  \ee
  \edf
  Note that that $C_\pm$ is a central object of  $\Loco$.

  \Lemma[{\cite[Proposition 4.9]{loc1}}] \label{lem:Loc0} 
  The subcategory $\\Loco$ is closed by taking
 subquotients, extensions and convolution products. 
Moreover, the canonical functor
  $\Loco[(C_\pm)^{-1}]\to\Loc_\pm$ is an equivalence of monoidal categories.
  \enlemma

  \Prop[{cf. \cite[\S\;8.2]{refl}}] \label{prop:Loc0}
  The functor $\Loco\to\Loc_\pm$ is fully faithful.
  \enprop

  \Proof
  By using \cite[Proposition 8.1]{refl}, we can prove the assertion similarly to the proof of
  \cite[Theorem 8.3 and Corollary 8.4]{refl}.
  \QED

\Th\label{th:qusirigid}
The monoidal category $\Loc_\pm$ is quasi-rigid.
\enth
\Proof
This follows from
Lemma~\ref{lem:Loc0}, Proposition~\ref{prop:Loc0} and the fact that $\Loco$ is quasi-rigid. 
\QED

\Lemma 
If $M\in R_+\gmod$ is  \afr simple,
then so is
$\Qtp(M)$ \ro in the sense of  {\rm Definition \ref{def:affinization}}\rf, provided $\Qtp(M)\neq 0$.
\enlemma
\Proof
Let $\Ma$ be an affinization of $M$. 
Since $\Qtp$ is an exact monoidal functor from $\Pro(R_+\gmod) $ to $\Pro(\Loc_+)$,
  $\Qtp(\Ma) \in \Pro(\Loc_+)$ is an affine object of $\Qtp(M)$ with 
$z_{\Qtp(\Ma) } = \Qtp(z_\Ma)$. 
Now $\Rmat_{\Qtp(\Ma)}(\Qtp(X))= \Qtp(\Runi_{\Ma}(X))$ 
for $X \in R_+\gmod$ gives a rational center structure on $\Qtp(\Ma)$ by the same proof in \cite[Lemma 8.8]{refl}.
\QED

\Lemma \label{lem:aff dual}\hfill
\bnum
\item
  The object $\qi$ has a right dual $\Qtp(\ang{i_+})\conv(\tC_+)^{\circ-1}$
  in $\Loc_+$,
  and $\Qtp(\ang{i_+})$ has a left dual $(\tC_+)^{\circ-1}\conv\Qi$
  in $\Locp$.
\item
  $\Qiz/z^n\Qiz$ has a right dual
  $\Qtp\bl\E_i\bl(C_+)^{\circ n}\br\conv (\tC_+)^{\circ-n}$ for any $n\in\Z_{\ge0}$.
  \item
   there exists a right dual
  $$\Diz\simeq\proolim[n] z^{n-1}\cor[z]\tens_{\cor[z]}\Qtp\bl\E_i\bl(C_+)^{\circ n}\br\br\conv (\tC_+)^{\circ-n}$$
  of $\qiz$ in $\Aff(\Loc_+)$,
  and the left dual
  $$\Daf^{-1}\bl\Qtp(\ang{i_+}_{z})\br\simeq\proolim[n]
  z^{n-1}\cor[z]\tens_{\cor[z]}\Bigl((\tC_+)^{\circ-n}\conv\Qtp\bl\Es_{i_+}\bl(C_+)^{\circ n}\br\br\Bigr)$$
  of $\Qtp(\ang{i_+}_{z})$ in $\Aff(\Loc_+)$.
  \ee
  \enlemma

  \Proof
  It follows from Proposition ~\ref{prop: right dual} and Proposition~\ref{prop:Dual}.
  \QED

  The proof of the following proposition will be given in  Appendix \ref{sec:appendixA}.
  \Prop\label{prop:affD}
The affine objects $\Diz$ and $\Dipz$ have rational center structure
$(\Diz, \Rmat_{\Diz})$ and $(\Dipz,\Rmat_{\sDipz})$, respectively \ro see
{\rm Definition~\ref{def:rational center}}\rf.
  \enprop
{\em In the sequel, we normalize $\Rmat_{\sDiz}$ and $\Rmat_{\sDipz}$ as in the following lemma.}

\Lemma\label{lem:normR}
We can normalize $\Rmat_{\sDiz}$ and $\Rmat_{\sDipz}$ so that
\bna
\item\label{itemha}
  $\Rmat_{\sDiz}(X)$ \ro resp.\ $\Rmat_{\sDipz}(X)$ \rf
  is a morphism in $\Proc(\cor[z],\Loc_+)$ and homogeneous of degree
  $-\la_+(\al_{i_+},\wt(X))$ \ro resp.\ $-\la_+(\al_{i},\wt(X))$\rf
\item
  $\Rmat_{\sDiz}(X)$ is a renormalized R-matrix if
  $X=\Qi[j]$ with $j\in I$ or $X=\tC_+^{\circ\pm1}$.
\ee
\enlemma
\Proof Since the proof of $\Rmat_{\sDipz}$ is similar, 
we only give the proof for $\Rmat_{\sDiz}$.

 For $X\in \Loc_+$ with $\wt(X )=\sum_{j\in I_+} c_j \al_j$, 
taking suitable $m_j$'s and by replacing $\Rmat_{\sDiz}(X)$ with
$z^{\sum_{j\in I} c_j m_j}\Rmat_{\sDiz}(X)$,
we may assume from the beginning that
\eq
&&\text{$\Rmat_{\sDiz}(\Qi[j])$ is a renormalized R-matrix for any $j\in  I_+$.}
\eneq

Then we may assume that $\Rmat_{\sDiz}(X)$ is a morphism in $\Proc(\cor[z],\Loc_+)$ 
homogeneous of degree $-\la_+(\al_{i_+},\wt(X))$ for any $X\in\Loc_+$.
Note that $\La(\Di,\Qtp(\ang{j})) =   \La(\Qtp(\ang{i_+})\conv \tC_+^{-1},\Qtp(\ang{j})) =  \La(\Qtp(\ang{i_+}),\Qtp(\ang{j}))  = \La(\ang{i_+},\ang{j})
=\la_+(\al_{i_+},\al_j)$.

Since $\La(\Di, \tC_+^{\circ\pm1})=0=
-\la_+(\al_{i_+},\wt(\tC_+^{\circ\pm1}))$, we conclude that
$\Rmat_{\sDiz}(\tC_+^{\circ\pm1}))$ is a renormalized R-matrix.
Hence it is an isomorphism in $\Proc(\cor[z],\Loc_+)$.
\QED

\Lemma \label{lem:LaQ}
Let $M,N$ be simple modules in $R_+\gmod$.
Assume that one of them is \afr.
If\/ $\Qt_+(M\hconv N)$ is simple,
then $\Qt_+(M\hconv N)\simeq\Qtp(M)\hconv\Qtp(N)$ and
$$\La\bl\Qt_+(M),\Qt_+(N)\br=\La(M,N).$$
\enlemma
\Proof
The composition
$\rmat{M,N}\cl M\conv N \epito M\hconv N\monoto N\conv M$
induces
$ \Qtp(\rmat{M,N})\cl \Qtp(M)\conv\Qtp(N) \epito \Qtp(M\hconv N)\monoto \Qtp(N)\conv \Qtp(M)$, which does not vanish by the assumption. 
Note that $\HOM_{\Loc_+}(\Qtp(M)\conv \Qtp(N), \Qtp(N)\conv \Qtp(M)) = \cor \Qtp(\rmat{M,N})$, since $\Loc_+$ is quasi-rigid.
\QED

\Lemma\label{lem:LaQ+2}
Let $M,N$ be simple modules in $R_+\gmod$.
Assume that one of them is \afr. 
If $M\in\Loco[+]$  and $\Qtp(N)$ is simple, then
$\Qtp(M\hconv N)$ is simple and
$$\La\bl\Qt_+(M),\Qt_+(N)\br=\La(M,N).$$
\enlemma
\Proof
Since $C_+$ and $M$ commute, we have
$\La(C_+ , M\hconv N)=\La(C_+ , M)+\La(C_+ ,N)=0+0=0$,
and hence
$\Qtp(M\hconv N)$ is simple by Lemma \ref{lem: Q simple}.
Hence the assertion follows  by Lemma  \ref{lem:LaQ}. \QED

\Lemma \label{lem:C+Kj}
For $j\in I\setminus\st{i}$,
we have
$$\La(C_+,\ang{i^{-\sfc_{i,j}}\,j})=
\La(\ang{i^{-\sfc_{i,j}}\,j},C_+)=0.$$
In particular, we have
$$\de(C_+,\ang{i^{-\sfc_{i,j}}\,j})=0.$$
  \enlemma
  \Proof
  Note that $\ang{i^{-\sfc_{i,j}}\,j}$ and $C_+$ commute with $\ang{i}$ (see, for example, \cite{KKOP18}).
  We have
  \eqn
  \La(C_+,\ang{i^{-\sfc_{i,j}}\,j})
    && =  \La(C_+,\ang{i^{-\sfc_{i,j}}})+\La(C_+,\ang{j})=0
      \eneqn
      and
    \eqn
  \La(\ang{i^{-\sfc_{i,j}}\,j},C_+)
  && = \La(\ang{i^{-\sfc_{i,j}}\,j},\ang{i}) +
  \La(\ang{i^{-\sfc_{i,j}}\,j},\ang{i_+})\\
&&=-  \La(\ang{i},\ang{i^{-\sfc_{i,j}}\,j}) +
\La(\ang{i^{-\sfc_{i,j}}},\ang{i_+})
+\La(\ang{j},\ang{i_+})\\
&&=-\La(\ang{i},\ang{j})+\La(\ang{j},\ang{i_+})
=(\al_i,\al_j)-(\al_i,\al_j)=0.
      \eneqn
      \QED

      \Cor

The $R_+$-modules
      $$\ang{i},\, \ang{i_+},\, \ang{i^{-\sfc_{i,j}}\,j}
        \qt{$(j\in I\setminus\st{i})$}$$ 
        belong to $\Loco[+]$.
      \encor

\section{Schur-Weyl duality functors} \label{Sec: SW duality}
In this section we will define Schur-Weyl duality functors
$\Psi_{\pm,\mp}\cl R_\mp \gmod\to \Loc_\pm$.

Recall that 
for $w,v$ in the Weyl group  with $w \ge v$ in the Bruhat order and $\la \in \wtl_+$, there exists a \afr simple module $\dM(w\la, v\la)$, called the \emph{determinantial module} (see \cite[Section 3.3]{loc1}).  
Let $j\in I\setminus\st{i}$. Then the determinantial module  $\dM(s_is_j\La_j,\La_j)$  is isomorphic to $\ang{i^cj}\seteq  \ang{i^c}\hconv \ang{j}$,
where $c=-\ang{h_i,\al_j}\in\Z_{\ge0}$.

We omit the proof of the lemma below since it is similar to \cite[Lemma 9.11]{refl}. 
\Lemma\label{lem:cusp} 
Let $j\in I\setminus\st{i}$.
Set $c=-\ang{h_i,\al_j}\in\Z_{\ge0}$.
Then we have
\bnum
\item
  $\dM(s_is_j\La_j,\La_j) \simeq \ang{i^cj}$.
  Setting $\beta=c \al_i+\al_j$ and $C=R(\beta)$, we have
  \eqn\ang{i^cj}&&\simeq
  \dfrac{C}
  {\sum\limits_{k=1}^{c-1} C\tau_k+
    \sum\limits_{\nu\in I^{\beta},\nu_{c+1}=i} Ce(\nu)C+
    C x_{c+1}e(i^cj)C}\\
  &&\simeq
  \dfrac{Ce(i^cj)}
  {\sum\limits_{k=1}^{c} C\tau_ke(i^cj)+C x_{c+1}e(i^cj)}.
  \eneqn
Hence 
$\ang{i^cj}\simeq \ang{i^c}\tens\ang{j}$
  as an $R(c\al_i)\tens R(\al_j)$-module, 
\item its affinization $\ang{i^cj}_{z_j}\seteq\tdM(s_is_j\La_j,\La_j)_{z_j}$
  is given by
\eqn
  &&\dfrac{C}{\sum_{k=1}^{c-1} C\tau_k+
    \sum\limits_{\nu\in I^{\beta},\nu_{c+1}=i} Ce(\nu)C+
    C(z_j- x_{c+1})e(i^cj)C}\\
 &&\hs{10ex}\simeq
  \dfrac{Ce(i^cj)}
  {\sum\limits_{k=1}^{c} C\tau_ke(i^cj)+C (z_j-x_{c+1})e(i^cj)},
  \eneqn
  where $\beta=c\al_i+\al_j$ and $C=\cor[z_j]\otimes R(\beta)$. 
  Note that $\tau_c^2e(i^cj)=Q_{i,j}(x_c,x_{c+1})$. 
  \item  $\ang{i^cj}_{z_j}$ is isomorphic to
  $$\ang{i^c}_{z_j}\tens_{\cor[z_j]} R(\al_j)$$
  as an $\cor[z_j]\tens R(c\al_i)\tens R(\al_j)$-module.
  Here $\ang{i^c}_{z_j}=\Es_j\bl \ang{i^cj}_{z_j}\br$, 
\item $\ang{i^c}_{z_j}\simeq\dfrac{\cor[z_j]\tens R(c\al_i)}
  {\sum\limits_{k=1}^{c-1}\bl\cor[z_j]\tens R(c\al_i)\br\tau_k+
    \bl\cor[z_j]\tens R(c\al_i)\br Q_{i,j}(x_c,z_j)}$.
\item
  $\ang{i^c\,j}_{z_j}$ is isomorphic to the image of
  the morphism in $\Proc(\cor[z_j],R\gmod)$
  $$\ang{i^c}_{z_j}\conv[z_j]\ang{j}_{z_j}\To\ang{j}_{z_j}\conv[z_j]\
\ang{i^c}_{z_j}.$$
  \ee
\enlemma
\vskip 1em 

Let us define the Schur-Weyl duality datum $\st{K^+_j}_{j\in I_-}$ in
$\Loc_+$ as follows:
\begin{align} \label{Eq: Kj+}
K^+_j=\bc
q_i^{-1}  \D^{-1}\bl\Qt_+(\ang{i_+})\br\simeq  q_i^{-1} \tC_+^{\circ-1}\circ\Qt_+(\ang{i})&\text{if $j=i_-$,}\\
\Qt_+\bl\dM(s_is_j\La_j,\La_j)\br\simeq\Qt_+(\ang{i^{-\sfc_{i,j}}\,j})&\text{if $j\in I\setminus\st{i}$,}\\
  q_i   \D \bl\Qt_+(\ang{i})\br\simeq   q_i  \Qt_+(\ang{i_+})\conv \tC_+^{\circ-1}&\text{if $j=i$.}
      \ec
      \end{align}
(recall that $q_i=q^{(\al_i,\al_i)/2)}$)     and its affinization $\st{\tK^+_j}_{j\in I_-}$
\begin{align} \label{Eq: aKj+}
\tK^+_j=\bc
q_i^{-1}  \rDa^{-1}\bl\Qtp(\ang{i_+}_{z_j})\br
&\text{if $j=i_-$,}\\
\Qt_+(\ang{i^{-\sfc_{i,j}}\,j}_{z_j})&\text{if $j\in I\setminus\st{i}$,}\\
          q_i  \rDa\bl\Qtp(\ang{i}_{z_j})\br&\text{if $j=i$.}
      \ec
      \end{align}
      Here $z_j$ is a homogeneous indeterminate with degree $(\al_j,\al_j)$.
Indeed, $\tK^+_j$ is an affinization of $K^+_j$ by Lemma \ref{lem:Loc0},  Lemma \ref{lem:aff dual} and Proposition \ref{prop:rational dual}. 
      Similarly we define
      the Schur-Weyl duality datum $\st{K^-_j}_{j\in I_+}$ in
$\Loc_-$ as follows:
\begin{align} \label{Eq: Kj-}
K^-_j=\bc
 q_i^{-1} \D\bl\Qt_-(\ang{i_-})\br\simeq q_i^{-1}\Qt_-(\ang{i})\circ \tC_-^{\circ-1}
&\text{if $j=i_+$,}\\
\Qt_-(\ang{j\,i^{-\sfc_{i,j}}})&\text{if $j\in I\setminus\st{i}$,}\\
  q_i      \D^{-1}(\Qt_-\ang{i})\simeq  q_i \tC_-^{\circ-1}\conv\Qt_-(\ang{i_-})&\text{if $j=i$.}
      \ec
      \end{align}
      and its affinization $\st{\tK^-_j}_{j\in I_+}$
\begin{align} \label{Eq: aKj-}
\tK^-_j=\bc
 q_i^{-1}   \rDa\bl\Qt_-(\ang{i_-}_{z_{j}})\br
&\text{if $j=i_+$,}\\
\Qt_-(\ang{j\,i^{-\sfc_{i,j}}}_{z_j})&\text{if $j\in I\setminus\st{i}$,}\\
 q_i   \rDa^{-1}\bl\Qt_-(\ang{i}_{z_j})\br&\text{if $j=i$.}
      \ec
      \end{align}

      \Prop\label{prop:LaSW}
For any $j,k\in I_\mp$ such that $j\not=k$, we have
$$\La(K^\pm_j,K^\pm_k)=\la_{\mp}(\al_j,\al_k).$$
\enprop
\Proof
We prove it only for $K^+$.

\snoi
{\bf Case} $j=i_-$ and $k\in I\setminus\st{i}$.\quad
Let us show
$$\La\bl\D^{-1}\Qt_+(\ang{i_+}),\Qt_+(\ang{i^{-c_{i,k}}\,k})\br
=\la_-(\al_{i_-},\al_k).$$
We have
\eqn
\La\bl\D^{-1}\Qt_+(\ang{i_+}),\Qt_+(\ang{i^{-c_{i,k}}\,k})\br
&&=\La\bl\Qt_+(\ang{i^{-c_{i,k}}\,k}),\Qt_+(\ang{i_+})\br\\
&&=\La\bl\ang{i^{-c_{i,k}}\,k},\ang{i_+}\br
=\La\bl\ang{i^{-c_{i,k}}},\ang{i_+}\br
+\La\bl\ang{k},\ang{i_+}\br\\
&&=
-(\al_i,\al_k)=\la_-(\al_{i_-},\al_k),
\eneqn
where the first equality comes from \cite[Lemma 6.30]{refl}, and the second comes form Lemma \ref{lem:LaQ+2}.  

\snoi
{\bf Case} $k=i_-$ and $j\in I\setminus\st{i}$.\quad
Let us show
$$\La\bl\Qt_+(\ang{i^{-c_{i,j}}\,j}),\D^{-1}\Qt_+(\ang{i_+}))\br
=\la_-(\al_j,\al_{i_-}).$$
Since $\D^{-1}\Qt_+(\ang{i_+})\simeq \tC_+^{\circ-1}\conv\Qt_+(\ang{i})$, we have
\eqn
\La\bl\Qt_+(\ang{i^{-c_{i,j}}\,j}),\D^{-1}\Qt_+(\ang{i_+}))\br
&&=-\La\bl\Qt_+(\ang{i^{-c_{i,j}}\,j}),\tC_+\br
+\La\bl\Qt_+(\ang{i^{-c_{i,j}}\,j}),\Qt_+(\ang{i})\br\\
&&=\La\bl\tC_+,\Qt_+(\ang{i^{-c_{i,j}}\,j})\br+
\La\bl\Qt_+(\ang{i^{-c_{i,j}}\,j}),\Qt_+(\ang{i})\br\\
&&=\La\bl C_+,\ang{i^{-c_{i,j}}\,j}\br+
\La\bl \ang{i^{-c_{i,j}}\,j},\ang{i}\br\\
&&=-\La\bl\ang{i},\ang{j})=(\al_i,\al_j)=\la_-(\al_j,\al_{i_-}),
\eneqn
where the second equality and the third last equality follow from Lemma \ref{lem:C+Kj}.

\snoi
{\bf Case} $j=i$ and $k\in I\setminus\st{i}$.\quad
Let us show
$$\La\bl\D\Qt_+(\ang{i}),\Qt_+(\ang{i^{-c_{i,k}}\,k})\br
=\La(\ang{i},\ang{k}).$$
We have
\eqn
\La\bl\D\Qt_+(\ang{i}),\Qt_+(\ang{i^{-c_{i,k}}\,k})\br
&&=\La\bl \Qt_+(\ang{i_+})\conv\tC_+^{\circ-1},\Qt_+(\ang{i^{-c_{i,k}}\,k})\br\\
&&=\La\bl \Qt_+(\ang{i_+}),\Qt_+(\ang{i^{-c_{i,k}}\,k})\br
-\La\bl\tC_+,\Qt_+(\ang{i^{-c_{i,k}}\,k})\br\\
&&=\La\bl \ang{i_+},\ang{i^{-c_{i,k}}\,k}\br\\
&&=-c_{i,k}\La(\ang{i_+},\ang{i})+\La(\ang{i_+},\ang{k})\\
&&=-c_{i,k}\cdot 2\sfd_i+(\al_i,\al_k)=-(\al_i,\al_k).
\eneqn

\snoi
{\bf Case} $j\in I\setminus\st{i}$ and $k=i$.\quad
Let us show
$$\La\bl\Qt_+(\ang{i^{-c_{i,j}}\,j}),\D\Qt_+(\ang{i})\br
=\La(\ang{j},\ang{i}).$$

By Lemma~\ref{lem:LaXY}, we have
\eqn
\La\bl\Qt_+(\ang{i^{-c_{i,j}}\,j}),\D\Qt_+(\ang{i})\br
&&=\La\bl \Qtp(\ang{i}),\,\Qt_+(\ang{i^{-c_{i,j}}\,j})\br\\
&&=\La\bl \ang{i},\ang{i^{-c_{i,j}}\,j}\br
=\La(\ang{i},\ang{j}).
\eneqn

\snoi
{\bf Case} $j,k\in I\setminus\st{i}$ and $j\not=k$.\quad
Set $c=-\sfc_{i,j}$ and $c'=-\sfc_{i,k}$.
Then, since $\ang{i^c\,j}$ commutes with $i$, we have
\eqn
\La\bl\Qtp(\ang{i^c\,j}),\Qtp(\ang{i^{c'}\,k})\br&&=
\La(\ang{i^c\,j},\ang{i^{c'}\,k})\\
 && =\La(\ang{i^c\,j},\ang{i^{c'}})+\La(\ang{i^c\,j},\ang{k})\\
 && =-\La(\ang{i^{c'}},\ang{i^c\,j})+\La(\ang{i^c},\ang{k})+\La(\ang{j},\ang{k})
\\
&& =-c'\La(\ang{i},\ang{j})+c\La(\ang{i},\ang{k})+ \La(\ang{j},\ang{k}) \\
&&= c'(\al_i,\al_j)-c(\al_i,\al_k)+ \La(\ang{j},\ang{k})     =    \La(\ang{j},\ang{k}).
    \eneqn

    \snoi
    {\bf Case} $j=i$ and $k=i_-$.\quad
    Let us show $$\La\bl\D\Qt_+(\ang{i}),\D^{-1}\Qt_+(\ang{i_+})\br=\La(\ang{i},\ang{i_-}).$$
    We have
    \eqn
    \La\bl\D\Qt_+(\ang{i}),\D^{-1}\Qt_+(\ang{i_+})\br
     &&=\La\bl\Qt_+(\ang{i_+})\conv (\tC_+)^{\circ-1},
    \,(\tC_+)^{\circ-1}\conv\Qt_+(\ang{i}) \br\\
    &&  =\La\bl(\Qt_+(\ang{i_+}),\Qt_+(\ang{i})\br=\La(\ang{i_+},\ang{i})
    =\La(\ang{i},\ang{i_-}).
    \eneqn
\snoi
        {\bf Case} $j=i_-$ and $k=i$.\quad
        Let us show $$\La\bl\D^{-1}\Qt_+(\ang{i_+}),\D\Qt_+(\ang{i})\br
        =\La(\ang{i_-},\ang{i}).$$
    We have
    \eqn
\La\bl\D^{-1}\Qt_+(\ang{i_+}),\D\Qt_+(\ang{i})\br
    &&  =\La\bl\Qt_+(\ang{i}),\Qt_+(\ang{i_+}\br=\La(\ang{i},\ang{i_+})
    =\La(\ang{i_-},\ang{i}).
    \eneqn
    
    \QED

Recall that   $A \equiv B$ means that $A =c B$ for some $c\in \cor^\times$.
    \Th \label{thm:DeSW}
    For $j,k\in I_\mp$ such that $j\not=k$, we have
$$\De(\tK^\pm_j,\tK^\pm_k)  \equiv Q^{\mp}_{j,k}(z_j,z_k).$$
\enth
\Proof

We prove only the case $\tK^+$.

\snoi
{\bf Case} $j,k\in I\setminus\st{i}$.

We set $\tK=\tdM(s_is_j\La_j,\La_j)_{z_j}\in \Pro(R\gmod)$
and $\tK'=\tdM(s_is_k\La_k,\La_k)_{z_k}\in \Pro(R\gmod)$
Set $c=-\sfc_{i,j}$ and $c'=-\sfc_{i,k}$.
Let $R\cl \tK\conv \tK'\to \tK'\conv\tK$
and $R'\cl \tK'\conv \tK\to \tK\conv\tK'$
be the renormalized R-matrices (in $\Pro(R\gmod)$).
Note that
$\E_i^{(c)}\tK\simeq\ang{j}_{z_j}$ and
$\E_i^{(c')}\tK'\simeq\ang{k}_{z_k}$.

Applying $\E_i^{(c+c')}$ to $R$ and $R'$, we obtain
  $$\ang{j}_{z_j}\conv\ang{k}_{z_k}\To[\E_i^{(c+c')}(R)]\ang{k}_{z_k}\conv\ang{j}_{z_j}\To[\E_i^{(c+c')}(R')]
\ang{j}_{z_j}\conv\ang{k}_{z_k}.$$
Then, the composition coincides with
$\E_i^{(c+c')}\bl \De(\tK^+_j,\tK^+_k)\br = \De(\tK^+_j,\tK^+_k)$.
Remark that 
  $\E_i^{(c+c')}(R)$ and $\E_i^{(c+c')}(R')$ are renormalized R-matrices
  since they have the same degree by Proposition~\ref{prop:LaSW}.
  Hence the composition $\E_i^{(c+c')}(R')\cdot\E_i^{(c+c')}(R)$
  is equal to
  $\De(\ang{j}_{z_j},\ang{k}_{z_k})=Q_{j,k}(z_j,z_k)$.

    \mnoi
  {\bf Case}\ $j=i_-$ and $k\in I\setminus\st{i}$.\quad
  In this case, 
$$\de(K^+_j,K^+_k)= \dfrac{1}{2} \bl \dg_-(\al_{i_-}, \al_k) +\dg_-(\al_k,\al_{i_-} ) \br = (\al_{i_-},\al_k)=
0,$$ 
and hence
  $\De(\tK^+_j,\tK^+_k)\equiv1$.
  
  \mnoi
  {\bf Case}\  $j=i$ and $k=i_-$.\quad 
  Let us show
  $$\De\bl\Diz,\Daf^{-1}(\Qtp\ang{i_+}_w)\br\equiv z-w.$$

  We have a  monomorphism
  $$\ang{i\;i_+}_z\monoto\ang{i_+}_z\conv[z]\ang{i}_z.$$
  Hence we have
  $$\Daf^{-1}\bl\Qtp(\ang{i_+}_z)\br\conv[z]\Qtp(\ang{i\;i_+}_z)\monoto
  \Daf^{-1}\bl\Qtp(\ang{i_+}_z)\br\conv[z]\Qtp(\ang{i_+}_z)
  \conv[z]\Qtp(\ang{i}_z)
  \To[\ev]\Qtp(\ang{i}_z).$$
  It is an epimorphism since so is its specialization at $z=0$.
  
  Since  $\D(\Qtp\ang{i}) \simeq \Qt_+(\ang{i_+})\conv \tC_+^{\circ-1}$ commutes with $\Qtp(\ang{i_+ \:i})$,  $\Diz$ and $\Qtp(\ang{i_+\;i}_z)$ commute.
Hence  we have
  \eqn
  &&\Daf^{-1}\bl\Qtp(\ang{i_+}_z)\br\conv[z]\Diz\conv[z]\Qtp(\ang{i\;i_+}_{z})\\
 &&\hs{3ex}\simeq
  \Daf^{-1}\bl\Qtp(\ang{i_+}_z)\br\conv[z]\Qtp(\ang{i\;i_+}_z)\conv[z]\Diz
\epito \Qtp(\ang{i}_z)\conv[z]\Diz\epito \cor[z].
  \eneqn
  Hence, Proposition \ref{prop:defactor} implies that
  $z-w$ divides
  $$\De\bl\Daf^{-1}(\Qtp(\ang{i_+}_w)),\Diz\conv[z]\Qtp(\ang{i\;i_+}_z)\br
 \equiv \De\bl\Daf^{-1}(\Qtp(\ang{i_+}_w)),\Diz\br.$$
  Since they have the same homogeneous degree $2\sfd_i$, we obtain
  $$\De\bl\Daf^{-1}(\Qtp(\ang{i_+}_w)),\Diz\br\equiv z-w.$$
  
\mnoi
{\bf Case}\  $j\in I\setminus\st{i}$ and $k=i$.\quad
Let us show $$\De\bl\Diz,\tK^+_j\br\equiv Q_{i,j}(z,z_j),$$
where $z_j$ is a spectral parameter for $\tK^+_j$.

By Theorem~\ref{th:DiEi} below  (note that  Section \ref{Sec: gen for iRgmod} does not rely on the results in Section \ref{Sec: SW duality}), we have an exact sequence
\eqn
0\To\Daf(\Qtp\ang{i}_{z})\conv \tK_j^+
\To[\Rmat_{\sDiz}(\tK_j^+)] \tK_j^+\conv\Daf(\Qtp\ang{i}_{z})
\To\Qt_+(\E_i\tK_j^+)\To0
\eneqn
in $\Proc(\cor[z,z_j],\Loc_+)$.

Since $Q_{i,j}(z,z_j)\id_{\Qt_+(\E_i\tK_j^+)}$ vanishes, 
we have a commutative diagram
$$\xymatrix@C=20ex{
&  \tK_j^+\conv\Daf(\Qtp\ang{i}_{z})\ar[d]^{Q_{i,j}(z,z_j)}\ar@{.>}[dl]\\
\Daf(\Qtp\ang{i}_{z})\conv \tK_j^+\akew
\ar@{>->}[r]_{\Rmat_{\sDiz}(\tK_j^+)}& \tK_j^+\conv\Daf(\Qtp\ang{i}_{z})\;.
}$$
By Theorem \ref{th:ren_r_matrix} (a), the dashed morphism is a constant multiple of the renormalized R-matrix. 
Hence $\De\bl\Daf(\Qtp\ang{i}_{z}),\tK_j^+\br$ divides $Q_{i,j}(z,z_j)$.
Comparing their homogeneous degrees, we obtain
$$\De\bl\Daf(\Qtp\ang{i}_{z}),\tK_j^+\br\equiv Q_{i,j}(z,z_j).$$
\QED

      \Th\label{th:SW}
      The family $(\tK_j^\pm)_{j\in I_\mp}$ in $\Pro(\Loc_{\pm})$
      is a Schur-Weyl datum for $R_\mp$.
      \enth
\Proof
Since $\Loc_+$ is quasi-rigid, it is enough to show \eqref{eq:DeQ} of  \eqref{eq: duality daum}, which is nothing but Theorem \ref{thm:DeSW}. \QED

      Let $\Psi_{\pm,\mp}\cl R_\mp \gmod\to \Loc_\pm$
      be the associated Schur-Weyl duality functor.

      Then $\Psi_{\pm,\mp}$ induces an isomorphism
      $\psi_{\pm,\mp}\cl\rtl(R_\mp)\isoto\rtl(R_\pm)$ given by
      $$\psi_{\pm,\mp}(\al_j)=\bc
      s_i(\al_j)&\text{if $j\in I$,}\\
       -\al_{i_\pm}
       &\text{if $j=\al_{i_\mp}$.}
      \ec$$

\Lemma  \label{lem:Psi C}
We have $\Psi_{\pm,\mp}(C_\mp)\simeq(\tC_\pm)^{\circ-1}$.
\enlemma
\Proof We have
\eqn\Psi_{+,-}(C_-)\simeq\D^{-1}\Qtp(\ang{i_+})\hconv\Di
&&\simeq \bl\tC_+^{\circ-1}\conv\Qtp(\ang{i})\br\hconv\bl\Qtp(\ang{i_+})\conv\tC_+^{\circ-1}\br\\
&&\simeq\tC_+^{\circ-1}\conv\tC_+\conv \tC_+ ^{\circ-1}\simeq\tC_+^{\circ-1}.
\eneqn
\QED

\Prop
$\Psi_{+,-}$ factors as
$R_-\gmod\to\Loc_-\To[\F_i]\Loc_+$.
\enprop
\Proof
By Lemma \ref{lem:Psi C} and Theorem \ref{Thm: graded localization},
it remains to show that $\Psi_{+,-}(\coRr_{C_-}(M))\cl\Psi_{+,-}(M) \conv \Psi_{+,-}(C_-)  \to  \Psi_{+,-}(C_-) \conv \Psi_{+,-}(M)$
is an isomorphism for any $M\in R_-\gmod$. 

\medskip
We shall first show that
$$\text{$\Psi_{+,-}(\coRr_{C_-}(\ang{j}_z)) \cl \Psi_{+,-}(\ang{j}_z) \conv \tC_+  \to  \tC_+ \conv \Psi_{+,-}(\ang{j}_z)$ is an isomorphism for any $j\in I_-$.}$$
Since $\coRr_{C_-}(\ang{j}_z)$ is a rational isomorphism of homogeneous degree $0$,
so is $\Psi_{+,-}(\coRr_{C_-}(\ang{j}_z)) \cl \tK_j^+\conv \tC_+  \to  \tC_+ \conv \tK^+_j$. Since  $\La(K_j^+,\tC_+)=0$,
 $\Psi_{+,-}(\coRr_{C_-}(\ang{j}_z))$ is a renormalized R-matrix.
 Since $\tC_+$ is central in $\Loc_+$,
 the morphism $\Psi_{+,-}(\coRr_{C_-}(\ang{j}_z))$ is an isomorphism.
 
  \medskip
Since $\Psi_{+,-}$ is a monoidal functor,
$\Psi_{+,-}(\coRr_{C_-}(R_-(\beta)))$ is an isomorphism for any
$\beta\in\prtl(R_-)$.
Hence, 
the morphism $\Psi_{+,-}(\coRr_{C_-}(M))\simeq
\Psi_{+,-}(\coRr_{C_-}(R_-(\beta)))\tens_{R_-(\beta)} M$
is an isomorphism  for any $M\in R_-(\beta)\gmod$. 
\QED

Note that  we only know here  that $\Psi_{\pm,\mp}$ is right exact.
(The exactness of $\Psi_{\pm,\mp}$ follows from Proposition~\ref{prop:FF*}
below.) 
Let $\Ld_n\Psi_{\pm,\mp}$ be the left derived functor of  $\Psi_\pm\cl R_\mp
\gmod\to\Pro(\Locp)$. 
Since the family of $M\in R_\mp\gmod$ with finite projective dimension
  is $\Psip$-projective by Proposition~\ref{prop:F proj},  
  $\Ld_n\Psi_{\pm,\mp}(M)\simeq 0$ for any $n\not=0$  and  $M$ with finite projective dimension.

We extend $\Ld_n\Psi_{\pm,\mp}\cl\Pro(R_\mp\gmod)\to\Pro(\Loc_\pm)$
so that it commutes with codirected projective limits.

\Lemma\label{lem:bijSW}
We have
\bnum
\item
  $\Psi_{+,-}(\ang{j\,i^{-c_{i,j}}})\simeq \Qtp(\ang{j})$
  for any $j\in I\setminus\st{i}$,

\item
  $  \F_i\bl\D^{-1}(\Qt_-(\ang{i}))\br\simeq  q_i^{-1} \Qtp(\ang{i})$,
  \item
$  \F_i\bl\D(\Qt_-\ang{i_-})\br\simeq q_i \Qtp(\ang{i_+})$.
\ee
\enlemma

\Proof
Set $c=-c_{i,j}$.
\noi
(i)\
Note first that
$\ang{j\,i^c}$ is the image of
$\rmat{\ang{j},\ang{i^c}}\cl\ang{j}\conv\ang{i^c}\to q_i^{c^2}  \ang{i^c}\conv\ang{j}$.
Note also that
\eqn\Psi_{+-}(\ang{i^c})&&\simeq\Psi_{+-}(q_i^{c(c-1)/2}\ang{i}^{\circ c})
\simeq q_i^{c(c-1)/2}\bl  q_i \D\bl\Qt_+(\ang{i})\br^{\circ c}
\simeq q_i^{c(c-1)/2}q_i^c\D\bl \Qt_+(\ang{i}^{\circ c})\br\\
&&\simeq q_i^{c(c-1)/2}q_i^c\D\bl\Qt_+(q_i^{-c(c-1)/2}\ang{i^c})
\simeq q_i^{c^2}\D\bl\Qt_+(\ang{i^c})\br.
\eneqn

By applying the right exact functor $\Psip$
to $\ang{j}\conv\ang{i^c}\epito\ang{j\,i^c}$,
we obtain
$$\Psip(\ang{j})\conv\Psip(\ang{i^c})\epiTo\Psip(\ang{j\,i^c}).$$
It does not vanish by Lemma~\ref{lem:ijc}. 
Let us consider an exact sequence
$$0\To\ang{j\,i^c}\To  q_i^{-c^2} \ang{i^c}\conv\ang{j}\To N\To0.$$
Since $\ang{j\,i^c}$ and $\ang{i^c}\conv\ang{j}$ have finite projective dimension (see \cite{Murata}), 
so does $N$.
Hence $\Ld_1\Psip(N)\simeq0$,
and we obtain an exact sequence
$$0\To\Psip(\ang{j\,i^c})\To  q_i^{-c^2} \Psip(\ang{i^c})\conv\Psip(\ang{j}).$$
  Hence $\Psip(\ang{j\,i^c})$ is the image of
    the R-matrix $\Psip(\ang{j})\conv\Psip(\ang{i^c})
    \To  q_i^{-c^2}  \Psip(\ang{i^c})\conv\Psip(\ang{j})$,
    i.e.
    $$\Psip(\ang{j\,i^c}) \simeq  
 \Psip(\ang{j} \hconv \ang{i^c}) 
\simeq
\Psip(\ang{j})\hconv\Psip(\ang{i^c})\simeq
  q_i^{c^2} \Qtp(\ang{i^c\,j})\hconv \D(\Qtp(\ang{i^c})).$$

 Since $\ang{i^cj}\simeq q_i^{-c^2}\ang{j}\sconv\ang{i^c}$
(note that $\La(\ang{i^c},\ang{j})=c^2\sfd_i$), 
we have
$$\Qt_+(\ang{i^cj})\conv\D\bl\Qt_+(\ang{i^c})
\monoto q_i^{-c^2}\Qt(\ang{j})\conv\Qt_+(\ang{i^c})\conv\D\bl\Qt_+(\ang{i^c})).
\epito q_i^{-c^2}\Qt(\ang{j}).$$
Since the composition does not vanish,
we obtain
$$\Qt_+(\ang{i^cj})\hconv\D\bl\Qt_+(\ang{i^c})\br\simeq q_i^{-c^2}\Qt(\ang{j}).$$
It follows that 
\eqn
&&\Psi_{+-}(\ang{ji^c})\simeq q_i^{c^2}q_i^{-c^2}\Qt(\ang{j})
\simeq \Qt(\ang{j}).
\eneqn

    \snoi
    (ii)\ We have
    \eqn
    \F_i\bl\D^{-1}(\Qt_-(\ang{i}))\br&&\simeq\D^{-1}\F_i\bl\Qt_-(\ang{i})\br
    \simeq    \D^{-1}(\Psip(\ang{i}))\\
    &&\simeq
  \D^{-1}\bl q_i \D(\Qtp\ang{i})\br\simeq      q_i^{-1}   \Qtp(\ang{i}).\eneqn

    \snoi(iii) can be proved similarly.
\QED

\Lemma\label{lem:SW} We have
\bnum
\item 
  $\Psi_{+-}(\ang{j\,i^{-c_{i,j}}}_{z_j})\simeq  \Qtp(\ang{j}_{z_j})$ 
  for any $j\in I\setminus\st{i}$,
\item
  $  \F_i\bl\Daf^{-1}(\Qt_-(\ang{i}_z))\br\simeq   q_i^{-1} \Qtp(\ang{i}_z)$,
  \item
$  \F_i\bl\Daf(\Qt_-\ang{i_-}_z)\br\simeq   q_i  \Qtp(\ang{i_+}_z)$.
\ee
\enlemma
\Proof
In the proof below, we ignore the grading shifts.

(i)\ Set $c=-\sfc_{i,j}$.
We have an exact sequence
$$0\To\ang{j\,i^c}_{z_j}\To[z_j]\ang{j\,i^c}_{z_j}\To\ang{j\,i^c}\To0.$$
Since $\ang{j\,i^c}$ has a finite projective dimension (see \cite{Murata}),
        $\Ld_1\Psip(\ang{j\,i^c})\simeq0$ and hence
          we have an exact sequence
          $$0\To\Psip(\ang{j\,i^c}_{z_j})\To[z_j]\Psip(\ang{j\,i^c}_{z_j})
          \To\Psip(\ang{j\,i^c})\To0.$$
          Hence
          $\Psip(\ang{j\,i^c}_{z_j})\in\Aff(\Locp)$ and
          $\Psip(\ang{j\,i^c}_{z_j} )/z_j\Psip(\ang{j\,i^c}_{z_j})\simeq
          \Qtp(\ang{j} )$.

By Lemma~\ref{lem:cusp}, we have
$$\ang{j}_{z_j}\conv[z_j]\ang{i^c}_{z_j}\epiTo \ang{j\,i^{c}}_{z_j}
\monoTo \ang{i^c}_{z_j}\conv[z_j]\ang{j}_{z_j}.$$
Applying the right exact monoidal functor $\Psi_{+-}$, we obtain
\eq
\Psi_{+-}(\ang{j}_{z_j})\conv[z_j]\Psi_{+-}(\ang{i^c}_{z_j})
\epito \Psi_{+-}(\ang{j\,i^{c}}_{z_j}).\label{eq:epiPsij}
\eneq
On the other hand, we have  
\eqn
&&\Psi_{+-}(\ang{j}_{z_j})\conv[z_j]\Psi_{+-}(\ang{i^c}_{z_j})
\simeq \Qtp(\ang{i^c\,j}_{z_j} ) \conv[z_j] \D_{\cor[z_j]}  ( \Qtp(\ang{i^c}_{z_j}) )\\
&&\hs{8ex}\To
\Qtp(\ang{j}_{z_j})\conv[z_j]\Qtp(\ang{i^c}_{z_j})
\conv[z_j] \D_{\cor[z_j]}  (   \Qtp(\ang{i^c}_{z_j}))
\To\Qtp(\ang{j}_{z_j})\conv[z_j]\cor[z_j]
\simeq\Qtp(\ang{j}_{z_j}).
\eneqn
See Appendix \ref{sec:appendixC5} for the isomorphism
$$\Psip(\ang{i^c}_{z_j}) \simeq  \D_{\cor[z_j]} ( \Qtp(\ang{i^c}_{z_j})).$$

Since its specialization at $z_j=0$ is an epimorphism,
we obtain an epimorphism
$$\Psi_{+-}(\ang{j}_{z_j})\conv[z_j]\Psi_{+-}(\ang{i^c}_{z_j})
\epito\Qtp(\ang{j}_{z_j}).$$
Together with \eqref{eq:epiPsij}, Corollary~\ref{cor:qutaff} implies that
$$ \Psi_{+-}(\ang{j\,i^{c}}_{z_j})\simeq\Qtp(\ang{j}_{z_j}).$$

\snoi
(ii)\ We have
\eqn
\F_i\bl\Daf^{-1}(\Qt_-(\ang{i}_z))\br\simeq
\Daf^{-1}\F_i(\Qt_-(\ang{i}_z))&&\simeq
\Daf^{-1}\Psi_{+,-}(\ang{i}_z)\\
&&\simeq
\Daf^{-1}\Daf(\Qt_+(\ang{i}_z))\simeq\Qt_+(\ang{i}_z).
\eneqn

\snoi
The proof of (iii) is similar.
\QED

\Prop\label{prop:FF*}
$\F_i$ and $\F^*_i$ are  quasi-inverses to each other.

\enprop
\Proof
In order to see $\F_i\circ\F_i^*\simeq\id_{\Loc_+}$,
it is enough to show that
$\F_i\circ\Psi_{-,+}\simeq\Qt_+$, which follows from
Lemma~\ref{lem:SW}.
Indeed,  
Lemma ~\ref{lem:SW} implies that the functor $\F_i \circ \Psi_{-,+}$ is the restriction of the Schur-Weyl duality functor associated with the duality datum
$$\bl \{\Qtp(\ang{j}_{z_j}) \}_{j\in I_+}, \{\Rre_{\Qtp(\ang{j}_{z_j}),\Qtp(\ang{k}_{z_k})} \}_{j,k\in I_+}\br.$$
Now the claim follows from \cite[Lemma 7.7]{refl} on the uniqueness of Schur-Weyl duality functors. 
\QED

\Rem
 We have defined the
 functors $\Psi_{\pm,\mp}\cl R_{\mp}\gmod\to \Loc_\pm$
 by the duality data
 $\bl\st{(\tK_j^\pm,z_j)}_{j\in I_{\mp}}, \st{\Rmat_{\tK_j^\pm,\tK_{j'}^\pm}}_{j,j'\in I_\mp}\br$,
 and we have
 $$\xymatrix@C=15ex{ R_{+}\gmod\ar[r]_{\Qt_+}\ar@/^1.5pc/[rr]^{\Psi_{-+}}
   &\Loc_{+}\ar[r]_{\F^*_i}&\Loc_{-}\ar[r]_{\F_i}&\Loc_{+}.
   }
   $$
  We have proved
$\F_i\circ\Psi_{-+}\bl x_1\vert_{R_+(\al_j)}\br
=\Qt_+(x_1\vert_{R_+(\al_j)})$, i.e.\  the commutativity of
   \eqn
   &&  \xymatrix@C=18ex{
     \F_i\circ\Psi_{-+}\bl R_+(\al_j)\br\ar[r]^-\sim \ar[d]^{\F_i\circ\Psi_{-+}(x_1)}& \Qt_+(R_+(\al_j))\ar[d]^{\Qt_+(x_1)}\\
          \F_i\circ\Psi_{-+}\bl R_+(\al_j)\br\ar[r]^-\sim& \Qt_+(R_+(\al_j))
}
  \eneqn
  for any $j\in I_+$.
  Here $x_1$ is regarded as an element of $\END_{R_+(\al_j)}(R_+(\al_j))$.
  However we can only see that
$$\F_i\circ\Psi_{-+}\bl\tau_1\vert_{R_+(\al_j+\al_{j'})e(j,j')}\br
=c_{j,j'}\Qt_+(\tau_1\vert_{R_+(\al_j+\al_{j'})e(j,j')})$$
for some $c_{j,j'}\in \cor^\times$ ($j,j'\in I_+$).
(Note that $c_{j,j'}=1$ if $j=j'$ and $c_{j,j'}c_{j',j}=1$.)
Hence, in order to have an isomorphism 
$\F_i\circ \F^*_i\simeq\id_{\Loc_+}$ of functors of graded monoidal categories,
we have to normalize
$\Rmat_{\tK^-_j,\tK^-_{j'}}$ by replacing it
with $c_{j,j'}^{-1}\Rmat_{\tK^-_j,\tK^-_{j'}}$.
Note that, in the beginning,  $\Rmat_{\tK^\pm_j,\tK^\pm_{j'}}$ ($j\not=j'$)
are renormalized R-matrices chosen only by the condition
$\Rmat_{\tK^\pm_{j'},\tK^\pm_{j}}\circ\Rmat_{\tK^\pm_j,\tK^\pm_{j'}}=Q_{j,j'}(z_j,z_{j'})$ in $\END_{\Pro(\Loc_\pm)}(\tK^\pm_j\tens\tK^\pm_{j'})$.
\enrem

Summarizing, we have the following main result of this paper. 
      \Th\label{th:FLoc}
      The functor
      $\Psi_{\pm,\mp}\cl R_\mp\gmod\to \Loc_\pm$
      factors as
 $$R_-\gmod\to \Loc_-\To[\F_i]\Loc_+\qtq
      R_+\gmod\to \Loc_+\To[\F^*_i]\Loc_-$$ 
      Moreover, $\F_i\cl \Loc_-\to\Loc_+$ and
      $\F^*_i\cl\Loc_+\to\Loc_-$ are  quasi-inverses  to each other.
      \enth

\section{Generators of $(R\gMod) {}_i $ } \label{Sec: gen for iRgmod}

Let $R\gMod$ be the subcategory of $\Modg(R)$ consisting of  finitely generated graded $R$-modules. 
Set
\begin{equation} \label{Eq: iR and Ri}
\begin{aligned} 
(R\gMod)_i&\seteq\st{M\in R\gMod\mid \E_iM\simeq0}\qtq\\
{}_i (R\gMod)&\seteq \st{M\in R\gMod\mid \Es_iM\simeq0}.
\end{aligned}
\end{equation}

For $\nu=(\nu_1,\ldots,\nu_n)\in I^n$, we set
$$\dM(\nu)=R u(\nu),$$ where
$u(\nu)$ has the defining relations:
\begin{equation} \label{Eq: Mnu}
\begin{aligned}
&e(\nu)u(\nu)=u(\nu),\\
&\tau_1\cdots\tau_{k-1}u(\nu)=0\qt{if $1\le k\le n$ and
  $\nu_k=i$.}
\end{aligned}
\end{equation}
In particular, if $\nu_1=i$, then $\dM(\nu)=0$. 
\Lemma\label{lem:gen2}
$$\dM(\nu*j)\simeq\bc
\dM(\nu)\conv \ang{j}_{z_j}&\text{if $j\not=i$,}\\
\Coker\bl\ang{i}_z\conv\dM(\nu)\to\dM(\nu)\conv\ang{i}_z\br
&\text{if $j=i$.}
\ec
$$
Here $\nu*j$ is the concatenation of $\nu$ and $j$.
\enlemma

\Proof
By the definition, there exists a surjective $R$-module homomorphism
$$\dM(\nu) \conv\ang{j}_{z_j} \To[g] \dM(\nu *j)\quad \text{given by}\  u(\nu) \tens   u(j)_z  \mapsto u(\nu*j),$$
\ where $u(j)_z$ denotes a generating vector of $\ang{j}_{z_j}$.
Moreover, if $j\neq i$, then it is an isomorphism. 

Assume that $j=i$.
Then $\bl\ang{i}_z,\dM(\nu)\br$ is unmixed  and hence
we have a canonical morphism
$\rmat{}\cl\ang{i}_z\conv \dM(\nu)\to\dM(\nu)\conv\ang{i}_z$
of degree $\bl\al_i,\wt(\dM(\nu))\br$  which is given by 
$$\rmat{} (x\tens y) = \tau_1\ldots\tau_n   \;(y\tens x) \quad \text{for} \  x \in \ang{i}_z, \ y \in \dM(\nu).$$
Hence
there is an $R$-module homomorphism
$$\dM(\nu*j) \To[h] \Coker(\rmat{}),\quad  \text{given by} \ u(\nu*j) \mapsto  u(\nu) \tens  u(i)_z + \Im(\rmat{}).$$
Since  $g\circ \rmat{}=0$, there is an $R$-module homomorphism $\Coker(\rmat{}) \to \dM(\nu*j)$, which is an inverse of $h$, as desired.
\QED
\Lemma\label{lem:gen1}
For each $\nu \in I^n$, we have $\dM(\nu)\in(R\gMod)_i$,
and we have an isomorphism
$$\HOM_{R\gMod_i}(\dM(\nu),X)\isoto \HOM_{R\gMod}(Re(\nu),X)$$ for any $X\in(R\gMod)_i$, which is functorial in $X$.
In particular $\dM(\nu)$ is a projective object of
$(R\gMod)_i$, and any object of $(R\gMod)_i$ is a quotient of
a direct sum of $\dM(\nu)$'s. 
\enlemma
\Proof
Let $\pi_\n\cl Re(\nu) \epito \dM(\nu)$ be the epimorphism given by $e(\nu) \mapsto u(\nu)$. Then $-\circ \pi_\nu\cl \HOM_{R\gMod_i}(\dM(\nu),X)\monoto \HOM_{R\gMod}(Re(\nu),X)$ is a monomorphism,  which is  functorial in $X \in R\gMod$. 
Let $X\in (R\gMod)_i$. Then $\tau_1\cdots \tau_k e(\nu) X =0$  for $k$ such that  $\nu_k=i$, and hence every homomorphism $Re(\nu) \to X$ factors through $\pi_\nu$. Hence $-\circ \pi_\nu$ is surjective.
The last assertion follows from the analogous properties of  $Re(\nu)$ in $R\gMod$.
\QED

For $M\in R\gmod$, 
let $r\cl e(i,*)M\to \E_iM$
be the natural  $\cor$-linear map so that
\eq
&&x_kr(u)=r(x_{k+1}u)\qtq\tau_lr(u)=r(\tau_{l+1}u)\qt{for $u\in e(i,*)M$.}
\label{eq:r}
\eneq

Recall that, for $M,N\in \Mod(\cor[z],\Mod(R))$, $M\conv[z]N$
denotes the cokernel of $z\conv\id_N-\id_M\conv z\in\END(M\conv N)$. 

The proof of the following theorem is given in Appendix  \ref{sec:appendixB}
\Th\label{th:J}
 Let $M$ be an $R(\beta)$-module with $\beta\in\prtl$.
 Set $n=\height(\beta)$.
 \bnum
 \item

Then the $\cor$-linear map
$$J_M\cl M\to(\E_iM)\conv[z]\ang{i}_z$$
defined by
 $$J_M(u)=\sum_{a\in[1,n]; \nu_a=i}
 \tau_a\cdots\tau_{n-1}\bl r(\tau_1\cdots\tau_{a-1}u)\tens \ang{i}_z\br
 \qt{for $\nu\in I^\beta$ and $u\in e(\nu)M$.}$$
is $R$-linear and
homogeneous of weight $\bl\al_i,\wt(M)+\al_i\br$.
Here $\cor[z]$ acts on $\E_i(M)$ via $z\cdot r(w)\seteq r(x_1 w)$ for $w \in e(i,*)M$.
\item
  The composition
  $$\E_i(M)\To[\E_i(J_M)]\E_i\bl(\E_iM)\conv[z]\ang{i}_z\br
  \To(\E_iM)\conv[z](\E_i\ang{i}_z)\simeq\E_i(M)$$
  coincides with the identity,  where the second homomorphism comes from the \emph{shuffle lemma (\cite[Proposition 2.18]{KL09}).}
\ee
\enth

\Lemma
For any $\beta\in\rtl$, we have
$$\la_+(\al_{i_+},\beta)
=(\al_i,\beta).$$
  \enlemma
  \Proof
  For any $j\in I$, we have
  $\la_+(\al_{i_+},\al_j)=(\al_i,\al_j) $.
  \QED

\Lemma\label{lem:DiLa}
For any simple $M\in R\gmod$ such that $\Qtp(M)$ is simple, we have
\eqn
&&\La\bl \Qt_+(M), \Di \br=\La(\ang{i},M),\\
&&\La\bl\Di, \Qt_+(M)\br=-\bl \al_i,\wt(M)\br,\\
&&\de\bl\Di, \Qt_+(M)\br=\tLa(\ang{i},M).
\eneqn
Similarly, we have
$$\La\bl\Qt_-(M),\D^{-1}\Qt_-(\ang{i})\br=-\bl \al_i,\wt(M)\br.$$
\enlemma
\Proof
We have

\eqn
\La\bl \Qt_+(M), \Di \br=
\La\bl\Qtp(\ang{i}),\Qtp(M)\br=
\La(\ang{i},M),
\eneqn
where the last equality follows from $\ang{i}\in \Loco[+]$.
We have
\eqn
\La\bl\Di, \Qt_+(M)\br&&=
\La\bl\tC_+^{\circ-1}\conv\Qt_+(\ang{i_+}), \Qt_+(M)\br\\
&&=\La\bl\Qtp(\ang{i_+}), \Qt_+(M)\br
=\La\bl\ang{i_+},M\br\\
&&=\la_+\bl-\al_{i_+},\wt(M)\br
=-\bl\al_i,\wt(M)\br.
\eneqn
Note that
$\La\bl\Qtp(\ang{i_+}), \Qt_+(M)\br=\La\bl\ang{i_+},M\br$
by Lemma~\ref{lem:LaQ+2}, and
$\La\bl\ang{i_+},M\br=\la_+\bl-\al_{i_+},\wt(M)\br$
since $(\ang{i_+},M)$ is unmixed.
\QED

\Lemma\label{lem:Dizqiz} 
$\HOM_{\Proc( \cor[z],\Locp)}\bl\Diz  \conv[z]\qiz[z],\cor[z] \br \simeq0$.
\enlemma
\Proof
If it did not vanish, then $\HOM_{\Loc_+}(\Di\conv\qi,\one)$,
its specialization at  $z=0$,  would not vanish
by \cite[Proposition~3.7 (ii)]{refl} (set $H_0=0$). 
It contradicts $\Di\hconv\qi\simeq
\Qt_+(\ang{i_+\;i})\conv\tC_+^{\circ-1}\not\simeq\one$.
\QED

\Lemma
$\De(\Qiz[w],\Diz)\equiv z-w$.
\enlemma
\Proof
Since there is an epimorphism
$\Qiz\convz\Diz\epito\cor[z]$, Proposition  \ref{prop:defactor} 
implies that $z-w$ divides $\De(\Qiz[w],\Diz)$.
Since $\De(\Qiz[w],\Diz)$ has homogeneous degree
$2\de(\Qi,\Di)=2\sfd_i$, we obtain
the desired result.
\QED
\Lemma\label{lem:Dii}
We have an exact sequence in $\Proc(\cor[z,w], \Loc_+)$:
$$
\xymatrix@C=3ex@R=3ex{
 0\ar[r]&\Diz\conv\Qt_+(\ang{i}_w)\ar[r]&
 \Qt_+(\ang{i}_w)\conv \Diz\ar[r]&\dfrac{\cor[z,w]}{\cor[z,w](z-w)}\ar[r]&0\;.}$$
\enlemma
\Proof
Let $f\cl \Diz\conv\Qt_+(\ang{i}_w)\to\Qiz[w]\conv \Diz$ be the R-matrix.
 Note that $f$ is a monomorphism. 
Since $\De(\Diz,\Qiz[w])\equiv z-w$, we have a commutative diagram
$$
\xymatrix@C=4ex@R=3ex{&&\Qt_+(\ang{i}_w)\conv \Diz\ar[ld]\ar[d]^{z-w}\\
 0\ar[r]&\Diz\conv\Qt_+(\ang{i}_w)\ar[r]^f&
 \Qt_+(\ang{i}_w)\conv \Diz\ar[r]&\Coker(f)\ar[r]&0\;.}$$
Hence $\Coker(f)$ belongs to
$\Proc(\cor[z,w]/\cor[z,w](z-w),\Locp)$.
Since the composition
\eqn
\Diz\conv\Qt_+(\ang{i}_w)\To[f]\Qiz[w]\conv \Diz\To&&\Qiz\convz\Diz\\
&&\hs{3ex}\to \cor[z,w]/\cor[z,w](z-w)
\eneqn
vanishes by Lemma~\ref{lem:Dizqiz}, we obtain
$$\Coker(f)\to\cor[z,w]/\cor[z,w](z-w).$$

Applying
$\cor[z,w]/(\cor[z,w]z+\cor[z,w]w)\tens_{\cor[z,w]}\scbul$,
we obtain an exact sequence
$$\Di\conv\Qi\to\Qi\conv\Di\to \Coker(f)/z\Coker(f)\To0.$$
Since  there is an exact sequence 
$\ang{i_+} \conv \ang{i} \to \ang{i} \conv \ang{i_+} \to C_+ \to 0,$
we have an exact sequence
$$\Di\conv \Qi\To\Qi\conv\Di\To\one\To0,$$
and hence we obtain
$$\Coker(f)/z\Coker(f)\simeq\cor.$$
Hence $\Coker(f)\to\cor[z]$ is an epimorphism.
Let $N$ be the kernel of $\Coker(f)\to\cor[z]$.

Then applying $(\cor[z]/\cor[z]z)\tens_{\cor[z]}$
to the exact sequence
$0\to N\to\Coker(f)\to\cor[z]\to0$,
we obtain an exact sequence
$0\to N/zN\to\Coker(f)/z\Coker(f)\to \cor\to0$.
Hence $N/zN\simeq0$, which implies that $N\simeq0$.
Thus, we obtain the desired result.
\QED

Recall the choice of $\Rmat_{\sDiz}$ in Lemma~\ref{lem:normR}. 
\Th\label{th:DiEi}
There exists a short exact sequence in $\Proc(\cor[z],\Loc_+)$:
$$0\to\Daf(\Qtp\ang{i}_{z})\conv \Qt_+(M)
\To[{\Rmat_{\sDiz}(M)}]\Qt_+(M)\conv\Daf(\Qtp\ang{i}_{z})
\To[{J'_M}]\Qt_+(\E_iM)\to0.
$$
functorially in $M\in R\gmod$.
 Here $z$ acts on $\E_iM$ by $x_1\in\END\bl e(i,*)M\br$.
 \enth
 
\Proof
We write $\Rmat(M)$ for
$\Rmat_{\sDiz}(M)$. 
Then $R(M)$ is a monomorphism in $\Proc(\cor[z],\Loc_+)$ since it is a rational isomorphism. 

We have a morphism in $\Proc(\cor[z],\Locp)$
\eqn
J'_M\cl
\Qtp(M)\conv \Diz&&\To[\Qtp(J_M)]
\Qtp(\E_iM)\conv[z]\qiz\conv[z]\Diz\\
&&\hs{5ex}\To[\ev] \Qtp(\E_iM)
\eneqn
with homogeneous degree $\bl\al_i,\wt(M)+\al_i\br$, where $J_M$ is the morphism in Theorem \ref{th:J}.

Let us first show that the composition
\eq
\Diz\conv \Qt_+(M)
\To[\Rmat(M)] \Qt_+(M)\conv\Diz
\To[J'_M]\Qtp(\E_iM)
\label{eq:compM}
\eneq
vanishes.
We have
a commutative diagram

\scalebox{1}{
\parbox{\textwidth}{
\eqn
\xymatrix@C=8ex{
  \Diz\conv \Qtp(M)\ar[dd]^-{\Rmat(M)}\ar[r]^-{\Qtp(J_M)}&
    \Diz\conv[z]\Qtp(\E_i\,M)\conv[z]\qiz\ar[d]^{\Rmat(\E_iM)}\\
    &\Qtp(\E_i\,M)\conv[z]\Diz\conv[z]\qiz\ar[d]^-{\Rmat(\iz)}\ar@/^10pc/[dd]^-g\\
    \Qt_+(M)\conv\Diz\ar[r]^-{\Qtp(J_M)}&    \Qtp(\E_i\,M)\conv[z]\qiz\conv[z]\Diz
    \ar[d]^-{\ev}\\
&\Qtp(\E_iM).
}
\eneqn
}}
Since $\HOM_{\Aff(R_+\gmod)}(\Diz\conv[z]\qiz,\cor[z])\simeq0$
by Lemma~\ref{lem:Dizqiz},
$g$ vanishes: thus so does the composition ~\eqref{eq:compM}.

\smallskip
Let us set
$$\Phi(M)=\Coker\bl\Diz\conv \Qt_+(M)\To[\Rmat(M)] \Qt_+(M)\conv\Diz\br.$$
Note that $\Phi$ is an {\em exact} functor from
$R\gmod$ to $\Proc(\cor[z],\Loc_+)$.
Since
$\Rmat_{\sDiz}(M)$ is a rational isomorphism, $z^n\Phi(M)\simeq0$ for $n\gg0$,
which implies that $\Phi(M)\in\Loc_+$ by  \cite[Proposition 4.8 (iii-b)]{refl}. 
Thus,
$\Phi$ gives an exact functor
$$\Phi\cl R\gmod\to\Mod(\cor[z],\Loc_+).$$
The morphism $J'_M$ gives a morphism
$$I_M\cl \Phi(M)\to \Qtp(\E_i M)$$
in $\Mod(\cor[z],\Loc_+)$ functorial in $M\in R\gmod$.

Let us show that
\eq \text{$I_M$ is an epimorphism for any $M\in R\gmod$.}\label{eq:Iepi}
\eneq
First assume that $M$ is simple.
We may assume that $m\seteq\eps_i(M)>0$ so that $\E_i(M)\neq 0$. 
Hence $\tE_iM=\hd(\E_i M)$ is a simple module.
Then we have
$(\E_iM)\conv[z]\ang{i}_z\epito\tE_iM\conv\ang{i}$.
The composition
$$J_M^0\cl M\To[J_M] (\E_iM)\conv[z]\ang{i}_z\to \tE_iM\conv\ang{i}$$
does not vanish in $R\gmod$, which follows from
the commutative diagram
  $$\xymatrix@C=12ex{
    \E_iM\ar[r]^-{\E_i(J_M)}\ar@/_3pc/[ddr]|(.35){\id_{\E_iM}\akete[.5ex]\ake[-1ex]}&\E_i(\E_iM\convz\ang{i_z})\ar@{->>}[r]\ar[d]&\E_i(\tE_iM\conv\ang{i})\ar[d]\\
    &\E_iM\convz\E_i\ang{i_z}\ar[r]\ar[d]^\bwr
    &\tE_iM\conv\E_i\ang{i}\ar[d]^\bwr\\
    &\E_iM\ar[r]&\tE_iM.}
  $$

  \medskip
  The composition
  $$\Qtp(M)\conv\Di\To[\Qtp(J_M^0)] \Qtp(\tE_iM)\conv\Qtp(\ang{i})\conv \Di\To\Qtp(\tE_iM)$$
  is an epimorphism, since the left and right arrows are non-zero morphisms,
  and $\Qtp(\ang{i})$ is simple. 
By the commutative diagram
  $$\xymatrix{\Qtp(M)\conv\Diz\ar@{->>}[r]\ar@{->>}[d]&\Phi(M)\ar[r]^{I_M}&\Qtp(\E_iM)\ar[d]\\
   \Qtp(M)\conv\Di\ar[r]&\Qtp(\tE_iM)\conv\Qi\conv \Di\ar[r]&\Qtp(\tE_iM),
  }$$
  the composition
  $\Phi(M)\To[I_M]\Qtp(\E_iM)\epito \Qtp(\tE_iM)$ is an epimorphism.
      Since $\Qtp(\E_iM)$ has a simple head $\Qtp(\tE_iM)$ and the composition
  $\Phi(M)\To[I_M]\Qtp(\E_iM)\epito \Qtp(\tE_iM)$ is an epimorphism,
  we obtain that $I_M\cl\Phi(M)\to\Qtp(\E_iM)$ is an epimorphism
  for any simple $M$.
  Since $\Phi$ and $\E_i$ are exact functors, we conclude 
  \eqref{eq:Iepi} by induction on the length of $M$.

  \bigskip

For
$M,N\in R\gmod$,
let us consider a commutative diagram

\scalebox{.8}{\parbox{70ex}{
\eqn
\xymatrix@C=5ex@R=5ex{&&0\ar[d]&0\ar[d]\\
  0\ar[r]&\Diz\conv \Qtp(M)\conv \Qtp(N)\ar[r]^{R(M)}\ar[d]^-\bwr&\Qtp(M)\conv\Diz\conv \Qtp(N)
  \ar[r]\ar[d]_{R(N)}&\Phi(M)\conv \Qtp(N)\ar[r]\ar[d]&0\\
0\ar[r]&\Diz\conv \Qtp(M)\conv \Qtp(N)\ar[r]^{R(M\circ N)}&\Qtp(M)\conv \Qtp(N)\conv\Diz\ar[r]\ar[d]&\Phi(M\conv N)\ar[r]\ar[d]&0\\
&&\Qtp(M)\conv \Phi(N)\ar[r]^\sim\ar[d]&\Qtp(M)\conv \Phi(N)\ar[d]\\
&&0&0
  }
  \eneqn
}}

  Hence we obtain an exact sequence
  \eq
  0\To\Phi(M)\conv \Qtp(N)\To\Phi(M\conv N)\To \Qtp(M)\conv\Phi(N)\To0.\label{eq:exactMN}
  \eneq

  \medskip
  Now we shall show that
  $\Phi(M)\to \Qtp(\E_iM)$ is an isomorphism for any $\beta\in\prtl$ and
any  $M\in R(\beta)\gmod$ by induction on $\height{(\beta)}$.

  For $X\in\Locp$, we denote by $[X]$ the corresponding element in
  the Grothendieck group $K(\Locp)$.

  \snoi
  {\bf Step1}\ First $\Phi(M)\to \Qtp(\E_iM)$ is an isomorphism for $M=\ang{j}$
  with $j\in I$   (see Lemma~\ref{lem:Dii} when $j=i$).  
Hence $\Phi(M)\to \Qtp(\E_iM)$ is an isomorphism
  for any $M\in R(\al_j)\gmod$.

  \snoi
  {\bf Step2}\ Assume that $\height(\beta)>1$.

  Now let $\beta=\gamma+\delta$ with $\gamma,\delta\in\prtl\setminus\st{0}$.
  Let $M\in R(\gamma)$ and $N\in R(\delta)\gmod$.
  Then by the induction hypothesis, we have
  $[\Phi(M)]=[\Qtp(\E_i M)]$ and  $[\Phi(N)]=[\Qtp(\E_i N)]$.
  Hence we have
  \eqn
[\Phi(M\conv N)]&&=[\Phi(M)\conv \Qtp(N)]+[\Qtp(M)\conv\Phi(N)]\\
&&  =[\Qtp(\E_iM)\conv \Qtp(N)]+[\Qtp(M)\conv\Qtp(\E_iN)]=[\Qtp(\E_i(M\conv N))].
  \eneqn
  Therefore,
        $[\Phi(X)]=[\Qtp(\E_iX)]$ if $X\in R\gmod$ is isomorphic to
    $M\conv N$.
Hence $\Phi(X)\to\Qtp(\E_iX)$ is an isomorphism for such an $X$
by \eqref{eq:Iepi}.
Since any $M\in R(\beta)\gmod$ has a resolution 
$X_1\to X_0\to M\to 0$ such that $\Phi(X_k)\to\Qtp(\E_iX_k)$ is an isomorphism
for $k=0,1$,
we obtain a desired result: $\Phi(M)\to\Qtp(\E_iM)$ is an isomorphism
for any $M\in R\gmod$.
  \QED

\section{Reflection functor from $\Ri$ to $\iR$} \label{Sec: reflection}

\Th\label{th:FR}
We have  quasi-commutative diagrams:
\eqn
\xymatrix@C=10ex{
  (R\gmod){}_i\ar@{^{(}->}[d]\ar[r]_{\F_i}&{}_i(R\gmod)\ar@{^{(}->}[d]\\
  R_-\gmod\ar[d]_{\Qt_-}\ar[dr]|{\Psi_{+,-}}&R_+\gmod\ar[d]^{\Qt_+}\\
  \Loc_-\ar[r]_{\F_i}& \Loc_+ }&\qquad&
\xymatrix@C=10ex{
  {}_i(R\gmod)\ar@{^{(}->}[d]\ar[r]_{\F_i^*}&(R\gmod)_i\ar@{^{(}->}[d]\\
  R_+\gmod\ar[d]_{\Qt_+}\ar[dr]|{\Psi_{-,+}}&R_-\gmod\ar[d]^{\Qt_-}\\
  \Loc_+\ar[r]_{\F^*_i}&\; \Loc_-\;.}
\eneqn
Moreover, $\F_i\cl(R\gmod){}_i\to{}_i(R\gmod) $ and
$\F_i^*\cl{}_i(R\gmod)\to(R\gmod){}_i$ are quasi-inverses to each other.
\enth

In order to see Theorem~\ref{th:FR}, it is enough to show the following proposition (and its variant for $\F_i^*$) 
since $\Ri$ and $\iR$ are full subcategories of $\Loc_-$ and $\Loc_+$ via $\Qtm$ and $\Qtp$, respectively.

\Prop \label{prop:FiQtmRi}
We have $\F_i\bl\Qt_-(\Ri)\br\subset \Qt_+(\iR)$.
\enprop

\Proof
We shall first show that
$$\F_i\bl\Qt_-(\Ri)\br\subset\Qtp(R\gmod).$$
In order to see this it is enough to show the following statement:
\eq
\parbox{67ex}{if $M\in \Ri $  satisfies
  $\F_i\bl\Qt_-(M)\br  \in\Qtp(R\gmod)$,
  then
  $\F_i\bl\Qt_-(X)\br\in\Qtp(R\gmod)$, where
  $X=\Coker\bl\ang{i}_z\conv M\to M\conv\ang{i}_z\br$.}
\label{eq:FiR}
\eneq
Indeed  we have $\F_i(\Qt_-(\ang{j}_z)) \simeq \Psi_{+,-}(\ang{j}_z) \in \Qtp( \Pro(R\gmod))$ for $j\in I$,  $j \neq i$.
Assume \eqref{eq:FiR}.  Then by Lemma~\ref{lem:gen2}, we get $\F_i(\Qt_-(\dM(\nu))) \in \Qtp(\Pro(R\gmod))$ by induction on length of $\nu$. 
By Lemma~\ref{lem:gen1},  if $M \in\Ri$,  then $\F_i\bl\Qt_-(M))$ is a subquotient of $\soplus_k\F_i(\Qt_-(\dM(\nu_k)))$ for some $\nu_k$'s. 
Since $R\gmod$ is stable under taking subquotients in $\Pro(\Loc_+)$ by Proposition~\ref{prop:stab subq}, $\F_i\bl\Qt_-(M))$ belongs to $\Qtp(R\gmod)$.

\medskip
Let us prove \eqref{eq:FiR}. 
Let $N\in R\gmod$ be a module such that
$\F_i\bl\Qt_-(M)\br\simeq\Qtp(N)$.
Then, we have
\eqn
\F_i\bl\Qt_-(X)\br
&&\simeq\Coker\bl\Diz\conv \Qtp(N)\to\Qtp(N)\conv\Diz\br\\
&&\simeq\Qtp(\E_i N),
\eneqn
where the last isomorphism follows from Theorem~\ref{th:DiEi}.
Hence we obtain \eqref{eq:FiR}.

\vskip 1em
It remains to show that
\eq
\parbox{67ex}{if a simple $M\in\Ri$ and a simple $N\in R\gmod$ satisfy
  $\F_i\bl\Qt_-(M)\br\simeq\Qtp(N)$,
  then $N$ belongs to $\iR$.}
\label{eq:FiiR}
\eneq

We have
\eqn
\La(N,\ang{i})&&=\La\bl(\Qtp(N),\Qtp(\ang{i})\br
=\La(\Qt_-(M),\D^{-1}(\Qt_-\ang{i})\br
=-\bl\al_i,\wt(M)\br,
\eneqn
where the last equality follows from
Lemma~\ref{lem:DiLa}.
Hence we have
\eqn
2\tLa(N,\ang{i})&&=\La(N,\ang{i})-\bl\al_i,\wt(N)\br\\
&&=-\bl\al_i,\wt(M)\br-\bl\al_i,\;s_i\wt(M)\br=0.
\eneqn
Thus we obtain $N\in \iR$ by \cite[Lemma 2.15]{loc2}.
\QED

\Th\label{th:saito}
The algebra homomorphism
$K((R\gmod)_i)\to K({}_i(R\gmod))$ induced by $\F_i$
coincides with the braid symmetry homomorphism \ro see \cite{ext}\rf. 
\enth
\Proof
We have exact sequences (respecting the gradings)
$$0\to q_i^2\ang{i_+\,i}\to\ang{i}\conv\ang{i_+}\to C_+\to0
\qtq0\to C_+ \to\ang{i_+}\conv\ang{i}\to\ang{i_+\,i}\To0$$
in $R_+\gmod$.
Hence we have
$$[\Qtp\ang{i}]\cdot[\Qtp\ang{i_+}]
-q_i^2[\Qtp\ang{i_+}]\cdot[\Qtp\ang{i}]=(1-q_i^2)[\tC_+],$$
which implies that
\eq
&&[\Qtp\ang{i}]\cdot[\D \Qtp\ang{i}]
-q_i^2[\D\Qtp\ang{i}]\cdot[\Qtp\ang{i}]=1-q_i^2.\label{eq:bos}
\eneq
Since $\la_+(\al_{i_+},j)=(\al_i,\al_j)$ for $j\in I\setminus\st{i}$,
we have
$q^{(\al_i,\al_j)}\ang{i_+}\conv\ang{j}\simeq\ang{j}\conv\ang{i_+}$,
  which implies
  \eq
  [\Qtp\ang{j}]\cdot[\D\Qtp\ang{i}]=q^{(\al_i,\al_j)}  [\D\Qtp\ang{i}]\cdot[\Qtp\ang{j}].
  \label{eq:bos2}
  \eneq

Let  $ \g$ be the Kac-Moody algebra associate with the Cartan matrix $\cartan$ and let $\hcalA$ be the bosonic extension of type $\g$. As a $\Q(q)$-algebra, $\hcalA$ is generated by $\set{f_{j,m}}{j\in I, \ m \in \Z}$.  For the definition and properties of $\hcalA$ we refer \cite{ext}. 
Let us denote by $\B_i$ the $\Qq$-subalgebra of $\hcalA$
generated by $f_{i,1}$ and $\st{f_{j,0}\mid j\in I}$. 
Then \eqref{eq:bos} and \eqref{eq:bos2} imply that
there exists a $\Qq$-algebra homomorphism
$\B_i\to \Qq\tens_{\Zq} K(\Loc_+)$ given by
$$f_{j,0}\mapsto  q_j^{-1/2} [\Qtp\ang{j}] \quad (j\in I),\quad f_{i,1}\mapsto q_i^{-1/2} [\F_i(\Qt_-\ang{i})]=q_i^{1/2} [\D\Qtp\ang{i}].$$
Then we have a commutative diagram (see \cite[Section 3]{ext})
$$\xymatrix@C=6ex{
\hcalA[0]\ar[d]\ar[r]^{\TT_i}&\B_i\ar[d]\\
\Qq\tens_{\Zq}K(\Loc_-)\ar[r]^{\F_i}&\Qq\tens_{\Zq}K(\Loc_+)}
$$
Since $\TT_i\in\End(\hcalA)$ induces the braid symmetry homomorphism
$${\mathsf S}_i\cl A_{\Zq}(\mathfrak n)[i]\isoto A_{\Zq}(\mathfrak n)[i]^*$$
(for the definitions, see \cite[Section 1]{ext}), 
and $A_{\Zq}(\mathfrak n)[i] \simeq  K((R\gmod)_i)$,  $A_{\Zq}(\mathfrak n)[i]^* \simeq  K(_i(R\gmod))$,
we obtain Theorem~\ref{th:saito}.
\QED

\Conj
If $\g$ is finite type, then
$\F_i$ coincides with the reflection functor constructed in
\cite{refl}.
\enconj

\appendix

\section{Proof of Proposition~\ref{prop:affD}
} \label{sec:appendixA}
\subsection{The case of the right dual of $\Qi$}
In this subsection, we prove the following lemma.
\Lemma \label{lem:Qiz}
The affine object $\Diz$ has a rational center structure
  $\Rmat_{\sDiz}$.
\enlemma

The following lemma follows from the universal properties of the right duals.
\Lemma
Let $\dM$ be an object of $\Aff(\shc)$ which admits a right dual $\Da(\dM)$ in $\Aff(\shc)$.
Then there is an  isomorphism
$$\HOM_{\Proc(\cor[z],\shc)}(\dM \tens X, Y\tens \dM) \simeq \HOM_{\Proc(\cor[z],\shc)}( X\tens \Da(\dM),\Da(\dM) \tens Y) $$
which is functorial in $X,Y \in \shc$.
\enlemma
Hence $\Rmat_{\Qiz}(X)\cl \Qiz\conv X\to X\conv \Qiz$ induces a morphism
$$\Rmat'_{\sDiz}(X)\cl X\conv\Diz\to \Diz\conv X$$
functorial in $X\in\Loc_+$.
Hence, in order to prove Lemma~\ref{lem:Qiz}, it is enough to show that it is a rational isomorphism.
We may assume that $X=\Qtp(R(\nu))$ for $\nu\in I_+^n$, and then we can reduce to the case
$X=\Qtp(R(\al_j))$ and then $X=\Qi[j]$
for $j\in I_+$.
Hence it suffices to show that a non-zero morphism
$$\Rmat\cl \Qtp(\ang{j})\conv\Diz\to \Diz\conv\Qtp(\ang{j})$$
is a rational isomorphism for any $j\in I_+$.

Since $\Rmat$ does not vanish, replacing $\Rmat$ with $z^{-n}\Rmat$ for some $n$,
we may assume that
$\Rmat\vert_{z=0}\cl \Qi[j]\conv\Di\to \Di\conv \Qi[j]$
does not vanish.
If $j\not=i$, then $\Rmat\vert_{z=0}$ is an isomorphism because
$\Qtp(\ang{j})$ and $\Di$ commute. This implies that $\Rmat$ is an isomorphism.
Hence we may assume that $j=i$.

Thus we can assume that $I=\st{i}$, and hence $R_+$ is of type $A_2$.
Then we can embed $\Loc_+$ into the rigid monoidal category $(R_{A_2}\gmod)\,\widetilde{}$,
the localization of $R_{A_2}\gmod$ by the  \emph{quasi-centers,}; i.e., the real commuting family of left braiders associated with the simple modules $\ang{i i_+}$ and $\ang{i_+ i}$ (see \cite[Section 5.1]{loc1}).  
 More precisely, the functor
$\Loc_+\to(R_{A_2}\gmod)\,\widetilde{}$ is exact and fully faithful by \cite[Theorem 8.3]{refl}.

Since Proposition~\ref{prop:rational dual} implies that
the image of $\Diz$ in $(R_{A_2}\gmod)\,\widetilde{}$ is a rational center,
$\Diz$ is a rational center in $\Loc_+$.
Hence there exists a morphism 
$R'\cl\Qtp(\ang{j})\conv\Diz\to \Diz\conv\Qtp(\ang{j})$
in $\Proc(\cor[z],\Loc_+)$,
which is a rational isomorphism.

On the other hand,
since $\HOM(\Qi[j]\conv\Di, \Di\conv \Qi[j])=\cor(\Rmat\vert_{z=0})$,
applying \cite[Proposition 3.7]{refl} with $H_0=\cor[z]\Rmat$, we obtain
$$\HOM_{\cor[z]}\bl\Qi[j]\conv\Diz,\; \Diz\conv \Qi[j]\br=\cor[z]\Rmat.$$
Hence we have $R'\in\cor[z]\Rmat$, which implies that
$\Rmat$ is a rational isomorphism.

\subsection{The case of the left dual of $\Qtp(\ang{i_+})$}
Let us prove the following lemma.
\Lemma \label{lem:Dizp}
  There exists a rational center
  $(\Dipz,\Rmat_{\sDipz})$.
\enlemma

The morphism
$\Rmat_{\Qipz}(X)^{-1}\cl X\conv\Qipz \to \Qipz\conv X$ induces a morphism
$$\Rmat_{\sDipz}(X)\cl \sDipz \conv X\to X\conv \sDipz$$
functorial in $X\in\Loc_+$.
It is enough to show that it is a rational isomorphism.
By a similar argument as the one in the preceding subsection, it is enough to show that
a non-zero morphism
$$\Rmat\cl\Dipz\conv \Qtp(\ang{j})\to \Qtp(\ang{j})\conv\Dipz$$
is a rational isomorphism for any $j\in I_+$.
We may assume that
$\Rmat\vert_{z=0}\cl\Dip\conv\Qtp(\ang{j})\to \Qtp(\ang{j})\conv\Dip$
does not vanish.
If $j=i$ it is an isomorphism, since
$\Qi$ and $\Dip$ commute. This implies that $\Rmat$ is an isomorphism.
If $j=i_+$, we can assume that $I=\st{i}$, and hence $R_+$ is of type $A_2$.
Then we can embed $\Loc_+$ into the rigid monoidal category $(R_{A_2}\gmod)\,\widetilde{}$,
the localization of $R_{A_2}\gmod$ by the quasi-centers.
Hence the assertion follows from Proposition~\ref{prop:rational dual}
similarly to the case of $\Diz$. 

Thus we reduce to the case $j\in I\setminus\st{i}$.
If $(\al_i,\al_j)=0$, then $\ang{j}$ and $\D^{-1}\ang{i_+}$ commute and hence
$\Rmat$ is an isomorphism.
Hence we may assume that $(\al_i,\al_j)<0$,

We shall show the following lemma.
\Lemma\label{lem:extf}
For any $j\in I\setminus\st{i}$ such that $(\al_i,\al_j)<0$, there exist $n\in\Z_{\ge0}$
and
 a morphism
\eqn
&&f\cl \Ka\conv \Qtp(\ang{j})
\to \Qtp(\ang{j})\conv\Ka\eneqn
in $\Proc(\cor[z],\Loc_+)$
such that $\Im(f)$ contains
$\Qtp(\ang{j})\conv(z^n\Ka)\subset\Qtp(\ang{j})\conv\Ka$.
Here $\Ka\seteq\Dipz/z^{n+1}\Dipz\in \Proc(\cor[z],\Loc_+)$.
\enlemma

Admitting this lemma for a while, let us complete the proof of
Lemma~\ref{lem:Dizp}.
By \cite[Proposition 3.7]{refl},
$f=cz^k\Rmat\vert_{z^{n+1}=0}$ for some $c\in\cor^\times$ and $k\in[0,n]$.
Hence we have
$\Qi[j]\circ z^n\Dipz\subset\Im(\Rmat)+\Qi[j]\circ z^{n+1}\Dipz$.
Hence we have
$\Qi[j]\circ z^n\Dipz\subset\Im(\Rmat)$, which implies that
$\Rmat$ is a rational isomorphism. Thus we obtain Lemma~\ref{lem:Dizp}.

\bigskip
Finally let us prove Lemma~\ref{lem:extf}.
By the commutative diagram   (up to degree shifts)
\eqn
\xymatrix{
  z^n\Ka\akete\ar[r]^-\sim\ar@{^{(}->}[d]&\tC_+^{\circ(-n-1)}\conv
  \Qtp(C_+^{\circ(n)}\conv\E_{i_+}^*C_+)
  \ar@{^{(}->}[d]\\
\Ka\ar[r]^-\sim&\tC_+^{\circ(-n-1)}\conv\Qtp\bl\Es_{i_+}(C_+^{\circ(n+1)})\br,
}
\eneqn
it is enough to construct a morphism
$$g\cl \Qtp\bl\Es_{i_+}(C_+^{\circ(n+1)})\conv\ang{j}\br
\to\Qtp\bl\ang{j}\conv\Es_{i_+}(C_+^{\circ(n+1)})\br$$
in $\Loc_+$ such that
\eq
\Qtp\bl\ang{j}\conv C_+^{\circ n}\conv (\Es_{i_+}C_+)\br\subset\Im(g).
\label{eq:condg}
\eneq

By Theorem \ref{Thm: graded localization} (ii), 
 the morphism
$$ \Qtp(C_+\conv\ang{j})\To\Qtp(\ang{j}\conv C_+)$$
is an isomorphism in $\Loc_+$.
Hence,  its composition
$$\Qtp(C_+^{\circ n}\conv\ang{j})\To\Qtp(\ang{j}\conv C_+^{\circ n})$$
is an isomorphism.
Since the image $C_+^{\circ n}\conv\ang{j}\To\ang{j}\conv C_+^{\circ n}$ is the simple socle
$\ang{j}\sconv C_+^{\circ n}$ of $\ang{j}\conv C_+^{\circ n}$, we obtain
$$\Qtp\bl\ang{j}\conv C_+^{\circ n}\br=\Qtp\bl\ang{j}\sconv C_+^{\circ n}\br.$$
Hence \eqref{eq:condg} is equivalent to
\eq
\Qtp\bl(\ang{j}\sconv C_+^{\circ n})\conv \Es_{i_+}(C_+)\br\subset\Im(g).
\label{eq:condg1}
\eneq
Let
$$h\cl C_+^{\circ(n+1)}\conv\ang{j}\to \ang{j}\conv C_+^{\circ(n+1)}$$
be a non-zero morphism in $R_+\gmod$.
Then $\Im(h)$ coincides with the simple socle of
$\ang{j}\conv C_i^{\circ(n+1)}$.

Now let us take $n\gg0$ such that
$\ang{j}\sconv C_+^{\circ n}$ commutes with $C_+$ (see  \cite[Corollary 3.18]{loc1}).
Then, $(\ang{j}\sconv C_+^{\circ n})\conv C_+$ is a simple submodule of
$\ang{j}\conv C_i^{\circ(n+1)}$.
Hence we obtain
\eq\Im(h)=(\ang{j}\sconv C_+^{\circ n})\conv C_+.
\label{eq:Imh}
\eneq

Taking $\Es_{i_+}$, we obtain a morphism
$$h'\cl  \Es_{i_+}(C_+^{\circ(n+1)})\conv\ang{j}\simeq
\Es_{i_+}\bl C_+^{\circ(n+1)}\conv\ang{j}\br\To[\Es_{i_+}(h)]
\Es_{i_+}\bl\ang{j}\conv C_i^{\circ(n+1)}\br\simeq\ang{j}\conv\Es_{i_+}(C_i^{\circ(n+1)}).$$
By \eqref{eq:Imh}, we obtain
\eq\Im(h')=\Es_{i_+}\bl(\ang{j}\sconv C_+^{\circ n})\conv C_+\br
\supset (\ang{j}\sconv C_+^{\circ n})\conv  \Es_{i_+}(C_+).
\label{eq:Imh'}
\eneq
Thus
$$g=\Qtp(h')\cl  \Qtp\bl\Es_{i_+}(C_i\conv C_i)\conv\ang{j}\br
\to\Qtp\bl\ang{j}\conv\Es_{i_+}(C_i\conv C_i)\br$$
satisfies \eqref{eq:condg1}.
It shows Lemma~\ref{lem:extf}
and completes the proof of Lemma~\ref{lem:Dizp}.

\section{Proof of Theorem~\ref{th:J}} \label{sec:appendixB}

 \Lemma
 For $\nu\in I^n$ and $a,b\in\Z$ such that $1\le a\le b<n$
 and $k\in[1,n]$, we have
 \eqn
&& x_k\tau_a\tau_{a+1}\cdots \tau_b(\nu)\\
&&\hs{3ex}= \bc
  \tau_a\cdots \tau_bx_{b+1} e(\nu)
-\sum_{s=a}^{b}\delta(\nu_{s}=\nu_{b+1})(\tau_{a}\cdots
\tau_{s-1})(\tau_{s+1}\cdots\tau_b)e(\nu)&\text{if $k=a$,}\\
  \tau_a\cdots \tau_b  x_{k-1}e(\nu)+\delta(\nu_{k-1}=\nu_{b+1})
  (\tau_a\cdots \tau_{k-2}) (\tau_{k}\cdots \tau_b)e(\nu)
  &\text{if $a<k\le b+1$},\\
   \tau_a\cdots \tau_b  x_{k}e(\nu)&\text{otherwise.}
  \ec
  \eneqn
  \enlemma
  \Proof
 Let us first prove  the case $k=a$ by induction on $b-a$.
  We have
  \eqn
&&  x_a\tau_a\cdots \tau_b(\nu)\\
 && \hs{3ex}=\bl \tau_a x_{a+1}-\delta(\nu_a=\nu_{b+1})\br\tau_{a+1}\cdots \tau_be(\nu)\\
 &&\hs{3ex}=
 \tau_ax_{a+1}\tau_{a+1}\cdots \tau_b e(\nu)
 -\delta(\nu_a=\nu_{b+1})(\tau_{a+1}\cdots \tau_b)e(\nu)\\
 &&\hs{3ex}\udset{\ast}{=}
\tau_a \Bigl(\tau_{a+1}\cdots\tau_bx_{b+1}e(\nu)-\sum_{s=a+1}^b
\delta(\nu_s=\nu_{b+1})(\tau_{a+1}\cdots \tau_{s-1})(\tau_{s+1}\cdots\tau_b)
e(\nu)\Bigr)\\
&&\hs{30ex} -\delta(\nu_a=\nu_{b+1})(\tau_{a+1}\cdots \tau_b)e(\nu)\\
&&\hs{3ex} =  \tau_a\tau_{a+1} \cdots \tau_bx_{b+1} e(\nu)
-\sum_{s=a}^{b}\delta(\nu_{s}=\nu_{b+1})(\tau_{a}\cdots
\tau_{s-1})(\tau_{s+1}\cdots\tau_b)e(\nu).
\eneqn
Here $\udset{\ast}{=}$ follows from the induction hypothesis.

\snoi
 If $a<k\le b+1$, then we have
 \eqn
 x_k\tau_a\cdots \tau_b(\nu)&&=
 \tau_{a}\cdots\tau_{k-2}x_k\tau_{k-1}\cdots\tau_be(\nu)\\
 &&=  \tau_{a}\cdots\tau_{k-2}\bl \tau_{k-1}x_{k-1}+\delta(\nu_{k-1}=\nu_{b+1})\br
 \tau_k\cdots\tau_be(\nu)\\
 &&=  \tau_{a}\cdots\tau_bx_{k-1}e(\nu)
 +\delta(\nu_{k-1}=\nu_{b+1})(\tau_{a}\cdots\tau_{k-2})
 ( \tau_k\cdots\tau_b)e(\nu) .
  \eneqn
  \QED

  For a finite subset $A\subset[1,n]$, we write
  $\uprod_{s\in A}\tau_s=\tau_{j_1}\cdots\tau_{j_t}$,
  where $A=\st{j_1,\ldots,j_t}$ with $j_1<\cdots<j_t$.
  For $i,j\in I$ and $t\in\Z_{\ge1}$, we define
  $Q_{i,j}(u_1,\ldots,u_t;v)$ inductively in $t$ as follows:
    \eqn
    Q_{i,j}(u_1,\ldots,u_t;v)
    =\bc
    Q_{i,j}(u_1,v)&\text{if $t=1$,}\\
    \Bigl( Q_{i,j}(u_1,u_3,\ldots,u_t;v)
    -Q_{i,j}(u_2,u_3,\ldots,u_t;v)\Bigr)\Big/(u_1-u_2)
    &\text{if $t>1$.}
    \ec
    \eneqn
    Note that $Q_{i,j}(u_1,\ldots,u_t;v)$ is a symmetric polynomial in
    $(u_1,\ldots,u_t)$. For a non-empty finite set $A$ and
    a family of indeterminates $\st{u_a}_{a\in A}$,
    we write $Q_{i,j}(\st{u_a}_{a\in A};v)$ for
    $Q_{i,j}(u_{j_1},\ldots u_{j_{|A|}};v)$
    where $A=\st{j_1,\ldots, j_{|A|}}$.

    For $i\in I$, $\nu\in I^n$, and $a,b$ such that $1\le a\le b+1\le n$, we set
    $$A_i([a,b];\nu)\seteq\st{\sigma\mid \sigma\subset[a,b],\;
\text{$\nu_s=i$ for any $s\in \sigma$}}.$$
If there is no afraid of confusion, we write it
simply by $A_i[a,b]$.

    \Lemma\label{lem:tauQ} 
    Let $\nu\in I^n$, $a,b\in\Z$ and $\gamma\subset[1,n]$.
    Assume that $1\le a\le b+1\le n$ and $\st{k}\cup\gamma\subset
    [1,n]\setminus[a, b+1]$. Set $i=\nu_{b+1}$ and $j=\nu_{k}$.
    Then we have
    \eqn
    && \tau_a\cdots\tau_b\,
    Q_{i,j}(x_{b+1},\st{x_s}_{s\in\gamma};x_k)e(\nu)\\
      &&\hs{5ex}=\sum_{\sigma\in A_i[a,b]}
      Q_{i,j}(x_a,\st{x_{s+1}}_{s\in\sigma},\st{x_s}_{s\in\gamma};x_k)
      \bl\uprod_{s\in[a,b]\setminus\sigma}\tau_s\br e(\nu).      
    \eneqn
    \enlemma
    \Proof
    Let us prove it by the descending induction on $a$.
    By the induction hypothesis, we have
      \eqn
    && \tau_{a+1}\cdots\tau_b\,
    Q_{i,j}(x_{b+1},\st{x_s}_{s\in\gamma};x_k)e(\nu)\\
      &&\hs{5ex}=\sum_{\sigma\in A_i[a+1,b]}
      Q_{i,j}(x_{a+1},\st{x_{s+1}}_{s\in\sigma},\st{x_s}_{s\in\gamma};x_k)
      \bl\uprod_{s\in[a+1,b]\setminus\sigma}\tau_s\br e(\nu).      
      \eneqn

      Hence we have
          \eqn
    && \tau_a\cdots\tau_b\,
    Q_{i,j}(x_{b+1},\st{x_s}_{s\in\gamma};x_k)e(\nu)\\
      &&\hs{5ex}=\hs{-3ex}\sum_{\sigma\in A_i[a+1,b]}\tau_a
      Q_{i,j}(x_{a+1},\st{x_{s+1}}_{s\in\sigma},\st{x_s}_{s\in\gamma};x_k)
      \bl\uprod_{s\in[a+1,b]\setminus\sigma}\tau_s\br e(\nu)\\
            &&\hs{5ex}=\sum_{\sigma\in A_i[a+1,b]}\Bigl(
    Q_{i,j}(x_{a},\st{x_{s+1}}_{s\in\sigma},\st{x_s}_{s\in\gamma};x_k)\tau_a\\
&& \hs{20ex}           +\delta(\nu_a=i)
Q_{i,j}(x_{a},x_{a+1},\st{x_{s+1}}_{s\in\sigma},\st{x_s}_{s\in\gamma};x_k)\Bigr)
\bl\uprod_{s\in[a+1,b]\setminus\sigma}\tau_s\br e(\nu)\\[-.5ex]
            &&\hs{5ex}=\sum_{\sigma\in A_i[a,b]}
      Q_{i,j}(x_a,\st{x_{s+1}}_{s\in\sigma},\st{x_s}_{s\in\gamma};x_k)
      \bl\uprod_{s\in[a,b]\setminus\sigma}\tau_s\br e(\nu).      
\eneqn
\QED

  \Lemma
  For $\nu\in I^n$, $1\le a\le b<n$, and $k$ such that $1\le k<n$, we have
  \eqn
  &&(\tau_a\cdots\tau_b)\tau_ke(\nu)\\[2ex]
 &&\hs{3ex}=\left\{ \ba{ll}
  \tau_k(\tau_a\cdots\tau_b)e(\nu)&\text{if $k>b+1$,}\\
    (\tau_a\cdots\tau_b\tau_{b+1})e(\nu)&\text{if $k=b+1$,}\\[2ex]
    \sum\limits_{\sigma\in A_{\nu_b}[a,b-1]}Q_{\nu_b,\nu_{b+1}}
    (x_a,\st{x_{1+s}}_{s\in\sigma},;x_{b+1})
\bl\uprod_{s\in[a,b-1]\setminus\sigma}\tau_s\br e(\nu)\akew&\text{if $k=b$,}\\[4ex]
\tau_{k+1}(\tau_a\cdots\tau_b)e(\nu)\\
\hs{5ex}-  \delta(\nu_{k}=\nu_{b+1}) \hs{-3ex}
  \sum\limits_{\sigma\in A_{\nu_{b+1}}[a,k-1]}\\
\hs{15ex}  Q_{\nu_{b+1},\nu_{k+1}}(x_a,\st{x_{s+1}}_{s\in\sigma}, x_{k+2};x_{k+1})\\
\hs{35ex}
\bl\uprod_{s\in ([a,k-1]\setminus\sigma)\cup[k+2,b]}\tau_s\br e(\nu)\akew
&\text{if $a\le k<b$,}\\[3ex]
\tau_a(\tau_{a-1}\tau_{a+1}\cdots\tau_b)e(\nu)&\text{if $k=a-1$,}\\
\tau_k(\tau_a\cdots\tau_b)e(\nu)&\text{if $k<a-1$.}  \ea\right.
  \eneqn
  \enlemma
  \Proof
  Let us show by descending induction on $a$.
  
\snoi
{\bf Case $k=b$}.
We have
\eqn
&&(\tau_{a}\cdots\tau_b)\tau_be(\nu)\\
&&=(\tau_{a}\cdots\tau_{b-1})Q_{\nu_b,\nu_{b+1}}(x_b;x_{b+1})e(\nu)\\
&&=\sum\limits_{\sigma\in A_{\nu_b}[a,b-1]}Q_{\nu_b,\nu_{b+1}}
    (x_{a},\st{x_{s+1}}_{s\in\sigma};x_{b+1})
 \bl\uprod_{s\in[a,b-1]\setminus\sigma}\tau_s\br e(\nu).
 \eneqn
 Here the second equality follows from Lemma~\ref{lem:tauQ}.
 
\snoi
{\bf Case $a\le k<b$}.
We have
  \eqn
  &&(\tau_{a}\cdots\tau_b)\tau_ke(\nu)\\
  &&\hs{3ex}  =(\tau_a\cdots\tau_{k-1})
  (\tau_{k}\tau_{k+1}\tau_k)(\tau_{k+2}\cdots\tau_b)e(\nu)\\
  &&\hs{3ex}= (\tau_a\cdots\tau_{k-1}) \Bigl(\tau_{k+1}\tau_{k}\tau_{k+1}
  -\delta(\nu_{k}=\nu_{b+1})Q_{\nu_{b+1},\nu_{k+1}}
  (x_k,x_{k+2};x_{k+1})\Bigr)
  (\tau_{k+2}\cdots\tau_b)e(\nu)\\
  &&\hs{3ex}  =\tau_{k+1}(\tau_{a}\cdots\tau_b)e(\nu)
  -\delta(\nu_{k}=\nu_{b+1})\hs{-3ex}
  \sum_{\sigma\in A_{\nu_{b+1}}[a,k-1]}\hs{-3ex}Q_{\nu_{b+1},\nu_{k+1}}
  (x_a,\st{x_{s+1}}_{s\in\sigma};x_{k+2};x_{k+1})\\
&&\hs{55ex}(\uprod_{s\in ( ([a,k-1]\setminus\sigma ) \cup[k+2 ,b]}\tau_s)e(\nu).
  \eneqn
\QED

 \Proof[\bf Proof of Theorem \ref{th:J}]

 Let us show that $J_M$ is $R(\beta)$-linear.

 We set
 $$h_a(u)=
 \tau_a\cdots\tau_{n-1}\bl r(\tau_1\cdots\tau_{a-1}u)\tens \ang{i}_z\br.$$

 \snoi
 (A)\  Let us first show
 \eq e(\nu)J_M(u)=J_M(u)\qt{for any $u\in e(\nu)M$.}
 \eneq
 Set $v=r(\tau_1\cdots\tau_{a-1}u)\tens \ang{i}_z$.
 Then, we have
 $e(\nu_1,\ldots,\nu_{a-1},\nu_{a+1},\ldots, \nu_n,i)
 v=v$.
 Hence
 $e(\nu_1,\ldots,\nu_{a-1},i,\nu_{a+1},\ldots, \nu_n)
 \tau_a\cdots\tau_{n-1}v=
 \tau_a\cdots\tau_{n-1}v.$
 Thus we obtain $e(\nu)h_a(u)=h_a(u)$ if $\nu_a=i$.

 \mnoi
 (B)\  Assume that $u \in e(\nu) M$.  Let us show 
 \eq x_kJ_M(u)=J_M(x_ku)\qt{for $k\in[1,n]$.}
 \eneq

 Let $\mu\in I^n$ such that
$e(\mu)r(\tau_1\cdots\tau_{a-1}u)\tens \ang{i}_z=
r(\tau_1\cdots\tau_{a-1}u)\tens \ang{i}_z$.
Then $\mu_{n}=i$ and $\mu_s=\nu_{s+1}$ for $s\in[a,n-1]$,

 \snoi
 (Ba)\ 
 For $k>a$, we have
 \eqn
 && x_kh_a(u)=x_k\tau_a\cdots\tau_{n-1}\bl r(\tau_1\cdots\tau_{a-1}u)\tens \ang{i}_z\br\\
&&\hs{3ex} = \bl\tau_a\cdots\tau_{n-1}x_{k-1}+\delta(\mu_{k-1}=\mu_n)
(\tau_a\cdots \tau_{k-2}) (\tau_{k}\cdots \tau_{n-1})\br
  \bl r(\tau_1\cdots\tau_{a-1}u)\tens \ang{i}_z\br\\
&&\hs{3ex}=h_a(x_{k}u) +\delta(\nu_{k}=i)
(\tau_{a}\cdots \tau_{k-2})(\tau_{k}\cdots \tau_{n-1})
   \bl r(\tau_1\cdots\tau_{a-1}u)\tens
   \ang{i}_z\br\\
   &&\hs{3ex}=h_a(x_{k}u) +\delta(\nu_{k}=i)
(\tau_{k}\cdots \tau_{n-1})
   \bl r\bl(\tau_1\cdots\tau_{a-1})(\tau_{a+1}\cdots \tau_{k-1})u)\tens
   \ang{i}_z\br.
 \eneqn

 \snoi(Bb)\ 
 For $1\le k<a$, we have
 \eqn
x_kh_a(u)&&=
x_k\tau_a\cdots\tau_{n-1}\bl r(\tau_1\cdots\tau_{a-1}u)\tens \ang{i}_z\br
 =
 \tau_a\cdots\tau_{n-1}\bl r(x_{k+1}\tau_1\cdots\tau_{a-1}u)\tens \ang{i}_z\br,
\eneqn
 and we have
$$
x_{k+1}\tau_1\cdots\tau_{a-1}u
 =  \tau_1\cdots \tau_{a-1}  x_{k}u+\delta(\nu_{k}=\nu_a)
  (\tau_1\cdots \tau_{k-1}) (\tau_{k+1}\cdots \tau_{a-1} )u.
  $$
  Hence we obtain for $k<a$
  \eqn
x_kh_a(u)&&=
  \tau_a\cdots\tau_{n-1}\bl r(\tau_1\cdots\tau_{a-1}x_ku)\tens \ang{i}_z\br\\
&&\hs{20ex}+\delta(\nu_{k}=i)
\tau_a\cdots\tau_{n-1}\bl r\bl (\tau_1\cdots\tau_{k-1})(\tau_{k+1}\cdots \tau_{a-1})
u\br\tens \ang{i}_z\br\\
&&=h_a(x_ku)
+\delta(\nu_{k}=i)
(\tau_{k}\cdots \tau_{a-2})
(\tau_a\cdots\tau_{n-1})\bl r(\tau_1\cdots\tau_{k-1}u)\tens \ang{i}_z\br.
\eneqn

\snoi
(Bc)\ 
Finally assume that $k=a$.
Hence
\eqn
x_kh_a(u)&&=x_k\tau_k\cdots\tau_{n-1}\bl r(\tau_1\cdots\tau_{k-1}u)\tens \ang{i}_z\br\\
&&\hs{5ex}=\tau_k\tau_{k+1} \cdots \tau_{n-1}x_{n}
\bl r(\tau_1\cdots\tau_{k-1}u)\tens \ang{i}_z\br\\
&&\hs{10ex}-\sum_{s=k}^{n-1}\delta(\nu_{s+1}=i)(\tau_{k}\cdots
\tau_{s-1})(\tau_{s+1}\cdots\tau_{n-1})\bl r(\tau_1\cdots\tau_{k-1}u)\tens \ang{i}_z\br\\
&&\hs{5ex}=\tau_k\tau_{k+1} \cdots \tau_{n-1}x_{n}
\bl r(\tau_1\cdots\tau_{k-1}u)\tens \ang{i}_z\br\\
&&\hs{10ex}-\sum_{s=k+1}^{n}\delta(\nu_{s}=i)(\tau_{k}\cdots
\tau_{s-2})(\tau_{s}\cdots\tau_{n-1})\bl r(\tau_1\cdots\tau_{k-1}u)\tens \ang{i}_z\br.
\eneqn
Since we have
\eqn\tau_k\cdots \tau_{n-1}x_{n}
\bl r(\tau_1\cdots\tau_{k-1}u)\tens \ang{i}_z\br
&&=\tau_k\tau_{k+1} \cdots \tau_{n-1}
\bl r(\tau_1\cdots\tau_{k-1}u)\tens z\ang{i}_z\br\\
&&=\tau_k\tau_{k+1} \cdots \tau_{n-1}
\bl r(x_1\tau_1\cdots\tau_{k-1}u)\tens \ang{i}_z\br
\eneqn
and
\eqn
x_1\tau_1\cdots\tau_{k-1}u
=  \tau_1\cdots \tau_{k-1}x_{k}u
-\sum_{s=1}^{k-1}\delta(\nu_{s}=\nu_{k})(\tau_{1}\cdots
\tau_{s-1})(\tau_{s+1}\cdots\tau_{k-1})u,
\eneqn
we obtain
\eqn
&&\tau_k \cdots \tau_{n-1}x_{n}
\bl r(\tau_1\cdots\tau_{k-1}u)\tens \ang{i}_z\br\\
&&\hs{5ex}=\tau_k \cdots \tau_{n-1}
\bl r(\tau_1\cdots\tau_{k-1}x_ku)\tens \ang{i}_z\br\\
&&\hs{10ex}-\sum_{s=1}^{k-1}\delta(\nu_{s}=i)\tau_k\cdots \tau_{n-1}
\bl r\bl\tau_{1}\cdots
\tau_{s-1})(\tau_{s+1}\cdots\tau_{k-1})u\br\tens \ang{i}_z\br\\
&&\hs{5ex}=h_k(x_ku)
-\sum_{s=1}^{k-1}\delta(\nu_{s}=i)
(\tau_{s}\cdots\tau_{k-2})(\tau_k\cdots \tau_{n-1})
\bl r(\tau_{1}\cdots\tau_{s-1}u)\tens \ang{i}_z\br.
\eneqn

Hence
\eqn
x_kh_k(u)=
h_a(x_au)
&&-\sum_{s=1}^{k-1}\delta(\nu_{s}=i)
(\tau_{s}\cdots\tau_{k-2})(\tau_k\cdots \tau_{n-1})
\bl r(\tau_{1}\cdots\tau_{s-1})u)\tens \ang{i}_z\br\\
&&-\sum_{s=k+1}^{n}\delta(\nu_{s}=i)(\tau_{k}\cdots
\tau_{s-2})(\tau_{s}\cdots\tau_{n-1})\bl r(\tau_1\cdots\tau_{k-1}u)\tens \ang{i}_z\br.
\eneqn

\snoi
(Bd)\ 
Now, for any $k\in[1,n]$, we have
\eqn
&&x_kJ_M(u)-J_M(x_ku)
=\sum_{\nu_a=i}\bl
x_kh_a(u)-h_a(x_ku)\br\\
&&\hs{2ex}=
\delta(\nu_{k}=i)\sum_{a<k,\; \nu_a=i}
   (\tau_{k}\cdots \tau_{n-1})
  \bl r\bl(\tau_1\cdots\tau_{a-1})(\tau_{a+1}\cdots \tau_{k-1})u\br\tens \ang{i}_z\br\\
&&\hs{5ex}+\delta(\nu_{k}=i)\sum_{a>k,\;\nu_a=i}
\tau_a\cdots\tau_{n-1}\bl r\bl (\tau_1\cdots\tau_{k-1})(\tau_{k+1}\cdots \tau_{a-1})
u\br\tens \ang{i}_z\\
&&\hs{8ex}-\delta(\nu_k=i)\sum_{s<k}\delta(\nu_{s}=i)
(\tau_k\cdots \tau_{n-1})
\bl r\bl\tau_{1}\cdots\tau_{s-1})(\tau_{s+1}\cdots\tau_{k-1})u\br\tens \ang{i}_z\br\\
&&\hs{10ex}-\delta(\nu_k=i)\sum_{s>k}\delta(\nu_{s}=i)
(\tau_{s}\cdots\tau_{n-1})\bl r\bl(\tau_1\cdots\tau_{k-1})(\tau_{k+1}\cdots
\tau_{s-1})u\br\tens \ang{i}_z\br
\\
&&\hs{2ex}=0.
\eneqn

 \snoi
(C)\ Assume that $u \in e(\nu) M$. Let us show 
 \eq  J_M(\tau_ku)=\tau_kJ_M(u)\qt{for $k\in[1,n-1]$.}
 \eneq

 Since $\tau_ku\in e(s_k\nu)M$, we have
 \eqn
 J_M(\tau_ku)=\sum_{a\not=k,k+1,\;\nu_a=i}h_a(\tau_ku)
 +\delta(\nu_k=i)h_{k+1}(\tau_ku)
 +\delta(\nu_{k+1}=i)h_{k}(\tau_ku).
 \eneqn

 First we shall prove
\Sublemma\label{sub:Qq}
   If $\mu\in I^n$ satisfies $\mu_n=i$ and $\mu_s=\nu_{s+1}$
   for $s\in [a,n-1]$ and if
   $A\subset[1,a-1]$, $k\in[1,a-1]$ and $j=\mu_k$, then we have
   \eqn
   &&\tau_a\cdots\tau_{n-1}Q_{i,j}(x_n,\st{x_s}_{s\in A};x_k)e(\mu)
   =\hs{-2ex}\sum_{\sigma\in A_i([a+1,n],\nu)}\hs{-3ex}
   Q_{i,j}(\st{x_s}_{s\in\sigma\cup\st{a}},\st{x_s}_{s\in A};x_k)
(\hs{-2ex}\uprod_{s\in[a+1,n]\setminus\sigma}\hs{-2ex}\tau_{s-1})e(\mu).\eneqn
   \ensublemma
   \Proof
   We have
   \eqn
   \tau_a\cdots\tau_{n-1}Q_{i,j}(x_n,\st{x_s}_{s\in A};x_k)e(\mu)
=\hs{-2ex}\sum_{\gamma\in A_i([a,n-1];\mu)}\hs{-5ex}
   Q_{i,j}(x_a,\st{x_{s+1}}_{s\in\gamma},\st{x_s}_{s\in A};x_k)
   (\hs{-2ex}\uprod_{s\in[a,n-1]\setminus\gamma}\tau_{s})e(\mu).\eneqn
   We have a bijection $A_i([a,n-1];\mu)\isoto A_i([a+1,n];\nu)$ by
   $A_i([a,n-1];\mu)\ni\gamma\mapsto
   \sigma=\st{s+1\mid s\in\gamma}\in A_i([a+1,n];\nu)$.
     Hence, we have
   \eqn
   \tau_a\cdots\tau_{n-1}Q_{i,j}(x_n,\st{x_s}_{s\in A};x_k)e(\mu)
=\hs{-3ex}\sum_{\sigma\in A_i([a+1,n];\nu)}\hs{-5ex}
   Q_{i,j}(x_a,\st{x_{s}}_{s\in\sigma},\st{x_s}_{s\in A};x_k)
   (\hs{-2ex}\uprod_{s\in[a+1,n]\setminus\sigma}\hs{-2ex}\tau_{s-1})e(\mu).\eneqn
   
   \QED

   For $\sigma\in A_i[1,k-1]$ and $\gamma\in A_i[k+2,n]$, we set
   \eq
   &&\hs{5ex}\ba{l}
   X(\sigma,\gamma)\seteq
   Q_{i,\nu_{k+1}}(x_{k+1},\st{x_{s}}_{s\in\gamma},
\st{x_{s}}_{s\in\sigma},;x_{k})\\[.5ex]
\hs{30ex}\bl\uprod_{s\in[k+2,n]\setminus\gamma}\tau_{s-1}\br 
\Bigl( r\bl
(\uprod_{s\in[1,k-1]\setminus\sigma}\tau_s)u\br\tens\ang{i}_z\Bigr).   
\ea   \label{eq:X}
\eneq

\mnoi
 (Ca)\ Assume first that $a+1\le k\le n-1$ and $\nu_a=i$. Then we have
 \eqn
 h_a(\tau_ku)&&=
 \tau_a\cdots\tau_{n-1}\bl
 r(\tau_1\cdots\tau_{a-1}\tau_ku)\tens\ang{i}_z\br\\
&&= \tau_a\cdots\tau_{n-1}\tau_{k-1}\bl
r(\tau_1\cdots\tau_{a-1}u)\tens\ang{i}_z\br.
\eneqn
Take $\mu\in I^n$ such that
$e(\mu)\bl r(\tau_1\cdots\tau_{a-1} u  ) \tens\ang{i}_z\br 
 =\bl r(\tau_1\cdots\tau_{a-1} u)  \tens\ang{i}_z\br$. Then, we have
 $\mu_n=i$ and $\mu_s=\nu_{s+1}$ for $s\in[a,n-1]$.

We have
\eqn
&&\tau_a\cdots\tau_{n-1}\tau_{k-1}e(\mu)\\
&&=\tau_k\tau_a\cdots\tau_{n-1}e(\mu)-
\delta(\mu_{k-1}=i)\sum_{\sigma\in A_i([a,k-2];\mu)}
  Q_{i,\mu_{k}}(x_a,\st{x_{s+1}}_{s\in\sigma}, x_{k+1};x_{k})\\
&&\hs{35ex}
\bl\uprod_{s\in ([a,k-2]\setminus\sigma)\cup[k+1,n-1]}\tau_s\br\Bigr)\,e(\mu).
\eneqn
Note that $A_i([a,k-2];\mu)\isoto A_i([a+1,k-1];\nu)$ by
$A_i([a,k-2];\mu)\ni\sigma\mapsto\st{s+1\mid s\in\sigma}\in A_i([a+1,k-1];\nu)$.
Hence we have
\eqn
&&\tau_a\cdots\tau_{n-1}\tau_{k-1}e(\mu)\\
&&=\tau_k(\tau_a\cdots\tau_{n-1})  e(\mu) -\delta(\nu_k=i)
\sum_{\sigma\in A_i([a+1,k-1];\nu)}
  Q_{i,\nu_{k+1}}(x_a,\st{x_{s}}_{s\in\sigma}, x_{k+1};x_{k})\\
&&\hs{35ex}
\bl\uprod_{s\in [a+1,k-1]\setminus\sigma}\tau_{s-1}\br\tau_{k+1}\cdots\tau_{n-1}
  e(\mu).
  \eneqn

  Hence we have  
  \eqn
&&h_a(\tau_ku)-\tau_kh_a(u)
=-  \delta(\nu_{k}=i) 
  \sum\limits_{\sigma\in A_i[a+1,k-1]}
Q_{i,\nu_{k+1}}(x_a,\st{x_{s}}_{s\in\sigma}, x_{k+1};x_{k})\\
&&\hs{30ex}
\bl\uprod_{s\in ([a+1,k-1]\setminus\sigma)}\tau_{s-1}\br
\tau_{k+1}\cdots\tau_{n-1}\bl r(\tau_1\cdots\tau_{a-1}u)\tens\ang{i}_z\br\\
&&=-  \delta(\nu_{k}=i) 
  \sum\limits_{\sigma\in A_i[a+1,k-1]}
Q_{i,\nu_{k+1}}(\st{x_{s}}_{s\in\sigma\cup\st{a}}, x_{k+1};x_{k})
\tau_{k+1}\cdots\tau_{n-1}\\
&&\hs{40ex}\Bigl( r\bl(\uprod_{s\in [a+1,k-1]\setminus\sigma}\tau_s)\tau_1\cdots\tau_{a-1}u\br\tens\ang{i}_z\Bigr)\\
&&=-  \delta(\nu_{k}=i) 
  \sum\limits_{\sigma\in A_i[a+1,k-1]}
Q_{i,\nu_{k+1}}(\st{x_{s}}_{s\in\sigma\cup\st{a}}, x_{k+1};x_{k})
\tau_{k+1}\cdots\tau_{n-1}\\
&&\hs{40ex}\Bigl( r\bl(\uprod_{s\in [1,k-1]\setminus(\sigma\cup\st{a})}\tau_s)u)
  \tens\ang{i}_z\Bigr).
 \eneqn

 Hence we obtain
  \eqn
&&h_a(\tau_ku)-\tau_kh_a(u)
=-  \delta(\nu_{k}=i)
\sum_{\sigma\in A_i[a+1,k-1]}X(\sigma\cup\st{a},\emptyset).
\eneqn
 \snoi
 (Cb)\ We have
 \eqn
\delta(\nu_{k+1}=i) h_k(\tau_ku)&&=
 \delta(\nu_{k+1}=i) \tau_k\cdots\tau_{n-1}\bl
 r(\tau_1\cdots\tau_{k-1}\tau_ku)\tens\ang{i}_z\br\\
&&= \delta(\nu_{k+1}=i) \tau_{k}h_{k+1}(u).
\eneqn

  \snoi
 (Cc)\ We have
 \eqn
&&\hs{-3ex} \delta(\nu_{k}=i)h_{k+1}(\tau_ku)=\delta(\nu_{k}=i)
 \tau_{k+1}\cdots\tau_{n-1}\bl
 r(\tau_1\cdots\tau_{k}\tau_{k}u)\tens\ang{i}_z\br\\
 && =\delta(\nu_{k}=i)\hs{-2ex}   \sum\limits_{\sigma\in A_{i}[1,k-1]}
 \hs{-2ex}
\tau_{k+1}\cdots\tau_{n-1}\Bigl( r\bl
 Q_{i,\nu_{k+1}}(x_1,\st{x_{1+s}}_{s\in\sigma},;x_{k+1})
 (\uprod_{s\in[1,k-1]\setminus\sigma}\tau_s)u\br\tens\ang{i}_z\Bigr)\\
 && =\delta(\nu_{k}=i)   \sum\limits_{\sigma\in A_{i}[1,k-1]}
\tau_{k+1}\cdots\tau_{n-1}
  Q_{i,\nu_{k+1}}(x_n,\st{x_{s}}_{s\in\sigma},;x_{k})
\Bigl( r\bl
(\uprod_{s\in[1,k-1]\setminus\sigma}\tau_s)u\br\tens\ang{i}_z\Bigr)\,.
\eneqn
Here the second equality follows from Lemma~\ref{lem:tauQ}.

Take $\mu\in I^n$ such that $e(\mu)  \bl r  \bl
(\uprod_{s\in[1,k-1]\setminus\sigma}\tau_s)u\br\tens\ang{i}_z\br
 =r\bl
(\uprod_{s\in[1,k-1]\setminus\sigma}\tau_s)u\br\tens\ang{i}_z $.
Then, we have
$\mu_n=i$ and $\mu_s=\nu_{s+1}$ for $s\in[k+1,n-1]$.
Hence Sublemma~\ref{sub:Qq} implies that
\eqn
&&\delta(\nu_{k}=i)\tau_{k+1}\cdots\tau_{n-1}
  Q_{i,\nu_{k+1}}(x_n,\st{x_{s}}_{s\in\sigma},;x_{k})e(\mu)\\
&&\hs{10ex}=\delta(\nu_{k}=i)\sum_{\gamma\in A_i[k+2,n]}
  Q_{i,\nu_{k+1}}(x_{k+1},\st{x_{s}}_{s\in\gamma},
  \st{x_{s}}_{s\in\sigma};x_{k})
  \bl\uprod_{s\in[k+2,n]\setminus\gamma}\tau_{s-1}\br  e(\mu).
\eneqn

Thus we obtain
\eqn
&&\delta(\nu_{k}=i)h_{k+1}(\tau_ku)\\
&&\hs{15ex}=\delta(\nu_{k}=i)\hs{-2ex}\sum_{\substack{\gamma\in A_i[k+2,n]\\\sigma\in  A_{i}[1,k-1]}}\hs{-2ex}
Q_{i,\nu_{k+1}}(x_{k+1},\st{x_{s}}_{s\in\gamma},
\st{x_{s}}_{s\in\sigma},;x_{k})\\
&&\hs{30ex}\bl\uprod_{s\in[k+2,n]\setminus\gamma}\tau_{s-1}\br 
\Bigl( r\bl
(\uprod_{s\in[1,k-1]\setminus\sigma}\tau_s)u\br\tens\ang{i}_z\Bigr).
\eneqn
In terms of $X(\sigma,\gamma)$ in \eqref{eq:X},
we have
\eqn
&&\delta(\nu_{k}=i)h_{k+1}(\tau_ku)
=\delta(\nu_{k}=i)\hs{-2ex}\sum_{\substack{\gamma\in A_i[k+2,n]\\\sigma\in  A_{i}[1,k-1]}}X(\sigma,\gamma).
\eneqn

\snoi
(Cd)\ We have
  \eqn
\delta(\nu_{k}=i) \tau_{k}h_k(u)&&=
 \delta(\nu_{k}=i) \tau_k\tau_k\cdots\tau_{n-1}\bl
 r(\tau_1\cdots\tau_{k-1}u)\tens\ang{i}_z\br\\
 &&= \delta(\nu_{k}=i) Q_{i,\nu_{k+1}}(x_{k+1},x_{k})\tau_{k+1}
 \cdots\tau_{n-1}
 \bl r(\tau_1\cdots\tau_{k-1}u)\tens\ang{i}_z\br\\
&& =\delta(\nu_{k}=i)X(\emptyset,\emptyset).
\eneqn

\snoi
(Ce)\ Assume that $1\le k<a-1$, then we have
\eqn
&&\tau_1\cdots\tau_{a-1}\tau_{k}u\\
&&\hs{3ex}=\tau_{k+1}(\tau_1\cdots\tau_{a-1})u
-  \delta(\nu_{k}=i) \hs{-3ex}
  \sum\limits_{\sigma\in A_{i}[1,k-1]}\hs{-3ex}
  Q_{i,\nu_{k+1}}(x_1,\st{x_{s+1}}_{s\in\sigma}, x_{k+2};x_{k+1})
   \\&&\hs{35ex}
  \bl\hs{-3ex}
\uprod_{s\in([1,k-1]\setminus\sigma)\cup[k+2,a-1]}\tau_s\br u.
\eneqn
Hence we have
\eqn
&&r(\tau_1\cdots\tau_{a-1}\tau_{k}u)\tens\ang{i}_z\\
&&=\tau_{k}\bl r(\tau_1\cdots\tau_{a-1})u\tens\ang{i}_z\br
-\delta(\nu_{k}=i) \hs{-3ex}
\sum\limits_{\sigma\in A_{i}[1,k-1]}\hs{-3ex}
Q_{i,\nu_{k+1}}(x_n,\st{x_{s}}_{s\in\sigma}, x_{k+1};x_{k})
\\&&
\hs{30ex}
 \bl r(\hs{-3ex}
  \uprod_{s\in ([1,k-1]\setminus\sigma)\cup[k+2,a-1]}\tau_s u)
  \tens\ang{i}_z\br.
  \eneqn
  Take $\mu\in I^n$ such that
  $e(\mu)\bl  r(\hs{-3ex}
  \uprod_{s\in ([1,k-1]\setminus\sigma)\cup[k+2,a-1]}\tau_s u  ) \tens \ang{i}_z ) =r(\hs{-3ex}
  \uprod_{s\in ([1,k-1]\setminus\sigma)\cup[k+2,a-1]}\tau_s u)  \tens \ang{i}_z$.
  Then, we have
  $\mu_{n}=i$, $\mu_s=\nu_{s+1}$ for $s\in[a+1,n-1]$.
  Hence, Lemma~\ref{lem:tauQ} implies that
  \eqn
  &&\delta(\nu_{k}=i)\tau_a\cdots\tau_{n-1}
  Q_{i,\nu_{k+1}}(x_n,\st{x_s}_{s\in\sigma}, x_{k+1};x_k)e(\mu)\\
  &&=\delta(\nu_{k}=i)\sum\limits_{\gamma\in A_{i}[a+1,n]}
  Q_{i,\nu_{k+1}}(x_a,\st{x_{s}}_{s\in\gamma},
\st{x_s}_{s\in\sigma}, x_{k+1};x_k)
\bl\uprod_{s\in[a+1,n]\setminus\gamma}\tau_{s-1}\br e(\mu).
\eneqn
Hence
\eqn
&&h_a(\tau_ku)=\\
&&\hs{5ex}=(\tau_a\cdots\tau_{n-1}\tau_k)r(\tau_1\cdots\tau_{a-1}u)\\
&&\hs{15ex}-\delta(\nu_{k}=i)\sum\limits_{\substack{\gamma\in A_{i}[a+1,n]\\\sigma\in A_i[1,k-1]}}
   Q_{i,\nu_{k+1}}(x_a,\st{x_{s}}_{s\in\gamma},
\st{x_s}_{s\in\sigma}, x_{k+1};x_k)\\
&&\hs{15ex}\bl\uprod_{s\in[a+1,n]\setminus\gamma}\tau_{s-1}\br
\bl r(\hs{-3ex}
  \uprod_{s\in ([1,k-1]\setminus\sigma)\cup[k+2,a-1]}\tau_s u)
  \tens\ang{i}_z\br\\
  &&\hs{5ex}  =\tau_kh_a(u)\\
  &&\hs{15ex}-\delta(\nu_{k}=i)
  \sum\limits_{\substack{\gamma\in A_{i}[a+1,n]\\\sigma\in A_i[1,k-1]}}
   Q_{i,\nu_{k+1}}(\st{x_{s}}_{s\in\gamma\cup\st{a}},
\st{x_s}_{s\in\sigma}, x_{k+1};x_k)\\
&&\hs{25ex}\bl\uprod_{s\in[a+1,n]\setminus\gamma}\tau_{s-1}\br
\bl\uprod_{s\in[k+2,a-1]}\tau_{s-1}\br
\bl r(\hs{-3ex}
  \uprod_{s\in [1,k-1]\setminus\sigma}\tau_s u)
  \tens\ang{i}_z\br\\
    && \hs{5ex} =\tau_kh_a(u)\\
  &&\hs{15ex}-\delta(\nu_{k}=i)\sum\limits_{\substack{\gamma\in A_{i}[a+1,n]\\\sigma\in A_i[1,k-1]}}
   Q_{i,\nu_{k+1}}(\st{x_{s}}_{s\in\gamma\cup\st{a}},
\st{x_s}_{s\in\sigma}, x_{k+1};x_k)\\
&&\hs{25ex}\bl\uprod_{s\in[k+2,n]\setminus(\gamma\cup\st{a})}\tau_{s-1}\br\
\bl r(\hs{-3ex}
  \uprod_{s\in [1,k-1]\setminus\sigma}\tau_s u)
  \tens\ang{i}_z\br.
  \eneqn

  Hence if $1\le k< a-1$, then we obtain
  \eqn
  &&  h_a(\tau_ku)-\tau_kh_a(u)
  =-\delta(\nu_{k}=i)\sum\limits_{\substack{\gamma\in A_{i}[a+1,n]\\\sigma\in A_i[1,k-1]}}X(\sigma,\gamma\cup\st{a}).
  \eneqn
  
\mnoi
(Cf)\ 
Now, we shall calculate $J_M(\tau_ku)-\tau_kJ_M(u)$.
 \eqn
 &&J_M(\tau_ku)-\tau_kJ_M(u)\\
&&=\sum_{a=1}^{k-1}\delta(\nu_a=i)\bl h_a(\tau_ku)-\tau_kh_a(u)\br\\
&&\hs{3ex}+\delta(\nu_k=i)h_{k+1}(\tau_ku)+\delta(\nu_{k+1}=i)h_k(\tau_ku)
-\delta(\nu_k=i)\tau_kh_k(u)-\delta(\nu_{k+1}=i)\tau_kh_{k+1}(u)\\
&&\hs{10ex}+\sum_{a=k+2}^{n-1}\delta(\tau_a=i)\bl h_a(\tau_ku)-\tau_kh_a(u)\br.
\eneqn
We have
\eqn
\sum_{a=1}^{k-1}\delta(\nu_a=i)\bl h_a(\tau_ku)-\tau_kh_a(u)\br
&&=-  \delta(\nu_{k}=i) 
 \sum_{a=1}^{k-1}\delta(\nu_a=i) \sum\limits_{\sigma\in A_i[a+1,k-1]}
 X(\sigma\cup\st{a}, \emptyset )\\
&&=-\delta(\nu_k=i)
\sum\limits_{\sigma\in A_i[1,k-1]\setminus\st{\emptyset}}
X(\sigma, \emptyset ),
\eneqn
and
\eqn
\sum_{a=k+2}^{n-1}\delta(\nu_a=i )\bl h_a(\tau_ku)-\tau_kh_a(u)\br
&&= - \delta(\nu_{k}=i)   \sum_{a=k+2}^{n-1}\delta(\nu_a=i) \sum\limits_{\substack{\gamma\in A_{i}[a+1,n]\\\sigma\in A_i[1,k-1]}}
X(\sigma,\gamma\cup\st{a})\\
 &&= -\delta(\nu_{k}=i)\sum\limits_{\substack{\gamma\in A_{i}[k+2,n]\setminus\st{\emptyset}\\\sigma\in A_i[1,k-1]}}
X(\sigma,\gamma).  \eneqn

  Hence, we have
\eqn
&&J_M(\tau_ku)-\tau_kJ_M(u)\\
&&\hs{8ex}=- \delta(\nu_k=i)\hs{-4ex}\sum\limits_{\sigma\in A_i[1,k-1]\setminus\st{\emptyset}}\hs{-4ex}
X(\sigma,\emptyset)
+\delta(\nu_{k}=i)\hs{-2.5ex}\sum_{\substack{\gamma\in A_i[k+2,n]\\\sigma\in  A_{i}[1,k-1]}}\hs{-2.5ex}
X(\sigma,\gamma)\\
&&\hs{30ex}- \delta(\nu_{k}=i)X(\emptyset,\emptyset)
 - \delta(\nu_k=i)  \hs{-3ex} \sum\limits_{\substack{\gamma\in A_{i}[k+2,n]\setminus\st{\emptyset}\\\sigma\in A_i[1,k-1]}}\hs{-4ex}
 X(\sigma,\gamma)\\
 &&\hs{8ex}=0.
    \eneqn
\QED

\section{Complements to \cite{refl} } \label{sec:appendixC}


\subsection{Non-commutative rings}
Let $B$ be a graded $\cor$-algebra satisfying the following conditions:
\eq
&&\left\{\parbox{70ex}{
    there exist a commutative graded $\cor$-algebra $A$ and a
    graded $\cor$-algebra homomorphism $A\to B$ such that
    \bna
  \item $A$ satisfies (3.2) in \cite{refl},
    \item $B$ is finitely generated as a left $A$-module.
    \ee
  }
\right.
\label{eq:condB}
\eneq
Hence we have
$B_n=0$ for $n\ll0$ and $\dim B_n<\infty$ for any $n\in\Z$.
Moreover there exists $d\in\Z$ such that $BB_{\ge n}B\subset B_{\ge n+d}$
for any $n\in\Z$.  

Let $\shc$ be an abelian $\cor$-linear graded category which satisfies
(3.1) in \cite{refl}.

\Def Let $B$ be a $\cor$-algebra which satisfies \eqref{eq:condB}.
Let $\Proc(B,\shc)$ be the full subcategory of
 $\Modg(B,\Pro(\shc))$ consisting of $\Ma\in \Modg(B,\Pro(\shc))$ such that
\bna
\item $\Ma/B_{>0}\Ma \in\shc$
\item $\Ma\isoto\proolim[n]\Ma/B_{\ge n}\Ma$.
\ee
\edf
\Lemma
If $\Ma\in\Proc(B,\shc)$, then $\Ma/B_{\ge n}\Ma\in\shc$ for any $n\in\Z$.
\enlemma
\Proof
We have an epimorphism $B_n\tens(\Ma/B_{>0}\Ma)\epito B_{\ge n}\Ma/B_{>n}\Ma$.
\QED

\Lemma Let $\Ma\in \Modg(B,\Pro(\shc))$.
Then $\Ma\in\Proc(B,\shc)$ if and only if $\Ma\in\Proc(A,\shc)$.
\enlemma
\Proof
Let us take $u_j\in B_{d_j}$ ($j=1,\ldots,r$) such that
$ B=\sum_jA u_j$.
Then we have
\eqn
A_{\ge n}\Ma\subset B_{\ge n}\Ma\subset \sum_jA_{\ge n-d_j} u_j\Ma
\subset A_{\ge n-d}\Ma,
\eneqn
where $d=\max\st{d_j\mid j=1,\ldots,r}$.
It implies the lemma.
\QED

Hence, \cite[Proposition 3.3]{refl}, as well as other
results, still holds for a non-commutative $\cor$-algebra case after suitable modifications.
\Prop
The full subcategory $\Proc(B,\shc)$ of $\Modg(B,\Pro(\shc))$ has the following properties.
\bnum
\item The category $\Proc(B,\shc)$ is stable by taking subquotients 
and taking extensions
in $\Modg(B,\Pro(\shc))$. 
\label{item:affstab} 
\item Any object of $\Proc(B,\shc)$ is noetherian.\label{it:Noethe}
\item For any $\Ma\in\Proc(B,\shc)$, we have 
\bna
\label{item: Mmn}
\item for any $m\in\Z_{\ge0}$ and any $\Na\in\Proc(B,\shc)$ 
such that $\Na\subset \Ma$, there exists $n>0$ such that
$\Na\cap  B_{\ge n}\Ma\subset B_{\ge m}\Na$,\label{it:subK}
\item if $\Ma=B_{>0}\Ma$, then $\Ma\simeq 0$. \label{it:Nakayama}
\ee\label{item:mic}
\ee
\enprop

Note that if $\Ma=B_{>0}\Ma$,
then $\Ma=(B_{>0})^n\Ma\subset B_{\ge n}\Ma$ for any $n\in\Z_{>0}$.

\subsection{Duality}\label{subsec:dual}

Let $\shc$ and $\shc'$ be $\cor$-linear abelian categories which satisfy 
(3.1) in \cite{refl}.
Let $\Dual\cl \shc^\opp\to\shc'$ be  an equivalence of categories
such that $\Dual\circ q\simeq q^{-1}\circ\Dual$.
 
 Let $B$ be a graded $\cor$-algebra which satisfies \eqref{eq:condB}.

\Def\label{def:duala}
For $\Ma\in\Proc(B,\shc)$, set $\DA[B](\Ma)\seteq\proolim[n] \Dual\bl(B/BB_{\ge n}B)^*\tens[B]\Ma\br
\in\Modg\bl B^\opp,\Pro(\shc')\br$.
Similarly, for $\Na\in \Proc(B^\opp,\shc')$, set
$\DmA[B](\Na)=\proolim[n] \Dual^{-1}\bl \Na\tens[B](B/BB_{\ge n}B)^*\br
\in\Modg\bl B,\Pro(\shc)\br$.
\edf
The following proposition can be proved in a similar way as \cite[Proposition 4.5]{refl}.
\Prop\label{prop:dualflat}
Assume that $\Ma\in\Proc(B,\shc)$ is $B$-flat. Then we have:
\bnum
\item $\DA[B](\Ma)$ belongs to $\Proc(B^\opp,\shc')$,
\item $\DA[B](\Ma)$ is $B^\opp$-flat,
\item
for any finite-$\cor$-dimensional graded $B$-module $X$, we have
$$\DA[B](\Ma)\tens[B]X\simeq \Dual(X^*\tens[B]\Ma),$$
\label{item:MXY}
\item $\DmA[B]\DA[B](\Ma)\simeq \Ma$.
\ee
\enprop
\Lemma\label{lem:pbB}
Let $B\to B'$ be a homomorphism of graded $\cor$-algebras satisfying \eqref{eq:condB}.
The for a $B$-flat object $\Ma\in\Proc(B,\shc)$, we have
$$\DA[B](\Ma)\tens _BB'\simeq\DA[B'](B'\tens_B\Ma)\qt{in $\Proc(B'\;{}^\opp,\shc')$.}$$
\enlemma
\Proof
We have
\eqn
\DA[B](\Ma)\tens_BB'&&\simeq
\proolim[n]\DA[B](\Ma)\tens_B(B'/B'B'_{\ge n}B')\\
&&\underset{*}{\simeq}\proolim[n]\Dual\Bigl((B'/B'B'_{\ge n}B')^*\tens_B\Ma\Bigr)\\
&&\simeq\proolim[n]\Dual\bl(B'/B'B'_{\ge n}B')^*\tens_{B'}B'\tens_B\Ma)\\
&&\simeq\DA[B'](B'\tens_B\Ma).
\eneqn
Here $\underset{*}{\simeq}$ follows from Proposition~\ref{prop:dualflat} (iii).
\QED

\subsection{Morita equivalence}
Let $B$ and $C$ be graded $\cor$-algebras satisfying \eqref{eq:condB},
and let $P$ be a graded $(B\tens_\cor C^\opp)$-module.
\Def
We say that $B$ and $C$ are {\em graded Morita  equivalent by the kernel $P$}
if they satisfy
\eq
&&\parbox{70ex}{
  \bna\item
  $C^\opp\to \END_B(P)$ is an isomorphism,
\item
  $P$ is a finitely generated faithfully flat $B$-module.
\ee  }\label{eq:Morita}
\eneq
\edf

Then, by standard arguments, we see that
\bnum
\item $P$ is a finitely generated faithfully flat $C^\opp$-module,
  \item $B\to \END_{C^\opp}(P)$ is an isomorphism,
\item $\HOM_B(P,B)$ is a finitely generated faithfully flat $B^\opp$-module
  and a finitely generated faithfully flat $C$-module,
\item Two functors
  \eqn &&P\tens_C\scbul: \Modg(C,\Pro(\shc))\to\Modg(B,\Pro(\shc))\qtq\\
  && \ihom_{B}(P,\scbul)\simeq\HOM_B(P,B)\tens_B\scbul: \Modg(B,\Pro(\shc))\to\Modg(C,\Pro(\shc))
  \eneqn
  are equivalences of categories and quasi-inverses to each other.

\item\label{it:PBC}
  $\HOM_B(P,B)\simeq  \HOM_{C^\opp}(P,C)$ as $C\tens B^\opp$-modules.
\item
 $C^\opp$ and $B^\opp$ are graded  Morita  equivalent by the kernel $P$, and
    $C$ and $B$ are  graded Morita equivalent by the kernel $\HOM_B(P,B)$.

  \ee
The statement  \eqref{it:PBC}
  can be seen as follows:
  the functor $\HOM_B(P,\scbul)$ sends $B$ to $\HOM_B(P,B)$ and
  its quasi-inverse $P\tens_ C\scbul$ sends $\HOM_{C^\opp}(P,C)$ to
  $$P\tens_ C\HOM_{C^\opp}(P,C)\simeq\HOM_{C^\opp}(P,P\tens _CC)\simeq
  \HOM_{C^\opp}(P,P)\simeq B.$$

  \Lemma
  There exists $d\in \Z_{>0}$ such that
  \eqn B_{\ge n+d}P\subset P_{\ge n}\subset B_{\ge n-d}P,\\
  PC_{\ge n+d}\subset P_{\ge n}\subset PC_{\ge n-d}.
  \eneqn
  \enlemma
  \Proof
  There exists $c\in\Z$ such that
  $P=P_{\ge c}$. Hence $B_{\ge n}P\subset P_{\ge n+c}$.
  Let us take $u_j\in P_{r_j}$ ($j=1,\ldots, s$) such that
  $P=\sum_{j}Bu_j$.
  Then we have
  $$P_{\ge n}=\sum_j B_{\ge n-r_j}u_j\subset B_{\ge n-d} P$$
  for $d\ge r_j$.
  \QED
  \Prop\label{prop:MoritaD}
  Assume that $B$ and $C$ are Morita equivalent by the kernel $P$.
  Then for any $B$-flat $\Ma\in\Proc(B,\shc)$, we have
  $$\DA[B](\Ma)\tens_BP\simeq\DA[C](\HOM_B(P,B)\tens_B\Ma)
\qt{in $\Proc(C^\opp,\shc')$.}$$
\enprop
\Proof
We have
\eqn
\DA[B](\Ma)\tens_BP
&&\simeq\proolim[n]\DA[B](\Ma)\tens_B (P/BP_{\ge n}C)\\
&&\simeq\proolim[n]\Dual\bl(P/BP_{\ge n}C)^*\tens_B\Ma\br.
\eneqn
On the other hand, we have 
\eqn
\proolim[n](P/BP_{\ge n}C)^*
&&\simeq\proolim[n](P/PC_{\ge n}C)^*\\
&&\simeq\proolim[n]\HOM_\cor\bl P\tens_C(C/CC_{\ge n}C),\cor\br\\
&&\simeq\proolim[n]\HOM_{C^\opp}\bl P, \HOM_\cor(C/CC_{\ge n}C,\cor)\br\\\
&&\simeq\proolim[n](C/CC_{\ge n}C)^*\tens_C\HOM_{C^\opp}(P,C)\\
&&\simeq\proolim[n](C/CC_{\ge n}C)^*\tens_C\HOM_{B}(P,B).
\eneqn

Hence
\eqn
\DA[B](\Ma)\tens_BP
&&\simeq\proolim[n]\Dual\bl(C/CC_{\ge n}C)^*\tens_C\HOM_{B}(P,B)\tens_B
\Ma\br\\
&&\simeq\DA[C](\HOM_B(P,B)\tens_B\Ma).
\eneqn
\QED

\subsection{Monoidal category}
Let $\shc$ be a $\cor$-linear graded monoidal category satisfying
(5.1) in \cite{refl}.

Let $\shc^{\rd}$ (resp.\ $\shc^\ld$ ) be the full subcategory of $\shc$ consisting of objects with right duals (resp.\ left dual).

Then it is easy to see
\bnum
\item
  $\shc^{\rd}$ and $\shc^\ld$ are stable by taking kernels, cokernels and
  tensor products,
\item
  $\shc^{\rd}$ and $\shc^\ld$ also satisfy (5.1) in \cite{refl}.
\ee
Let $\D\cl \shc^{\rd}\to (\shc^{\ld})^{\rev,\opp}$ be the functor of taking
the right dual. Then, it is an equivalence of categories.
Then, we can apply \S\,\ref{subsec:dual} by defining
$$\DB(\Ma)\seteq\proolim[n] \D\bl(B/BB_{\ge n}B)^*\tens[B]\Ma\br.$$

\subsection{Proof of $\Psi_{+-}(\ang{i^c}_{z_j})\simeq\Da(\ang{i^c}_{z_j})$} \label{sec:appendixC5}
\hfill

{\em We ignore grading shifts in this section.}
Remark that we already proved that $\ang{i}_{z_i}\in\Proc(\cor[z_i],(\Loc_+)^\rd)$.

Let us take a graded $\cor$-algebra $A$ satisfying (3.2) in \cite{refl}.
Let $c\in\Z_{\ge0}$, and let $a(z_i)\in (A[z_i])_{2c\sfd_i}$ be a quasi-monic polynomial of
degree $c$ (in $z_i$).

Set 
\eqn
\Rh&&\seteq R(c\al_i)\qqt{the quiver Hecke algebra of degree $c\al_i$,}\\
P&&\seteq \cor[x_1,\ldots, x_c]\subset \Rh,\\
Z&&\seteq \cor[x_1,\ldots, x_c]^{\sym_c}\subset P\subset \Rh\qqt{the center of $\Rh$.}\\
P(c\al_i)&&\seteq \dfrac{\Rh}{\sum_{k=1}^{c-1}\Rh\tau_k}\qqt{the indecomposable projective $\Rh$-module,}\\
P(c\al_i)^\opp&&\seteq \dfrac{\Rh}{\sum_{k=1}^{c-1}\tau_k\Rh}\qqt{the indecomposable projective $\Rh^\opp$-module,}\\
\Rc&&\seteq\dfrac{A\tens_\cor\Rh}{(A\tens_\cor\Rh) a(x_c)(\At\Rh)}\qqt{the cyclotomic quiver Hecke algebra (see \cite{KK11}),}
\eneqn 
Hence the algebras $\Rh$, $P$, $Z$, $\At\Rh$ and $\Rc$  satisfy the condition \eqref{eq:condB}.

Set
$$\ang{i^c}_a\seteq \Rc\tens_\Rh P(c\al_i)
\simeq\dfrac{\At\Rh}{\sum_{1\le k<c}(\At\Rh)\tau_k+(\At\Rh) a(x_c)}\in \Proc(A,R\gmod).$$

Then we have
$$(A/A_{>0})\tens_A\ang{i^c}_a\simeq \ang{i^c}.$$
The following statements are well-known 
\eq\label{eq:Moritas}
&&\hs{2ex}\left\{\parbox{75ex}{
  \bnum
\item $\Rh$ and $Z$ are graded Morita equivalent by the kernel $P(c\al_i)$
(\cite{Rouquier08}).
  \item
  $\Rc$ and $A$ are graded Morita equivalent by the kernel $\Rc\tens_\Rh P(c\al_i)=\ang{ i^c}_a$
  (\cite{KK11}, cf.\  \cite{Rouquier08}).
\item $\Rh\simeq\HOM_{P}(\Rh,P)$ as an $(\Rh,P)$-bimodule (\cite{Rouquier08}).
  \item$\Rc\tens_\Rh P(c\al_i)\simeq\HOM_{(\Rc)^\opp}(P(c\al_i)^\opp\tens_\Rh\Rc,\Rc)$ as an $\Rc$-module.
  \ee
}\right.
\eneq

Note that $\ang{i^{-\ang{h_i,\al_j}}}_{z_j}$ is equal to $\ang{i^c}_a$ with
$c=-\ang{h_i,\al_j}$, $A=\cor[z_j]$ and $a(z_i)=Q_{i,j}(z_i,z_j)$.

  \Prop We have
$$  \Psi_{+-}(\ang{i^c}_a)
\simeq\DB[A]\bl\Qt_+(\ang{i^c}_a)\br
\qt{in $\Proc(A,\Loc_+)$.}$$
  \enprop
  \Proof In the course of the proof, we drop $\Qt_+$ to simplify the notations. 
    Note that $R(c\al_i)\in\Proc\bl(\Rh)^\opp,\Loc_-)$.

  We have
  $\Psi_{+-}(R(\al_i))\simeq\Da(\ang{i}_{z_i})$.

  Hence, we have
  \eqn\Psi_{+-}(R(c\al_i))
  &&\simeq\Psi_{+-}\bl\ang{i}_{x_1}\conv\cdots\conv\ang{i}_{x_c}\br\\
  &&\simeq\Psi_{+-}(\ang{i}_{x_1})\conv\cdots\conv\Psi_{+-}(\ang{i}_{x_c})\\
      &&\simeq\Da(\ang{i}_{x_1})\conv\cdots\conv\Da(\ang{i}_{x_c}).
      \eneqn
      We have

      \eqn
    && \Da(\ang{i}_{x_1})\conv\cdots\conv\Da(\ang{i}_{x_c})\\
    &&\hs{3ex}\simeq  \proolim[n]
    \D\bl\ang{i}_{x_1}\tens_{\cor[x_1]}(\cor[x_1]/\cor[x_1]x_1^n)^*\br
     \conv\cdots\conv \D\bl\ang{i}_{x_c}\tens_{\cor[x_c]}
     (\cor[x_c]/\cor[x_c]x_c^n)^*\br\\
     &&\hs{3ex} \simeq  \proolim[n]    \D\Bigl(
      \bl\ang{i}_{x_c}\tens_{\cor[x_c]}
      (\cor[x_c]/\cor[x_c]x_c^n)^*\br \conv\cdots\conv
      \bl\ang{i}_{x_1}\tens_{\cor[x_1]}
     (\cor[x_1]/\cor[x_1]x_1^n)^*\br\Bigr)\\
      &&\hs{3ex}\simeq 
        \proolim[n]    \D\bl
        \Rh\tens_P(P/P_{\ge n})^*\br.
        \eneqn
        On the other hand, we have
        \eqn
        \Rh\tens_P(P/PZ_{\ge n})^*
        &&\hs{3ex}\underset{*}{\simeq}\HOM_P(\Rh,P)\tens_P(P/PZ_{\ge n})^*\\
        &&\hs{3ex}\simeq\HOM_P\bl \Rh,P\tens_P(P/PZ_{\ge n})^*\br\\
                &&\hs{3ex}\simeq\HOM_P\bl \Rh,(P/PZ_{\ge n})^*\br\\
                &&\hs{3ex}\simeq\HOM_{\cor}\bl \Rh\tens_P (P/PZ_{\ge n}),\cor\br\\
                &&\hs{3ex}\simeq(\Rh/ Z_{\ge n}\Rh)^*.
                \eneqn
Here $\underset{*}{\simeq}$ follows from \eqref{eq:Moritas}.
                
                Hence we obtain

        \eqn&&\Psi_{+-}(\Rh )
        \simeq\proolim[n]\D\bl(\Rh/Z_{\ge n}\Rh)^*\br
        \simeq\proolim[n]\D\bl \Rh\tens_\Rh(\Rh/Z_{\ge n}\Rh)^*\br\simeq\DB[\Rh^\opp](\Rh).
        \eneqn
        Note that
        $\Psi_{+-}( \Rh )\in\Proc(\Rh^\opp,\Loc_+)$
        and $\DB[\Rh^\opp](\Rh)\in\Proc(\Rh,\Loc_+)$.
        Here we identify $\Modg(\Rh,\Pro(\Loc_+))$ with $\Modg(\Rh^\opp,\Pro(\Loc_+))$ by the isomorphism $\xi\cl\Rh\isoto\Rh ^\opp$ given by
        $x_k\mapsto x_{c+1-k}$, $\tau_k\mapsto -\tau_{c-k}$.
Since $\xi_*(\Rc)\simeq\Rc$, we have
                \eqn
                \Psi_{+-}(\Rc)&&\simeq \Psi_{+-}(\Rh\tens_\Rh\Rc)\\
                &&\simeq                \Psi_{+-}(\Rh)\tens_\Rh\Rc\\
                &&\simeq\DB[\Rh^\opp](\Rh)  \tens_{\Rh^\opp}
                \xi_*(\Rc)\\
                &&\underset{*}{\simeq}\DB[(\Rc)^\opp](\Rh\tens_\Rh\Rc)\\
                &&\simeq\DB[(\Rc)^\opp](\Rc).
                \eneqn
                Here $\underset{*}{\simeq}$ follows from Lemma~\ref{lem:pbB},
                Hence we have
                \eqn
                \Psi_{+-}(\ang{i^c}_a)
&&\simeq\Psi_{+-}\bl(\Rc)\tens_\Rc\ang{i^c}_a\br
\simeq \Psi_{+-}(\Rc)\tens_\Rc \ang{i^c}_a\\
&&\simeq\DB[(\Rc)^\opp]\bl(\Rc)\tens_{(\Rc)^\opp}\xi_*(\ang{i^c}_a ) \br\\
&&\underset{(1)}{\simeq}\DB[A]\Bigl(
\HOM_{(\Rc)^\opp}\bl\xi_*(\ang{i^c}_a),(\Rc)^\opp\br
\tens_{(\Rc)^\opp}\Rc\Bigr)\\
&&\underset{(2)}{\simeq}\DB[A]\bl\Rc
\tens_\Rc\HOM_{(\Rc)^\opp}(P(c\al_i)^\opp\tens_\Rh\Rc,\Rc)\br\\
                       &&\underset{(3)}{\simeq}\DB[A](\ang{i^c}_a).
                       \eneqn
                       Here $\underset{(1)}{\simeq}$ follows from Proposition~\ref{prop:MoritaD}
                       together with \eqref{eq:Moritas}, 
$\underset{(2)}{\simeq}$ follows from
$\xi_*(\ang{i^c}_a)\simeq P(c\al_i)^\opp\tens_{\Rh}\Rc$,
and $\underset{(3)}{\simeq}$ follows from \eqref{eq:Moritas}.
 \QED




\begin{thebibliography}{99}




\bibitem{BK09}
J.~Brundan and A.~Kleshchev, \emph{Blocks of cyclotomic Hecke algebras and Khovanov-Lauda algebras}, Invent. Math. \textbf{178} (2009), 451--484.


\bibitem{KK11}
S.-J. Kang and M. Kashiwara, \emph{Categorification of Highest Weight Modules via Khovanov-Lauda-Rouquier Algebras},
Invent. Math. \textbf{190} (2012), no. 3, 699--742.
%
\bibitem{K^3}
S.-J. Kang, M. Kashiwara and M. Kim, {\em Symmetric quiver
	Hecke algebras and R-matrices of quantum affine algebras},
Invent. math. {\bf 211} (2018), 591--685,
arXiv:1304.0323 v3.
%
\bibitem{KKKO15}
S.-J. Kang, M. Kashiwara,  M. Kim  and   S.-j. Oh,
\newblock{\em Simplicity of heads and socles of tensor products},
Compos. Math. \textbf{151} (2015), no. 2, 377--396.
%
\bibitem{KKKO18}
\bysame,
\newblock{\em Monoidal categorification of cluster algebras}, 
J. Amer. Math. Soc. \textbf{31} (2018), no. 2, 349--426.


  \bibitem{KKOP18} 
M.~Kashiwara, M. Kim, S.-j. Oh, and  E.~Park,
\newblock{\em Monoidal categories associated with strata of flag manifolds},
Adv. Math. \textbf{328} (2018), 959--1009.

\bibitem{loc1}
 \bysame,  
\newblock{\em Localizations for quiver Hecke algebras}, Pure Appl. Math. Q. \textbf{17} (2021), no. 4, 1465--1548.

\bibitem{loc2}
\bysame, 
 \newblock{\em Localizations for quiver Hecke algebras II},
Proc. London Math. Soc. (4) \textbf{127} (2023) 1134--1184.

\bibitem{loc3}\bysame,
 \newblock{\em Localizations for quiver Hecke algebras III}, 
Math. Ann. 390 (2024), no. 4, 5075--5108. 


\bibitem{refl} \bysame,
  \newblock{\em Affinizations, R-matrices and reflection functors},
  Adv. Math. \textbf{443} (2024), Paper No. 109598, 83 pp. 


\bibitem{ext} 
  \bysame,
  \newblock{\em Braid symmetries on bosonic extensions},
  arXiv:2408.07312v1.

  \bibitem{KS} M.~Kashiwara and P.~Schapira,
\emph{Categories and Sheaves},
Grundlehren der Mathematischen Wissenschaften, vol.~\textbf{332}, Springer-Verlag, Berlin (2006).
%
%
%
%
%
%
%
%

%
%
%
%
%




%
%
%
%
%
%
\bibitem{KP18}
M.~Kashiwara and E.~Park, \newblock{\em Affinizations and $R$-matrices for quiver Hecke algebras},
J. Eur. Math. Soc. \textbf{20}, (2018), 1161--1193.
%
%
%
\bibitem{Kato14}S.~Kato,
{\em Poincar\'e-Birkhoff-Witt bases and Khovanov-Lauda-Rouquier algebras},
Duke Math. J. {\bf163} (2014), no. 3, 619--663. 
%
\bibitem{Kato20}\bysame,
{\em On the monoidality of Saito reflection functors},
Int. Math. Res. Not. 2020, no. 22, 8600--8623. 
%
%
\bibitem{KL09}
M.~Khovanov and A. Lauda, \emph{A diagrammatic approach to
categorification of quantum groups
 {I}}, Represent. Theory \textbf{13} (2009), 309--347.
%
\bibitem{KL11}
\bysame, \emph{A diagrammatic approach to categorification of
  quantum groups {II}}, Trans. Amer. Math. Soc. \textbf{363} (2011),
  2685--2700.
%
\bibitem{LV11}
A.~Lauda and M.~Vazirani, \emph{Crystals from categorified quantum groups},  Adv. Math. \textbf{228} (2011), no.~2, 803--861.
%
  \bibitem{Lu93}
  G.~Lusztig, {\em Introduction to Quantum Groups}, Progr. Math. \textbf{110}  Birkh\"{a}user, 1993.
%
%
%
%

 \bibitem{Murata}
H. Murata, \emph{Affine highest weight structures on module categories over quiver Hecke algebras},  arxiv:2412.12903v3.

\bibitem{Rouquier08}
R. Rouquier,  \emph{2-Kac–Moody algebras},  arXiv:0812.5023v1.

\bibitem{R11}
\bysame, {\em Quiver Hecke algebras and 2-Lie algebras},
Algebra Colloq. {\bf 19} (2012), no. 2, 359--410.

%
%
\bibitem{TW16}
P. Tingley and B. Webster,
{\em Mirkovi\'c-Vilonen polytopes and Khovanov-Lauda-Rouquier algebras},
Compos. Math. {\bf 152} (2016), no. 8, 1648--1696.
%
%
\bibitem{VV09}
M. Varagnolo and E. Vasserot,
\emph{Canonical bases and KLR algebras},
J. Reine Angew. Math. \textbf{659} (2011), 67--100.
%
%
%
%
%
%
%
%
%
%
%
%
%
%
%
\end{thebibliography}
\end{document}